\documentclass[10pt,a4paper]{article}
\usepackage{fullpage,mathrsfs,graphicx,framed,color,amssymb,amsmath,amsthm}
\usepackage[boxed, noline, ruled, linesnumbered]{algorithm2e}
\usepackage{epsfig, graphicx}
\usepackage{latexsym,amsfonts,amsbsy,amssymb}
\usepackage{amsmath,amsthm}
\usepackage{enumerate}
\usepackage{url}
\usepackage{changes} % for revision

\makeatletter %'@' is now a normal "letter" for TeX

\@addtoreset{equation}{section}
\makeatother %'@' is restored as a "non-letter" character for TeX
\allowdisplaybreaks % For long formula
\newtheorem{theorem}{Theorem}[section]
\newtheorem{lemma}{Lemma}[section]
\newtheorem{corollary}{Corollary}[section]
\newtheorem{remark}{Remark}[section]
\newcommand{\comm}[1]{{\color{red}#1}}
\newcommand{\revise}[1]{{\color{blue}#1}}
\begin{document}
\title{Augmented Subspace Scheme for Eigenvalue Problem by Weak Galerkin Finite Element
Method}
\author{
Yue Feng\footnote{LSEC,
Academy of Mathematics and Systems Science,
Chinese Academy of
Sciences, No.55, Zhongguancun Donglu, Beijing 100190, China, and School of
Mathematical Sciences, University of Chinese Academy
of Sciences, Beijing 100049, China
({\tt fengyue@amss.ac.cn})},\ \ \ Zhijin Guan\footnote{LSEC,  Academy of Mathematics and Systems Science,
Chinese Academy of
Sciences, No.55, Zhongguancun Donglu, Beijing 100190, China, and School of
Mathematical Sciences, University of Chinese Academy
of Sciences, Beijing 100049, China ({\tt guanzhijin@lsec.cc.ac.cn})},\ \ \
 Hehu Xie\footnote{LSEC,
Academy of Mathematics and Systems Science,
Chinese Academy of
Sciences, No.55, Zhongguancun Donglu, Beijing 100190, China, and School of
Mathematical Sciences, University of Chinese Academy
of Sciences, Beijing 100049, China ({\tt hhxie@lsec.cc.ac.cn})} \  \  and\  \ 
Chenguang Zhou\footnote{Faculty of Science, Beijing University of Technology, 
Beijing 100124, China (Corresponding author: {\tt zhoucg@bjut.edu.cn})} }

\date{}
\maketitle
\begin{abstract}
This study proposes a class of augmented subspace schemes for the weak Galerkin (WG) 
finite element method used to solve eigenvalue problems.
The augmented subspace is built with the conforming linear finite element space
defined on the coarse mesh and the eigenfunction approximations in the WG finite element space 
defined on the fine mesh. Based on this augmented subspace, solving the eigenvalue 
problem in the fine WG finite element space can be reduced to the solution of 
the linear boundary value problem in the same WG finite element space and 
a low dimensional eigenvalue problem in the augmented subspace.
The proposed augmented subspace techniques have the second order convergence 
rate with respect to the coarse mesh size, as demonstrated by the accompanying 
error estimates. Finally, a few numerical examples are provided to validate 
the proposed numerical techniques.


\end{abstract}

\section{Introduction}
One of the most important tasks in contemporary scientific and engineering society 
is solving eigenvalue problems. The difficulty of solving eigenvalue problems 
is invariably higher than that of solving similar linear boundary value problems 
due to the increased computing and memory requirements.
Large-scale eigenvalue problem solving in particular will provide formidable obstacles 
to scientific computing. Numerous eigensolvers have been developed so far, including 
the Jacobi-Davidson type technique \cite{Bai}, the Preconditioned INVerse ITeration (PINVIT) method \cite{BramblePasciakKnyazev,PINVIT,Knyazev}, 
the Krylov subspace type method (Implicitly Restarted 
Lanczos/Arnoldi Method (IRLM/IRAM) \cite{Sorensen}), 
and the Generalized Conjugate Gradient Eigensolver (GCGE) 
\cite{LiWangXie, LiXieXuYouZhang, ZhangLiXieXuYou}.
The orthogonalization processes involved in solving Rayleigh-Ritz problems 
are a common bottleneck in the design of effective parallel techniques for 
identifying a large number of eigenpairs, and they are included in all of 
these widely used approaches.

A class of augmented subspace methods and their multilevel correction methods 
has been proposed recently in \cite{full,HongXieXu,LinXie_MultiLevel,Xie_IMA,Xie_JCP,Xie_BIT,XieZhangOwhadi,XuXieZhang} 
for the solution of eigenvalue problems.
This kind of technique creates an augmented subspace using the low dimensional 
finite element space generated on the coarse grid, which is employed in each correction step.
The notion of an augmented subspace gives rise to a class of augmented subspace techniques that need just the final finite element space on the finest mesh and the low dimension finite element space on the coarse mesh.
Using the augmented subspace methods, the solution of the eigenvalue problem on the final level of mesh can be transformed to the
solution of linear boundary value problems on the final level of mesh and the solution
of the eigenvalue problem on the low dimensional augmented subspace.
Even the coarse and finest meshes lack nested properties, these kinds of algorithms can still 
work  \cite{Dang}.
The multilevel correction methods, which are based on the augmented subspace methods, 
provide ways to construct multigrid methods for eigenvalue problems \cite{full,HongXieXu,Xie_IMA,Xie_JCP,XieZhangOwhadi}. 
In addition, the authors design an eigenpair-wise parallel eigensolver 
for the eigenvalue problems in \cite{XuXieZhang}. 
A significant amount of the wall time in the parallel computation 
is saved by using this kind of parallel approach, 
which avoids performing orthogonalization and inner-products in the high dimensional space.
However, the aforementioned references are mostly investigated using 
conforming finite element methods. There are few results on the augmented 
subspace approaches based on nonstandard finite element methods for solving eigenvalue problems.

The WG method, which was initially introduced and explored in \cite{WangYe}, 
concerns the finite element methods utilized to solve partial differential 
equations in which the differential operators, such as gradient operator, 
divergence operator, curl operator, and so on, are approximated as distributions by weak forms.
The WG approach employs generalized discrete weak derivatives and parameter-free 
stabilizers to weakly enforce continuity in the approximation space, 
in contrast to the standard finite element technique.
Consequently, it ought to be more convenient to create high order 
precision discretization than the conforming finite element approach.
Additionally, the WG approach can be easily implemented on polygonal meshes 
thanks to the relaxation of the continuity constraint, which also 
gives additional freedom for $h$- and $p$-adaptation.
So far, the WG method has been applied to various partial differential equations,
such as the parabolic equation \cite{LiWang_WG,ZhouGaoLiSun}, the biharmonic
equation \cite{MuWangYe,WangWang,ZhangZhai}, the Brinkman equation 
\cite{MuWangYe_WG,ZhaiZhangMu}, 
the Helmholtz equation \cite{MuWangYe_H,MuWangYeZhao} and the Maxwell 
equation \cite{MuWangYeZhang_WG}.
The convergence analysis and several lower bound findings are produced 
in \cite{ZhaiXieZhangZhang_WG}, 
where the WG approximation to the eigenvalue problems is firstly studied.
Then, using the WG approach, the authors create a kind of two-grid or two-level schemes \cite{ZhaiXieZhangZhang_TwoGrid}, and in \cite{ZhaiHuZhang}, the shifted-inverse 
power technique is taken into consideration under the two-grid schemes.
Based on the theoretical analysis presented in \cite{ZhaiXieZhangZhang_TwoGrid}, 
it can be inferred that there is no independent relationship between 
the coarse and fine mesh sizes. As a result, the approaches cannot 
be used to develop an eigensolver for algebraic eigenvalue problems 
resulting from differential operator eigenvalue problems discretized by WG.

This paper's contribution is the augmented subspace methods for 
eigenvalue problem that are based on the WG approximation. 
To the best of our knowledge, this is the first work aimed at the numerical 
analysis of the WG finite element discretization-based augmented subspace approaches 
for eigenvalue problems.
In contrast to the findings in \cite{ZhaiXieZhangZhang_TwoGrid}, our approaches' 
selections for the coarse and fine mesh sizes are independent of one another. 
The algebraic eigenvalue problems that result from the WG approximation to 
the differential eigenvalue problems can then be solved by designing 
an eigensolver using the proposed techniques.
Furthermore, we demonstrate the algebraic error estimate for the WG augmented 
subspace approaches that follows 
\begin{eqnarray*}
\left\|\bar u_h-u_h^{(\ell+1)}\right\|_{a,h} \leq CH^2 
\left\|\bar u_h-u_h^{(\ell)}\right\|_{a,h},
\end{eqnarray*}
when the computing domain is convex.

This paper is organized as follows. We provide the WG approaches for the eigenvalue 
problems and deduce the associated error estimates in Section \ref{Section_2}. 
These results give explicit dependence of the error estimates on the eigenvalue 
distribution which is another contribution of this paper. The majority of this work, 
Section \ref{Section_3},  contains the augmented subspace techniques and 
the associated error estimates. A few numerical examples are given in 
Section \ref{Section_4} to validate the suggested augmented subspace 
algorithms' convergence rates. Lastly, the final section has a few closing thoughts.

\section{Discretization by WG finite element method}\label{Section_2}
The WG finite element approach for the second order elliptic eigenvalue problem 
is presented in this section. Additionally, the associated error estimates are offered. The letter $C$, with or without subscripts, symbolizes a generic positive constant for this purpose that may vary at various places in this work.

Here, we consider the numerical method to solve the following second order elliptic eigenvalue problem:
Find $(\lambda, u) \in \mathbb R\times H_0^1(\Omega)$ such that
\begin{eqnarray}\label{eigenvalue_problem}
\left\{
\begin{array}{rcl}
-\nabla\cdot(A\nabla u) &=&\lambda u,\ \ \ {\rm in}\ \Omega,\\
u &=&0,\ \ \ \ \ {\rm on}\ \partial\Omega,\\
(A\nabla u, \nabla u) &=&1,
\end{array}
\right.
\end{eqnarray}
where $\Omega$ denotes a convex bounded polygonal or polyhedral domain in $\mathbb R^d$, $d=2,3$,
and $A \in\left[L^{\infty}(\Omega)\right]^{d \times d}$ is a symmetric 
matrix-valued function on $\Omega$ with
suitable regularity.
Assume that  there exist positive constants $c$ and $C$ such that the matrix $A$ satisfies
the following property
\begin{eqnarray}
c \xi^T \xi \leq \xi^T A(x) \xi \leq C \xi^T \xi \ \ \text {       for all } 
\xi \in \mathbb R^d \ \text{and}\ \ x \in \Omega.
\end{eqnarray}

In order to define the WG finite element method for the eigenvalue problem, (\ref{eigenvalue_problem})
should  be written as the following variational form:
Find $(\lambda, u )\in \mathbb R\times V$ such that $a(u,u)=1$ and
\begin{eqnarray}\label{weak_eigenvalue_problem}
a(u,v)=\lambda b(u,v),\quad \forall v\in V,
\end{eqnarray}
where $V:=H_0^1(\Omega)$ \cite{Adams} and
\begin{eqnarray}
a(u,v) = (A\nabla u, \nabla v),\ \ \ \ b(u,v) = (u,v).
\end{eqnarray}
Furthermore,  based on the bilinear forms $a(\cdot,\cdot)$ and $b(\cdot,\cdot)$, we can define the
norms on the space $V$ as follows
\begin{eqnarray}\label{Nomr_a}
\left\|v\right\|_{a} = \sqrt{a(v,v)},\ \ \ \forall v\in V, \ \ \  \left\|w\right\|_{b} = \sqrt{b(w,w)},
\ \ \ \forall w\in L^2(\Omega).
\end{eqnarray}

It is well known that the eigenvalue problem (\ref{weak_eigenvalue_problem})
has an eigenvalue sequence $\{\lambda_j \}$ (cf. \cite{BabuskaOsborn_1989,Chatelin}),
$$0<\lambda_1\leq \lambda_2\leq\cdots\leq \lambda_k\leq \cdots,\ \ \ 
\lim_{k\rightarrow\infty}\lambda_k=\infty.$$
And the associated eigenfunctions are provided as
$$u_1, u_2, \cdots, u_k, \cdots.$$
Here $a(u_i,u_j)=\delta_{ij}$ ($\delta_{ij}$ denotes the Kronecker function).

%-----------------------------------------------------------------------------------------
Now, let us define the WG finite element space for the  eigenvalue problem
(\ref{weak_eigenvalue_problem}). First we generate a shape-regular, quasi-uniform mesh $\mathcal{T}_h$
of the computing domain $\Omega\subset \mathbb R^d\ (d=2,3)$.  Denote by $\mathcal E_h$
the set of all edges or faces of the mesh $\mathcal T_h$.
For simplicity, in this paper, we only consider
the triangle or tetrahedral  mesh.  The diameter of a cell $K\in\mathcal{T}_h$
is denoted by $h_K$ and the mesh size $h$ describes the maximal diameter of all cells
$K\in\mathcal{T}_h$.
For each cell $K \in \mathcal T_h$, we use $K_0$ and $\partial K$ to denote the interior and the
boundary of $K$.  In the sense of geometry, $K_0$ is identical to $K$. Then we identify them if no ambiguity.
Based on the mesh $\mathcal{T}_h$, we can construct the WG finite element space denoted by
 $V_h$  as follows
\begin{equation}\label{linear_fe_space}
  V_h = \Big\{ v: v|_{K_0} \in \mathcal P_r(K_0)\ {\rm for}\ K\in \mathcal T_h;
  v|_e \in \mathcal P_s(e) \ {\rm for}\ e\in\mathcal E_h, \ {\rm and}\ v|_e = 0\ {\rm for}\
  e\in \mathcal E_h\cap\partial\Omega\Big\},
\end{equation}
where  $\mathcal{P}_r(K_0)$ denotes the set of polynomials of degree no more than 
the integer $r\geq 0$, $\mathcal P_s(e)$ is the set of polynomials of 
degree no more than the integer $s\geq 0$.  In this paper, we are concerned with the
cases of $s = r$ or $r+1$.
From the definition of $V_h$, it is easy to know that the function in $V_h$ does not require any continuity
across interior edges/faces.  Actually, the function in $V_h$ can be characterized by its value on the interior of each element
and its value on edges/faces. Therefore, the functions in $V_h$ can be represented with two components, $v=\left\{v_0, v_b\right\}$,
where $v_0$ denotes the value of $v$ on all $K_0$ and $v_b$ denotes the value of $v$ on $\mathcal{E}_h$.
The polynomial space $\mathcal P_s(e)$ consists of two choices: $s=r$ or $r+1$ and the corresponding weak function
space will sometimes be abbreviated as $V_{r, r}$ or $V_{r, r+1}$, respectively.

In order to define the WG method for the eigenvalue problem (\ref{weak_eigenvalue_problem}),
we introduce the discrete weak gradient operator, which is defined on each element $K \in \mathcal{T}_h$.
For the choices of $V_h$ given above, i.e., using $V_{r, r}$ or $V_{r, r+1}$, suitable definitions of the weak
gradient involve the Raviart-Thomas (RT) element  or the Brezzi-Douglas-Marini (BDM) element \cite{FortinBrezzi},  respectively.
Let $K$ be either a triangle or a tetrahedron and denote by $\widehat{\mathcal P}_t(K)$ the set of homogeneous
polynomials of order $t$ in the variable $\mathbf{x}=\left(x_1, \ldots, x_d\right)^T$.
Define the BDM element by $G_r(K)=\left[\mathcal P_{r+1}(K)\right]^d$ and the RT element by
$G_r(K)=\left[\mathcal P_r(K)\right]^d+\widehat{\mathcal P}_r(K) \mathbf{x}$ for $r \geq 0$.
Then, we can define a discrete space
$$
\mathbf\Sigma_h=\left\{\mathbf{q} \in\left(L^2(\Omega)\right)^d:\left.\mathbf{q}\right|_K \in
G_r(K) \text { for } K \in \mathcal{T}_h\right\} .
$$
In the definitions of $V_h$ and $\mathbf\Sigma_h$, the RT element is coupled with $V_{r, r}$
while the BDM element is coupled with $V_{r, r+1}$. We should point out
that $\mathbf\Sigma_h$ is not necessarily a subspace of $H(\operatorname{div}, \Omega)$,
since it does not require any continuity in the normal direction across edges/faces.

The discrete weak gradient of $v_h\in V_h$ denoted by $\nabla_w v_h$ is defined as
the unique polynomial $\left.\left(\nabla_w v_h\right)\right|_K \in G_r(K)$ satisfying the following equation
\begin{eqnarray}\label{Definition_Weak_Gradient}
\left(\nabla_w v_h, \mathbf{q}\right)_K=-\left(v_0, \nabla \cdot \mathbf{q}\right)_K
+\left\langle v_b, \mathbf{q} \cdot \mathbf{n}\right\rangle_{\partial K} \quad \text { for all } \mathbf{q} \in G_r(K),
\end{eqnarray}
where $\mathbf{n}$ is the unit outward normal vector on $\partial K$.
Clearly, such a discrete weak gradient is always well-defined. Furthermore, if $v \in H^1(K)$, 
i.e., $v_b=\left.v_0\right|_{\partial K}$, and $\nabla v \in G_r(K)$. Then we have $\nabla_w v=\nabla v$.
Here we only consider the $V_{r, r}$-RT and $V_{r, r+1}$-BDM pairs on simplicial elements. % in the discretization.
Of course, there are many other different choices of discrete spaces in the WG method,
defined on either simplicial meshes or general polytopal meshes \cite{MuWangYe,WangWang}.

In order to define an interpolation operator for the WG finite element space,
we define an $L^2$ projection from $V$ onto $V_h$ by setting $Q_h v \equiv\left\{Q_0 v, Q_b v\right\}$,
where $\left.Q_0 v\right|_{K_0}$ is the local $L^2$ projection of $v$ to $\mathcal P_r\left(K_0\right)$,
for $K \in \mathcal{T}_h$, and $\left.Q_b v\right|_e$ is the local $L^2$ projection to $\mathcal P_s(e)$,
for $e \in \mathcal{E}_h$. We also introduce $\mathbb{Q}_h$ the $L^2$ projection onto $\mathbf\Sigma_h$.
It is well known that the following operator identity holds \cite{WangYe}:
\begin{eqnarray}\label{Property_Weak_Gradient}
\mathbb{Q}_h \nabla v=\nabla_w Q_h v, \quad \text { for all } v \in V.
\end{eqnarray}
For the $V_{r, r}$-RT and $V_{r, r+1}$-BDM pairs, the identity  (\ref{Property_Weak_Gradient}) 
shows that the discrete weak gradient is a good approximation to the classical gradient  \cite{WangYe}.
%as summarized in the following lemma

Then, the WG finite element method for the eigenvalue
problem (\ref{weak_eigenvalue_problem}) can be defined as follows:
Find $(\bar{\lambda}_h, \bar{u}_h)\in \mathbb R \times V_h$
such that $a_h(\bar{u}_h,\bar{u}_h)=1$ and
\begin{eqnarray}\label{Weak_Eigenvalue_Discrete}
a_h(\bar{u}_h,v_h)=\bar{\lambda}_h b_h(\bar{u}_h,v_h),\quad\ \  \ \forall v_h\in V_h,
\end{eqnarray}
where
\begin{eqnarray}
a_h(u_h,v_h) &=& (A\nabla_wu_h,\nabla_wv_h)_{\mathcal{T}_h}=\sum_{K\in \mathcal{T}_h}(A\nabla_wu_h,\nabla_wv_h)_K, \\
b_h(u_h,v_h) &=& (u_0,v_0)_{\mathcal{T}_h} = \sum_{K\in \mathcal{T}_h}(u_0,v_0)_K.
\end{eqnarray}
Based on the bilinear form $a_h(\cdot,\cdot)$, we can define the following discrete 
norm on the space $V_h$ as follows
\begin{eqnarray}\label{Nomr_a}
\left\|v\right\|_{a,h} = \sqrt{a_h(v,v)},\ \ \ \ \forall v\in V_h.
\end{eqnarray}
We can also define the semi-norm $\left\|\cdot\right\|_{b,h}$ by the bilinear form $b_h(\cdot,\cdot)$
on the space $V_h$
\begin{eqnarray}\label{Nomr_b}
\left\|w\right\|_{b,h} = \sqrt{b_h(w,w)},\ \ \ \ \forall w\in V_h.
\end{eqnarray}
From \cite{BabuskaOsborn_1989,BabuskaOsborn_Book}, we obtain
$$0<\bar{\lambda}_{1,h}\leq \bar{\lambda}_{2,h}\leq \cdots \leq \bar{\lambda}_{k,h} \leq \cdots \leq \bar{\lambda}_{N_h,h},$$
and corresponding eigenfunctions
\begin{eqnarray}\label{Discrete_Eigenfunctions}
\bar{u}_{1,h}, \bar{u}_{2,h}, \cdots, \bar{u}_{k,h}, \cdots, \bar{u}_{N_h,h},
\end{eqnarray}
where $a_h(\bar{u}_{i,h},\bar{u}_{j,h})=\delta_{ij}$, $1\leq i,j\leq N_h$ ($N_h$ is
the dimension of the finite element space $V_h$).

For the following analysis in this paper, we define $\mu_i=1/\lambda_i$ for $i=1,2,\cdots$, and
$\bar\mu_{i, h} = 1/\bar\lambda_{i,h}$ for $i=1, \cdots, N_h$.

In order to state the error estimates for the
eigenpair approximation by the WG finite element method,
we define the WG finite element projection $\mathcal P_h: V\mapsto V_h$ as follows
\begin{eqnarray}\label{Energy_Projection}
a_h(\mathcal P_h u, v_h) = \lambda b_h(u,v_h),\ \ \ \ \forall v_h\in V_h.
\end{eqnarray}
It is obvious that the finite element projection operator $\mathcal P_h$ has following error estimates.

\begin{lemma}(\cite{WangYe})\label{Lemma_Error_Estimate_WG}
Assume the source equation corresponding to the eigenvalue problem has $H^{1+s}(\Omega)$ regularity and
the eigenfunction $u$ of (\ref{eigenvalue_problem}) belongs to $H^{m+1}(\Omega)$ and $0\leq m\leq r+1$.
Then the following error estimates hold
\begin{eqnarray}
\left\|Q_hu-\mathcal P_h u\right\|_{a,h} &\leq &C_1h^m\|u\|_{m+1},\label{Delta_V_h_P_h}\\
\left\|Q_hu-\mathcal P_h u\right\|_{b,h} &\leq& C_2h^{m+s}\|u\|_{m+1}.\label{Aubin_Nitsche_Estimate}
\end{eqnarray}
\end{lemma}
%-------------------------------------------------------------------------------------------------
Before stating error estimates of the WG finite element method for the eigenvalue problem,
we introduce the following lemma.
\begin{lemma}\label{Strang_Lemma}
For any eigenpair $(\lambda,u)$ of (\ref{weak_eigenvalue_problem}), the following equality holds
\begin{eqnarray*}\label{Strang_Equality}
(\bar\lambda_{j,h}-\lambda)b_h(\mathcal P_hu,\bar u_{j,h})
=\lambda b_h(u-\mathcal P_hu,\bar u_{j,h}),\ \ \ j = 1, \cdots, N_h.
\end{eqnarray*}
\end{lemma}
%-------------------------------------------------------------------------------------------------
\begin{proof}
Since $-\lambda b_h(\mathcal P_h u,\bar u_{j,h})$ appears on both sides, we only need to prove that
\begin{eqnarray*}
\bar\lambda_{j,h}b_h(\mathcal P_h u,\bar u_{j,h})=\lambda b_h(u,\bar u_{j,h}).
\end{eqnarray*}
From (\ref{weak_eigenvalue_problem}), (\ref{Weak_Eigenvalue_Discrete}) and (\ref{Energy_Projection}),
the following equalities hold
\begin{eqnarray*}
\bar\lambda_{j,h}b_h(\mathcal P_h u,\bar u_{j,h}) = a_h(\mathcal P_h u,\bar u_{j,h})= \lambda b(u,\bar u_{j,h}).
\end{eqnarray*}
Then the proof is completed.
\end{proof}

Now, let us consider the error estimates for the first
$k$ eigenpair approximations associated with $\bar\lambda_{1,h}\leq \cdots \leq \bar\lambda_{k,h}$.
\begin{theorem}\label{Error_Estimate_Theorem_k}
Let us define the spectral projection $\bar F_{k,h}: V_h\mapsto  {\rm span}\{\bar u_{1,h}, \cdots, \bar u_{k,h}\}$
as follows
\begin{eqnarray}\label{Definition_Spectral_Projection}
a_h(\bar F_{k,h}w_h, \bar u_{i,h}) = a_h(w_h, \bar u_{i,h}), \ \ \ i=1, \cdots, k\ \ {\rm for}\ w\in V_h.
\end{eqnarray}
Then the associated exact eigenfunctions $u_1, \cdots, u_k$ of eigenvalue problem (\ref{weak_eigenvalue_problem}) have the following error estimates
\begin{eqnarray}\label{Theo Energy_Error_Estimate_k}
\left\|Q_h u_i - \bar F_{k,h}Q_h u_i\right\|_{a,h} \leq 2\|Q_hu_i-\mathcal P_h u_i\|_{a,h}
+ \frac{\sqrt{\bar\mu_{k+1,h}}}{\delta_{k,i,h}}\left\|Q_hu_i-\mathcal P_hu_i\right\|_{b,h}, \ \  1\leq i\leq k,
\end{eqnarray}
where $\delta_{k,i,h}$ is defined as follows
\begin{eqnarray}\label{Definition_Delta_k_i_0}
\delta_{k,i,h} = \min_{k<j\leq N_h}\left|\frac{1}{\bar\lambda_{j,h}}-\frac{1}{\lambda_i}\right|.
\end{eqnarray}
Furthermore, these $k$ exact eigenfunctions have the following error estimate in $\left\|\cdot\right\|_{b,h}$-norm
\begin{eqnarray}\label{Theo L2_Error_Estimate_k_0}
\left\|Q_hu_i-\bar F_{k,h} Q_hu_i\right\|_{b,h} \leq \left(2+\frac{\bar\mu_{k+1,h}}{\delta_{k,i,h}}\right)
\left\|Q_hu_i-\mathcal P_hu_i\right\|_{b,h}, \ \ \ 1\leq i\leq k.
\end{eqnarray}
\end{theorem}
\begin{proof}
Since $(I-\bar F_{k,h})\mathcal P_hu_i\in V_h$ and
$(I-\bar F_{k,h})\mathcal P_hu_i\in {\rm span}\{\bar u_{k+1,h},\cdots, \bar u_{N_h,h}\}$,
the following orthogonal expansion holds
\begin{eqnarray}\label{Orthogonal_Decomposition_k}
(I-\bar F_{k,h})\mathcal P_hu_i=\sum_{j=k+1}^{N_h}\alpha_j\bar u_{j,h},
\end{eqnarray}
where $\alpha_j=a_h(\mathcal P_hu_i,\bar u_{j,h})$. From Lemma \ref{Strang_Lemma}, we have
\begin{eqnarray}\label{Alpha_Estimate}
\alpha_j&=&a_h(\mathcal P_hu_i,\bar u_{j,h}) = \bar\lambda_{j,h} b_h\big(\mathcal P_hu_i,\bar u_{j,h}\big)
=\frac{\bar\lambda_{j,h}\lambda_i}{\bar\lambda_{j,h}-\lambda_i}b_h\big(u_i-\mathcal P_hu_i,\bar u_{j,h}\big)\nonumber\\
&=&\frac{1}{\mu_i-\bar\mu_{j,h}} b_h\big(u_i-\mathcal P_hu_i,\bar u_{j,h}\big).
\end{eqnarray}
From the orthogonal property of eigenfunctions $\bar u_{1,h},\cdots, \bar u_{N_h,h}$, we acquire
\begin{eqnarray*}
1 = a_h(\bar u_{j,h},\bar u_{j,h}) = \bar\lambda_{j,h} b_h(\bar u_{j,h},\bar u_{j,h})
= \bar\lambda_{j,h}\left\|\bar u_{j,h}\right\|_{b,h}^2,
\end{eqnarray*}
which leads to the following property
\begin{eqnarray}\label{Equality_u_j}
\left\|\bar u_{j,h}\right\|_{b,h}^2=\frac{1}{\bar\lambda_{j,h}}=\bar\mu_{j,h}.
\end{eqnarray}
Because of (\ref{Weak_Eigenvalue_Discrete}) and the definitions of eigenfunctions $\bar u_{1,h},\cdots, \bar u_{N_h,h}$,
we obtain the following equalities
\begin{eqnarray}\label{Orthonormal_Basis}
a_h(\bar u_{j,h},\bar u_{k,h})=\delta_{jk},
\ \ \ \ \ b_h\left(\frac{\bar u_{j,h}}{\left\|\bar u_{j,h}\right\|_{b,h}},
\frac{\bar u_{k,h}}{\left\|\bar u_{k,h}\right\|_{b,h}}\right)=\delta_{jk},\ \ \ 1\leq j,k\leq N_h.
\end{eqnarray}
Then due to (\ref{Orthogonal_Decomposition_k}), (\ref{Alpha_Estimate}), (\ref{Equality_u_j}) and (\ref{Orthonormal_Basis}),
we have the following estimates
\begin{eqnarray}\label{Equality_4_i}
&&\left\|(I-\bar F_{k,h})\mathcal P_hu_i\right\|_{a,h}^2 = \left\|\sum_{j=k+1}^{N_h}\alpha_j\bar u_{j,h}\right\|_{a,h}^2
= \sum_{j=k+1}^{N_h}\alpha_j^2\nonumber\\
&&=\sum_{j=k+1}^{N_h} \left(\frac{1}{\mu_i-\bar\mu_{j,h}}\right)^2 b_h\big(u_i-\mathcal P_hu_i,\bar u_{j,h}\big)^2\nonumber\\
&&\leq\frac{1}{\delta_{k,i,h}^2}\sum_{j=k+1}^{N_h}\left\|\bar u_{j,h}\right\|_{b,h}^2
b_h\left(u_i-\mathcal P_hu_i,\frac{\bar u_{j,h}}{\left\|\bar u_{j,h}\right\|_{b,h}}\right)^2\nonumber\\
&&=\frac{1}{\delta_{k,i,h}^2}\sum_{j=k+1}^{N_h}\bar\mu_{j,h}
b_h\left(u_i-\mathcal P_hu_i,\frac{\bar u_{j,h}}{\left\|\bar u_{j,h}\right\|_{b,h}}\right)^2\nonumber\\
&&
\leq \frac{\bar\mu_{k+1,h}}{\delta_{k,i,h}^2}\sum_{j=k+1}^{N_h}
b_h\left(u_i-\mathcal P_hu_i,\frac{\bar u_{j,h}}{\left\|\bar u_{j,h}\right\|_{b,h}}\right)^2 \nonumber\\
&&= \frac{\bar\mu_{k+1,h}}{\delta_{k,i,h}^2}\sum_{j=k+1}^{N_h}
b_h\left(Q_hu_i-\mathcal P_hu_i,\frac{\bar u_{j,h}}{\left\|\bar u_{j,h}\right\|_{b,h}}\right)^2
\leq \frac{\bar\mu_{k+1,h}}{\delta_{k,i,h}^2}\left\|Q_hu_i-\mathcal P_hu_i\right\|_{b,h}^2,
\end{eqnarray}
where the last inequality holds since $\frac{\bar u_{1,h}}{\left\|\bar u_{1,h}\right\|_{b,h}}$, $\cdots$,
$\frac{\bar u_{j,h}}{\left\|\bar u_{j,h}\right\|_{b,h}}$ constitute an orthonormal basis for the space $V_h$
in the sense of the inner product $b_h(\cdot, \cdot)$.

From (\ref{Equality_4_i}),  the following inequality holds
\begin{eqnarray}\label{Equality_5_k}
\left\|(I-\bar F_{k,h})\mathcal P_hu_i\right\|_{a,h}
\leq\frac{\sqrt{\bar\mu_{k+1,h}}}{\delta_{k,i,h}}\left\|Q_hu_i-\mathcal P_hu_i\right\|_{b,h}.
\end{eqnarray}
From (\ref{Equality_5_k}), $\|\bar F_{k,h}\|_{a,h}\leq 1$ and the triangle inequality,  it follows that
\begin{eqnarray*}
&&\left\|Q_hu_i-\bar F_{k,h}Q_hu_i\right\|_{a,h}= \left\|Q_hu_i-\mathcal P_hu_i\right\|_{a,h}
+\left\|(I-\bar F_{k,h})\mathcal P_hu_i\right\|_{a,h} + \left\|\bar F_{k,h}(\mathcal P_h-Q_h)u_i\right\|_{a,h}\nonumber\\
&&\leq  \left\|Q_hu_i-\mathcal P_hu_i\right\|_{a,h} + \left\|(I-\bar F_{k,h})\mathcal P_hu_i\right\|_{a,h} +
\left\|\bar F_{k,h}\right\|_{a,h}\left\|(\mathcal P_h-Q_h)u_i\right\|_{a,h}\nonumber\\
&&\leq 2\|Q_hu_i-\mathcal P_h u_i\|_{a,h} + \frac{\sqrt{\bar\mu_{k+1,h}}}{\delta_{k,i,h}}\left\|Q_hu_i-\mathcal P_hu_i\right\|_{b,h}.
\end{eqnarray*}
This is the desired result (\ref{Theo Energy_Error_Estimate_k}).
%-----------------------------------------------------------------------------------------------------

Similarly, with the help of  (\ref{Orthogonal_Decomposition_k}),
(\ref{Alpha_Estimate}), (\ref{Equality_u_j}) and (\ref{Orthonormal_Basis}),
we have following estimates
\begin{eqnarray*}
&&\left\|(I-\bar F_{k,h})\mathcal P_hu_i\right\|_{b,h}^2 = \left\|\sum_{j=k+1}^{N_h}\alpha_j\bar u_{j,h}\right\|_{b,h}^2
= \sum_{j=k+1}^{N_h}\alpha_j^2\left\|\bar u_{j,h}\right\|_{b,h}^2\nonumber\\
&&=\sum_{j=k+1}^{N_h} \left(\frac{1}{\mu_i-\bar\mu_{j,h}}\right)^2
b_h\big(u_i-\mathcal P_hu_i,\bar u_{j,h}\big)^2\left\|\bar u_{j,h}\right\|_{b,h}^2\nonumber\\
&&\leq\frac{1}{\delta_{k,i,h}^2}\sum_{j=k+1}^{N_h}\left\|\bar u_{j,h}\right\|_{b,h}^4\
b_h\left(u_i-\mathcal P_hu_i,\frac{\bar u_{j,h}}{\left\|\bar u_{j,h}\right\|_{b,h}}\right)^2\nonumber\\
&&=\frac{1}{\delta_{k,i,h}^2}\sum_{j=k+1}^{N_h}\bar\mu_{j,h}^2
b_h\left(Q_hu_i-\mathcal P_hu_i,\frac{\bar u_{j,h}}{\left\|\bar u_{j,h}\right\|_{b,h}}\right)^2
\leq \frac{\bar\mu_{k+1,h}^2}{\delta_{k,i,h}^2}\left\|Q_hu_i-\mathcal P_hu_i\right\|_{b,h}^2,
\end{eqnarray*}
which leads to the inequality
\begin{eqnarray}\label{Equality_8_k}
\left\|(I-\bar F_{k,h})\mathcal P_hu_i\right\|_{b,h} \leq \frac{\bar\mu_{k+1,h}}{\delta_{k,i,h}}
\left\|Q_hu_i-\mathcal P_hu_i\right\|_{b,h}.
\end{eqnarray}

From the definition of spectral projection (\ref{Definition_Spectral_Projection}), for any $w\in V_h$,
we have
\begin{eqnarray*}
\bar\lambda_{i,h}b_h(\bar F_{k,h}w, \bar u_{i,h}) = a_h(\bar F_{k,h}w, \bar u_{i,h})
= a_h(w, \bar u_{i,h}) = \bar \lambda_{i,h}b_h(w, \bar u_{i,h}), \ \ \ i=1, \cdots, k.
\end{eqnarray*}
This means the following equation holds
\begin{eqnarray*}
b_h(\bar F_{k,h}w, \bar u_{i,h})  = b_h(w, \bar u_{i,h}), \ \ \ i=1, \cdots, k,\ \ \forall w\in V_h,
\end{eqnarray*}
which leads to $\|\bar F_{k,h}\|_{b,h}\leq 1$.

From (\ref{Equality_8_k}), $\|\bar F_{k,h}\|_{b,h}\leq 1$ and the triangle inequality, we find the
following error estimates for the eigenfunction approximations in the $\left\|\cdot\right\|_{b,h}$-norm
\begin{eqnarray*}\label{Inequality_11}
&&\left\|Q_hu_i-\bar F_{k,h}Q_hu_i\right\|_{b,h}\leq \left\|Q_hu_i-\mathcal P_hu_i\right\|_{b,h}
+ \left\|(I-\bar F_{k,h}) \mathcal P_hu_i\right\|_{b,h} + \left\|\bar F_{k,h}(\mathcal P_hu_i-Q_hu_i)\right\|_{b,h}\nonumber\\
&&\leq\left(1+\|\bar F_{k,h}\|_{b,h}\right) \left\|\mathcal P_hu_i-Q_hu_i \right\|_{b,h} +
\left\|(I-\bar F_{k,h})\mathcal P_hu_i\right\|_{b,h}\nonumber\\
&&\leq\left(2+\frac{\bar\mu_{k+1,h}}{\delta_{k,i,h}}\right)
\left\|Q_hu_i-\mathcal P_hu_i\right\|_{b,h}.
\end{eqnarray*}
This is the second desired result (\ref{Theo L2_Error_Estimate_k_0}) and the proof is completed.
\end{proof}
For the sake of simplicity in notation and to make sense of the estimates 
(\ref{Theo Energy_Error_Estimate_k}) and (\ref{Theo L2_Error_Estimate_k_0}), 
we assume that the eigenvalue gap $\delta_{k,i,h}$ has a uniform lower bound, 
which is represented by $\delta_{k,i}$ (which can be understood as the ``true" 
separation of the eigenvalues $\lambda_1, \cdots, \lambda_k$ from the unwanted eigenvalues) 
in the following sections of this paper.
When the mesh size is sufficiently small, this assumption makes sense.
Based on Theorem \ref{Error_Estimate_Theorem_k} and the convergence consequences of 
the WG finite element method for boundary value problems, we then acquire 
the following convergence order.
\begin{corollary}\label{Error_Estimate_Corollary_k}
Under the conditions of Lemma \ref{Lemma_Error_Estimate_WG}, Theorem \ref{Error_Estimate_Theorem_k}
and $\delta_{k,i,h}$ having a uniform lower bound $\delta_{k,i}$,  the following error estimates hold
\begin{eqnarray}
&&\left\|Q_hu_i - \bar F_{k,h} Q_hu_i\right\|_{a,h} \leq C_3h^m\|u\|_{m+1}, \ \ \ \  \ 1\leq i\leq k,\label{Energy_Error_Estimate_k}\\
&&\left\|Q_hu_i-\bar F_{k,h}Q_hu_i\right\|_{b,h} \leq C_4h^{m+s}\|u\|_{m+1},\ \ \ 1\leq i\leq k.\label{L2_Error_Estimate_k}
\end{eqnarray}
\end{corollary}
%-------------------------------------------------------------------------------------------------
The following theorem gives the error estimates for the one eigenpair approximation and the proof is similar
to that of Theorem \ref{Error_Estimate_Theorem_k}.
%-------------------------------------------------------------------------------------------------
\begin{theorem}\label{Error_Estimate_Theorem_Old}
Let  $(\lambda,u)$ denote an exact eigenpair of the eigenvalue problem (\ref{weak_eigenvalue_problem}).
Assume the eigenpair approximation $(\bar\lambda_{i,h},\bar u_{i,h})$ has the property that
$\bar\mu_{i,h}=1/\bar\lambda_{i,h}$ is the closest to $\mu=1/\lambda$.
The corresponding spectral projector $E_{i,h}: V_h\mapsto {\rm span}\{\bar u_{i,h}\}$
is  defined as follows
\begin{eqnarray*}
a_h(E_{i,h}w,\bar u_{i,h}) = a_h(w,\bar u_{i,h}),\ \ \ \ {\rm for}\  w\in V_h.
\end{eqnarray*}
Then the following error estimate holds
\begin{eqnarray}\label{Energy_Error_Estimate_Old}
\left\|Q_hu-E_{i,h}Q_hu\right\|_{a,h}&\leq& 2\|Q_hu-\mathcal P_h u\|_{a,h}
+ \frac{\sqrt{\bar\mu_{1,h}}}{\delta_{\lambda,h}}\left\|Q_hu-\mathcal P_hu\right\|_{b,h},
\end{eqnarray}
where  $\delta_{\lambda,h}$ is defined as follows
\begin{eqnarray}\label{Definition_Delta}
\delta_{\lambda,h} &:=& \min_{j\neq i}|\bar\mu_{j,h}-\mu|=\min_{j\neq i} \left|\frac{1}{\bar\lambda_{j,h}}-\frac{1}{\lambda}\right|.
\end{eqnarray}
Furthermore, the eigenfunction approximation $\bar u_{i,h}$ has the following
error estimate in $\left\|\cdot\right\|_{b,h}$-norm
\begin{eqnarray}\label{L2_Error_Estimate_Old}
\|Q_hu - E_{i,h}Q_hu\|_{b,h} &\leq &  \left(2+\frac{\bar\mu_{1,h}}{\delta_{\lambda,h}}\right)
\left\|Q_hu_i-\mathcal P_hu_i\right\|_{b,h}.
\end{eqnarray}
\end{theorem}
\begin{proof}
Since $(I-E_{i,h})\mathcal P_hu\in V_h$ and
$(I-E_{i,h})\mathcal P_hu\in {\rm span}\{\bar u_{1,h}, \cdots, \bar u_{i-1,h}, \bar u_{i+1,h},\cdots, \bar u_{N_h,h}\}$,
the following orthogonal expansion holds
\begin{eqnarray}\label{Orthogonal_Decomposition_1}
(I-E_{i,h})\mathcal P_hu=\sum_{j\neq i}\alpha_j\bar u_{j,h},
\end{eqnarray}
where $\alpha_j=a_h(\mathcal P_hu,\bar u_{j,h})$ has the same equality (\ref{Alpha_Estimate}).

Then due to (\ref{Alpha_Estimate}), (\ref{Equality_u_j}),  (\ref{Orthonormal_Basis}) and  (\ref{Orthogonal_Decomposition_1}), we have following estimates
\begin{eqnarray}\label{Equality_4_i_1}
&&\left\|(I-E_{i,h})\mathcal P_hu\right\|_{a,h}^2 = \left\|\sum_{j\neq i}\alpha_j\bar u_{j,h}\right\|_{a,h}^2
= \sum_{j\neq i}\alpha_j^2\nonumber\\
&&=\sum_{j\neq i} \left(\frac{1}{\mu-\bar\mu_{j,h}}\right)^2 b_h\big(u-\mathcal P_hu,\bar u_{j,h}\big)^2
\leq\frac{1}{\delta_{\lambda,h}^2}\sum_{j\neq i}\left\|\bar u_{j,h}\right\|_{b,h}^2
b_h\left(u-\mathcal P_hu,\frac{\bar u_{j,h}}{\left\|\bar u_{j,h}\right\|_{b,h}}\right)^2\nonumber\\
&&=\frac{1}{\delta_{\lambda,h}^2}\sum_{j\neq i}\bar\mu_{j,h}
b_h\left(u-\mathcal P_hu,\frac{\bar u_{j,h}}{\left\|\bar u_{j,h}\right\|_{b,h}}\right)^2\nonumber\\
&&
\leq \frac{\bar\mu_{1,h}}{\delta_{\lambda,h}^2}\sum_{j\neq i}
b_h\left(u-\mathcal P_hu,\frac{\bar u_{j,h}}{\left\|\bar u_{j,h}\right\|_{b,h}}\right)^2
= \frac{\bar\mu_{1,h}}{\delta_{\lambda,h}^2}\sum_{j\neq i}
b_h\left(Q_hu-\mathcal P_hu,\frac{\bar u_{j,h}}{\left\|\bar u_{j,h}\right\|_{b,h}}\right)^2 \nonumber\\
&&\leq \frac{\bar\mu_{1,h}}{\delta_{\lambda,h}^2}\left\|Q_hu-\mathcal P_hu\right\|_{b,h}^2,
\end{eqnarray}
where the last inequality holds since $\frac{\bar u_{1,h}}{\left\|\bar u_{1,h}\right\|_{b,h}}$, $\cdots$,
$\frac{\bar u_{j,h}}{\left\|\bar u_{j,h}\right\|_{b,h}}$ constitute an orthonormal 
basis for the space $V_h$
in the sense of the inner product $b_h(\cdot, \cdot)$.

From (\ref{Equality_4_i_1}),  the following inequality holds
\begin{eqnarray}\label{Equality_5_1}
\left\|(I-E_{i,h})\mathcal P_hu\right\|_{a,h}
\leq\frac{\sqrt{\bar\mu_{1,h}}}{\delta_{\lambda,h}}\left\|Q_hu-\mathcal P_hu\right\|_{b,h}.
\end{eqnarray}
From (\ref{Equality_5_1}), $\|E_{i,h}\|_{a,h}\leq 1$ and the triangle inequality,  it follows that
\begin{eqnarray*}
&&\left\|Q_hu-E_{i,h}Q_hu\right\|_{a,h}= \left\|Q_hu-\mathcal P_hu\right\|_{a,h}+
\left\|(I-E_{i,h})\mathcal P_hu\right\|_{a,h} + \left\|E_{i,h}(\mathcal P_h-Q_h)u\right\|_{a,h}\nonumber\\
&&\leq  \left\|Q_hu-\mathcal P_hu\right\|_{a,h} + \left\|(I-E_{i,h})\mathcal P_hu\right\|_{a,h} +
\left\|E_{i,h}\right\|_{a,h}\left\|(\mathcal P_h-Q_h)u\right\|_{a,h}\nonumber\\
&&\leq 2\|Q_hu-\mathcal P_h u\|_{a,h} + \frac{\sqrt{\bar\mu_{1,h}}}{\delta_{\lambda,h}}\left\|Q_hu-\mathcal P_hu\right\|_{b,h}.
\end{eqnarray*}
This is the desired result (\ref{Energy_Error_Estimate_Old}).
%-----------------------------------------------------------------------------------------------------

Similarly, with the help of  (\ref{Alpha_Estimate}), (\ref{Equality_u_j}), (\ref{Orthonormal_Basis})
and (\ref{Orthogonal_Decomposition_1}), we have the following estimates
\begin{eqnarray*}
&&\left\|(I-E_{i,h})\mathcal P_hu\right\|_{b,h}^2 = \left\|\sum_{j\neq i}\alpha_j\bar u_{j,h}\right\|_{b,h}^2
= \sum_{j\neq i}\alpha_j^2\left\|\bar u_{j,h}\right\|_{b,h}^2\nonumber\\
&&=\sum_{j\neq i} \left(\frac{1}{\mu-\bar\mu_{j,h}}\right)^2
b_h\big(u-\mathcal P_hu,\bar u_{j,h}\big)^2\left\|\bar u_{j,h}\right\|_{b,h}^2\nonumber\\
&&\leq\frac{1}{\delta_{\lambda, h}^2}\sum_{j\neq i}\left\|\bar u_{j,h}\right\|_{b,h}^4\
b_h\left(u-\mathcal P_hu,\frac{\bar u_{j,h}}{\left\|\bar u_{j,h}\right\|_{b,h}}\right)^2\nonumber\\
&&=\frac{1}{\delta_{\lambda,h}^2}\sum_{j\neq i}\bar\mu_{j,h}^2
b_h\left(Q_hu-\mathcal P_hu,\frac{\bar u_{j,h}}{\left\|\bar u_{j,h}\right\|_{b,h}}\right)^2
\leq \frac{\bar\mu_{1,h}^2}{\delta_{\lambda,h}^2}\left\|Q_hu-\mathcal P_hu\right\|_{b,h}^2,
\end{eqnarray*}
which leads to the inequality
\begin{eqnarray}\label{Equality_8_1}
\left\|(I-E_{i,h})\mathcal P_hu\right\|_{b,h} \leq \revise{\frac{\bar\mu_{1,h}}{\delta_{\lambda,h}}}
\left\|Q_hu-\mathcal P_hu\right\|_{b,h}.
\end{eqnarray}
Similarly to the proof of Theorem \ref{Error_Estimate_Theorem_k}, we also have $\|E_{i,h}\|_{b,h}\leq 1$.
Then from (\ref{Equality_8_1}) and the triangle inequality, we find the
following error estimates for the eigenfunction approximations in the $\left\|\cdot\right\|_{b,h}$-norm
\begin{eqnarray*}\label{Inequality_11}
&&\left\|Q_hu-E_{i,h}Q_hu\right\|_{b,h}\leq \left\|Q_hu-\mathcal P_hu\right\|_{b,h} + \left\|(I-E_{i,h})
\mathcal P_hu\right\|_{b,h} + \left\|E_{i,h}(\mathcal P_hu-Q_hu)\right\|_{b,h}\nonumber\\
&&\leq\left(1+\|E_{i,h}\|_{b,h}\right) \left\|\mathcal P_hu-Q_hu \right\|_{b,h} +
\left\|(I-E_{i,h})\mathcal P_hu\right\|_{b,h}\leq\left(2+\frac{\bar\mu_{1,h}}{\delta_{\lambda,h}}\right)
\left\|Q_hu-\mathcal P_hu\right\|_{b,h}.
\end{eqnarray*}
This is the second desired result (\ref{L2_Error_Estimate_Old}) and the proof is completed.
\end{proof}
%--------------------------------------------------------------------------------------------
Likewise, for the sake of simplicity in notation and to make sense of the estimates (\ref{Energy_Error_Estimate_Old}) and (\ref{L2_Error_Estimate_Old}), 
we assume that the eigenvalue gap $\delta_{\lambda,h}$ defined by (\ref{Definition_Delta}) 
equally has a uniform lower bound, indicated by $\delta_\lambda$, which can be understood 
as the ``true" separation of the eigenvalue $\lambda$ from others in the following sections 
of this paper. When the mesh size is small enough, this assumption is also reasonable.
Next, we have the following convergence result for the eigenvalue problems using the 
WG finite element method, which is based on Theorem \ref{Error_Estimate_Theorem_Old}.
\begin{corollary}\label{Error_Estimate_Corollary}
Under the conditions of Lemma \ref{Lemma_Error_Estimate_WG}, Theorem \ref{Error_Estimate_Theorem_Old} and
$\delta_{\lambda,h}$ having a uniform lower bound $\delta_\lambda$, the following error estimates hold
\begin{eqnarray}
\left\|Q_hu-E_{i,h}Q_hu\right\|_{a,h}&\leq& C_5h^m\|u\|_{m+1},\label{Energy_Error_Estimate}\\
\left\|Q_hu-E_{i,h}Q_hu\right\|_{b,h}&\leq& C_6h^{m+s}\|u\|_{m+1}. \label{L2_Error_Estimate}
\end{eqnarray}
\end{corollary}
\begin{remark}
The convergence analysis of the WG finite element method for eigenvalue problems
has been provided in \cite{ZhaiXieZhangZhang_WG}. Compared with the results there, the convergence results in
Theorems \ref{Error_Estimate_Theorem_k} and \ref{Error_Estimate_Theorem_Old}
are sharper and gives the explicit dependence of the included constants on the eigenvalue distributions.
\end{remark}

\section{Augmented subspace method and its error estimates}\label{Section_3}
The augmented subspace techniques for the WG eigenvalue problem (\ref{Weak_Eigenvalue_Discrete}) are first laid out in this section.
These schemes involve solving the eigenvalue problem on the augmented subspace $V_{H,h}$, which is generated by the coarse conforming linear finite element space $W_H$, and a WG finite element function in the fine finite element space $V_h$. They also involve solving the auxiliary linear boundary value problem in the fine finite element space $V_h$.
Next, the related analysis of convergence for these augmented subspace schemes is addressed.

As in \cite{Xie_BIT}, we first create a coarse mesh $\mathcal{T}_H$ with the mesh size $H$, and the corresponding conforming linear finite element space $W_H$ is defined on the mesh $\mathcal{T}_H$. This allows us to design the augmented subspace technique. The coarse conforming linear finite element space $W_H$ is a subspace of the fine \comm{WG} finite element space $V_h$, which is defined on the fine mesh $\mathcal T_h$. This is because, for the sake of simplicity, we assume in this paper that the coarse mesh $\mathcal T_H$ and the fine mesh $\mathcal T_h$ have nested properties.

For the positive integer $\ell$ and given eigenfunction approximations $u_{1,h}^{(\ell)},\cdots, u_{k,h}^{(\ell)}$
which are the approximations for the first $k$ eigenfunctions $\bar u_{1,h},\cdots, \bar u_{k,h}$
of (\ref{Weak_Eigenvalue_Discrete}), we can do the augmented subspace iteration step
which is defined by Algorithm \ref{Algorithm_k} to improve the accuracy of $u_{1,h}^{(\ell)},\cdots, u_{k,h}^{(\ell)}$.
\begin{algorithm}[hbt!]
\caption{Augmented subspace method for the first $k$ eigenpairs }\label{Algorithm_k}
\begin{enumerate}
\item For $\ell=1$, we define $\widehat u_{i,h}^{(\ell)}=u_{i,h}^{(\ell)}$, $i=1,\cdots, k$, and
the augmented subspace $V_{H,h}^{(\ell)} = W_H +{\rm span}\{\widehat u_{1,h}^{(\ell)},
\cdots, \widehat u_{k,h}^{(\ell)}\}$.
Then solve the following eigenvalue problem:
Find $(\lambda_{i,h}^{(\ell)},u_{i,h}^{(\ell)})\in \mathbb R\times V_{H,h}^{(\ell)}$
such that $a_h(u_{i,h}^{(\ell)},u_{i,h}^{(\ell)})=1$ and
\begin{equation}\label{Aug_Eigenvalue_Problem_k_1}
a_h(u_{i,h}^{(\ell)},v_{H,h}) = \lambda_{i,h}^{(\ell)}b_h(u_{i,h}^{(\ell)},v_{H,h}),
\ \ \ \ \ \forall v_{H,h}\in V_{H,h}^{(\ell)},\ \ \ i=1, \cdots, k.
\end{equation}

\item Solve the following linear boundary value problems:
Find $\widehat{u}_{i,h}^{(\ell+1)}\in V_h$ such that
\begin{equation}\label{Linear_Equation_k}
a_h(\widehat{u}_{i,h}^{(\ell+1)},v_h) = \lambda_{i,h}^{(\ell)}b_h(u_{i,h}^{(\ell)},v_h),
\ \  \forall v_h\in V_h,\ \ \ i=1, \cdots, k.
\end{equation}

\item Define the augmented subspace $V_{H,h}^{(\ell+1)} = W_H +
{\rm span}\{\widehat{u}_{1,h}^{(\ell+1)}, \cdots, \widehat u_{k,h}^{(\ell+1)}\}$ and solve the
following eigenvalue problem:
Find $(\lambda_{i,h}^{(\ell+1)},u_{i,h}^{(\ell+1)})\in\mathbb R\times V_{H,h}^{(\ell+1)}$
such that $a_h(u_{i,h}^{(\ell+1)},u_{i,h}^{(\ell+1)})=1$ and
\begin{equation}\label{Aug_Eigenvalue_Problem_k}
a_h(u_{i,h}^{(\ell+1)},v_{H,h}) = \lambda_{i,h}^{(\ell+1)}b_h(u_{i,h}^{(\ell+1)},v_{H,h}),
\ \ \ \ \ \forall v_{H,h}\in V_{H,h}^{(\ell+1)},\ \ \ i=1, \cdots, k.
\end{equation}
Solve (\ref{Aug_Eigenvalue_Problem_k})  to obtain $(\lambda_{1,h}^{(\ell+1)},u_{1,h}^{(\ell+1)}), \cdots,
(\lambda_{k,h}^{(\ell+1)},u_{k,h}^{(\ell+1)})$.
\item Set $\ell=\ell+1$ and go to Step 2 for the next iteration until convergence.
\end{enumerate}
%Summarize the above two steps by defining
%\begin{eqnarray*}
%(\lambda_h^{(\ell+1)},u_h^{(\ell+1)}) =
%{\tt Correction}(V_H,V_h,\lambda_{h}^{(\ell)},u_{h}^{(\ell)}).
% \end{eqnarray*}
\end{algorithm}
%====================================================================================

For each $\ell$, it is easy to know, the eigenvalue problems (\ref{Aug_Eigenvalue_Problem_k_1})
and (\ref{Aug_Eigenvalue_Problem_k}) has the following eigenvalues 
\cite{BabuskaOsborn_1989,BabuskaOsborn_Book},
$$0< \lambda_{1,h}^{(\ell)}\leq  \lambda_{2,h}^{(\ell)}
\leq \cdots \leq \lambda_{k,h}^{(\ell)} \leq \cdots \leq  \lambda_{N_{H,h},h}^{(\ell)},$$
and corresponding eigenfunctions
\begin{eqnarray}\label{Discrete_Eigenfunctions_aug}
u_{1,h}^{(\ell)},  u_{2,h}^{(\ell)}, \cdots,  u_{k,h}^{(\ell)}, \cdots,  u_{N_{H,h},h}^{(\ell)},
\end{eqnarray}
where $N_{H,h}={\rm dim}V_{H,h}^{(\ell)} = N_H+k$ and $a_h(u_{i,h}^{(\ell)}, u_{j,h}^{(\ell)})=\delta_{ij}$, $1\leq i,j\leq N_{H,h}$.

From the min-max principle \cite{BabuskaOsborn_1989,BabuskaOsborn_Book}
and $V_{H,h}^{(\ell)}\subset V_h$, the eigenvalues
$\lambda_{1,h}^{(\ell)}, \cdots, \lambda_{N_{H,h},h}^{(\ell)}$ provide upper bounds for the first
$N_{H,h}$ eigenvalues of (\ref{Weak_Eigenvalue_Discrete})
\begin{eqnarray}\label{Upper_Bound_Result}
\bar\lambda_{i,h} \leq \lambda_{i,h}^{(\ell)},\ \ \   \bar\mu_{i,h} \geq \mu_{i,h}^{(\ell)},
\ \ \ \ 1\leq i\leq N_{H,h}.
\end{eqnarray}

Since the low dimensional augmented subspace $V_{H,h}^{(\ell)}$ is a subspace of the WG finite
element space $V_h$, the error estimates of eigenfunction approximations $u_{1,h}^{(\ell)}$, $\cdots$,
$u_{k,h}^{(\ell)}$ to the exact eigenfunctions $\bar u_{1,h}$, $\cdots$, $\bar u_{k,h}$ can be deduced
from the similar way of the conforming finite element method for the eigenvalue problem.

In order to give the error estimates for the augmented subspace method defined by Algorithm \ref{Algorithm_k},
we define the subspace projection $\mathcal P_{H,h}^{(\ell)}: V_h \mapsto V_{H,h}^{(\ell)}$ as follows
\begin{eqnarray}\label{Projection_aug}
a_h\left(\mathcal P_{H,h}^{(\ell)} w_h, v_{H,h}\right) = a_h\left(w_h, v_{H,h}\right),
\ \ \ \ \forall v_{H,h}\in V_{H,h}^{(\ell)},\ \ {\rm for}\ w_h\in V_h.
\end{eqnarray}
In order to give the error estimate of $\|w_h-\mathcal P_{H,h}^{(\ell)}w_h\|_{b,h}$,
we define the following quantity for error analysis:
\begin{eqnarray}\label{eta_a_H_Def}
\eta_a(W_H)&=&\sup_{\substack{ f\in L^2(\Omega)\\ \left\|f\right\|_{b,h}=1}}\inf_{w_H\in W_H}\left\|T_hf-w_H\right\|_{a,h},
\end{eqnarray}
where $T_h: L^2(\Omega)\mapsto V_h$ is defined as
\begin{equation}\label{laplace_source_operator}
a_h(T_hf,v) = b_h(f,v), \ \ \ \ \  \forall v \in V_h\ \ {\rm for}\  f \in L^2(\Omega).
\end{equation}

Then the projection operator $\mathcal P_{H,h}^{(\ell)}$ has following error estimates
\begin{eqnarray}
\|w_h - \mathcal P_{H,h}^{(\ell)}w_h\|_{a,h} &=& \inf_{v_{H,h}\in V_{H,h}^{(\ell)}}\|w_h-v_{H,h}\|_{a,h},
\ \ \ \ {\rm for}\ w_h\in V_h,
\label{Energy_Error_Estimate}\\
\|w_h-\mathcal P_{H,h}^{(\ell)}w_h\|_{b,h}  &\leq&  \eta_a(W_H) \|w_h-\mathcal P_{H,h}^{(\ell)} w_h\|_{a,h},
\ \ {\rm for}\ w_h\in V_h.  \label{L2_Energy_Estimate}
\end{eqnarray}

\begin{lemma}\label{Lemma_Error_Estimate_Subspace_k}
Let us define the spectral projection $F_{k,h}^{(m)}: V_h\mapsto {\rm span}\{u_{1,h}^{(m)}, \cdots, u_{k,h}^{(m)}\}$
for any integer $m \geq 1$ as follows
\begin{eqnarray}
a_h(F_{k,h}^{(m)}w, u_{i,h}^{(m)}) = a_h(w, u_{i,h}^{(m)}), \ \ \ i=1, \cdots, k\ \ {\rm for}\ w\in V_h.
\end{eqnarray}
Then the exact eigenfunctions $\bar u_{1,h},\cdots, \bar u_{k,h}$ of (\ref{Weak_Eigenvalue_Discrete})
and the eigenfunction
approximations $u_{1,h}^{(\ell+1)}$, $\cdots$,  $u_{k,h}^{(\ell+1)}$ from Algorithm \ref{Algorithm_k}
with the integer $\ell \geq 1$ have the following error estimate
\begin{eqnarray}\label{Error_Estimate_Inverse_aug}
\left\|\bar u_{i,h} -F_{k,h}^{(\ell)}\bar u_{i,h} \right\|_{a,h} \leq
\sqrt{1+\frac{\bar\mu_{k+1,h}}{(\delta_{k,i,h}^{(\ell)})^2}\eta_a^2(W_H)}
\left\|(I-\mathcal P_{H,h}^{(\ell)})\bar u_{i,h}\right\|_{a,h},
\end{eqnarray}
where $\delta_{k,i,h}^{(\ell)}$ is defined as follows
\begin{eqnarray}\label{Definition_Delta_k_i_aug}
\delta_{k,i,h}^{(\ell)} = \min_{k<j\leq N_h}\left|\frac{1}{\lambda_{j,h}^{(\ell)}}-\frac{1}{\bar\lambda_{i,h}}\right|.
\end{eqnarray}

Furthermore, the following $\left\|\cdot\right\|_{b,h}$-norm error estimate holds
\begin{eqnarray}\label{L2_Error_Estimate_Algorithm_1_aug}
&&\left\|\bar u_{i,h} -F_{k,h}^{(\ell)}\bar u_{i,h} \right\|_{b,h}\leq \bar\eta_a(W_H)\left\|\bar u_{i,h}-F_{k,h}^{(\ell)}\bar u_{i,h}\right\|_{a,h}.
\end{eqnarray}
where
\begin{eqnarray}\label{Definition_eta_bar}
\bar\eta_a(W_{H}) = \left(1+\frac{\bar\mu_{k+1,h}}{\delta_{k,i,h}^{(\ell)}}\right)\eta_a(W_H).
\end{eqnarray}

\end{lemma}
\begin{proof}
Since $(I-F_{k,h}^{(\ell)})\mathcal P_{H,h}^{(\ell)}\bar u_{i,h}\in V_{H,h}^{(\ell)}$ and
$(I-F_{k,h}^{(\ell)})\mathcal P_{H,h}^{(\ell)}\bar u_{i,h}\in {\rm span}\{u_{k+1,h}^{(\ell)},
\cdots, u_{N_{H,h},h}^{(\ell)}\}$,
the following orthogonal expansion holds
\begin{eqnarray}\label{Orthogonal_Decomposition_k_aug}
(I-F_{k,h}^{(\ell)})\mathcal P_{H,h}^{(\ell)}\bar u_{i,h}=\sum_{j=k+1}^{N_{H,h}}\alpha_ju_{j,h}^{(\ell)},
\end{eqnarray}
where  $\alpha_j=a_h(\mathcal P_{H,h}^{(\ell)}\bar u_{i,h},u_{j,h}^{(\ell)})$. From Lemma \ref{Strang_Lemma},
we have
\begin{eqnarray}\label{Alpha_Estimate_aug}
\alpha_j&=&a_h(\mathcal P_{H,h}^{(\ell)}\bar u_{i,h},u_{j,h}^{(\ell)}) = \lambda_{j,h}^{(\ell)}
b_h\big(\mathcal P_{H,h}^{(\ell)}\bar u_{i,h},u_{j,h}^{(\ell)}\big)
=\frac{\lambda_{j,h}^{(\ell)}\bar\lambda_{i,h}}{\lambda_{j,h}^{(\ell)}-\bar\lambda_{i,h}}
b_h\big(\bar u_{i,h}-\mathcal P_{H,h}^{(\ell)}\bar u_{i,h},u_{j,h}^{(\ell)}\big)\nonumber\\
&=&\frac{1}{\bar \mu_{i,h}- \mu_{j,h}^{(\ell)}} b_h\big(\bar u_{i,h}-\mathcal P_{H,h}^{(\ell)}\bar u_{i,h},u_{j,h}^{(\ell)}\big).
\end{eqnarray}
From the orthogonal property of eigenfunctions $u_{1,h}^{(\ell)},\cdots, \bar u_{N_{H,h},h}^{(\ell)}$, we have
\begin{eqnarray*}
1 = a_h(u_{j,h}^{(\ell)},u_{j,h}^{(\ell)}) = \lambda_{j,h}^{(\ell)}b_h(u_{j,h}^{(\ell)},u_{j,h}^{(\ell)})
= \lambda_{j,h}^{(\ell)}\left\|u_{j,h}^{(\ell)}\right\|_{b,h}^2,
\end{eqnarray*}
which leads to the following property
\begin{eqnarray}\label{Equality_u_j_aug}
\left\|u_{j,h}^{(\ell)}\right\|_{b,h}^2=\frac{1}{\lambda_{j,h}^{(\ell)}}= \mu_{j,h}^{(\ell)}.
\end{eqnarray}
Because of (\ref{Aug_Eigenvalue_Problem_k_1}),  (\ref{Aug_Eigenvalue_Problem_k}) and the definitions of eigenfunctions
$u_{1,h}^{(\ell)},\cdots, u_{N_{H,h},h}^{(\ell)}$, we obtain the following equalities
\begin{eqnarray}\label{Orthonormal_Basis_aug}
a_h(u_{j,h}^{(\ell)}, u_{k,h}^{(\ell)})=\delta_{jk},
\ \ \ \ \ b_h\left(\frac{u_{j,h}^{(\ell)}}{\left\|u_{j,h}^{(\ell)}\right\|_{b,h}},
\frac{u_{k,h}^{(\ell)}}{\left\|u_{k,h}^{(\ell)}\right\|_{b,h}}\right)=\delta_{jk},\ \ \ 1\leq j,k\leq N_{H,h}.
\end{eqnarray}
Then due to  (\ref{Upper_Bound_Result}), (\ref{Orthogonal_Decomposition_k_aug}), (\ref{Alpha_Estimate_aug}),
(\ref{Equality_u_j_aug}) and (\ref{Orthonormal_Basis_aug}), we have following estimates
\begin{eqnarray}\label{Equality_4_i_aug}
&&\left\|(I-F_{k,h}^{(\ell)})\mathcal P_{H,h}^{(\ell)}\bar u_{i,h}\right\|_{a,h}^2
= \left\|\sum_{j=k+1}^{N_{H,h}}\alpha_ju_{j,h}^{(\ell)}\right\|_{a,h}^2
= \sum_{j=k+1}^{N_{H,h}}\alpha_j^2\nonumber\\
&&=\sum_{j=k+1}^{N_{H,h}} \left(\frac{1}{\bar \mu_{i,h}- \mu_{j,h}^{(\ell)}}\right)^2
b_h\big(\bar u_{i,h}-\mathcal P_{H,h}^{(\ell)}\bar u_{i,h},u_{j,h}^{(\ell)}\big)^2\nonumber\\
&&
\leq\frac{1}{(\delta_{k,i,h}^{(\ell)})^2}\sum_{j=k+1}^{N_{H,h}}\left\|u_{j,h}^{(\ell)}\right\|_{b,h}^2
b_h\left(\bar u_{i,h}-\mathcal P_{H,h}^{(\ell)}\bar u_{i,h},\frac{u_{j,h}^{(\ell)}}{\left\|u_{j,h}^{(\ell)}\right\|_{b,h}}\right)^2\nonumber\\
&&=\frac{1}{(\delta_{k,i,h}^{(\ell)})^2}\sum_{j=k+1}^{N_{H,h}}\mu_{j,h}^{(\ell)}
b_h\left(\bar u_{i,h}-\mathcal P_{H,h}^{(\ell)}\bar u_{i,h},\frac{u_{j,h}^{(\ell)}}{\left\|u_{j,h}^{(\ell)}\right\|_b}\right)^2\nonumber\\
&&
\leq \frac{\mu_{k+1,h}^{(\ell)}}{(\delta_{k,i,h}^{(\ell)})^2}\sum_{j=k+1}^{N_{H,h}}
b_h\left(\bar u_{i,h}-\mathcal P_{H,h}^{(\ell)}\bar u_{i,h},\frac{u_{j,h}^{(\ell)}}{\left\|u_{j,h}^{(\ell)}\right\|_{b,h}}\right)^2\nonumber\\
&&\leq \frac{\mu_{k+1,h}^{(\ell)}}{(\delta_{k,i,h}^{(\ell)})^2}\left\|\bar u_{i,h}-\mathcal P_{H,h}^{(\ell)}\bar u_{i,h}\right\|_{b,h}^2,
\end{eqnarray}
where the last inequality holds since $\frac{ u_{1,h}^{(\ell)}}{\left\|u_{1,h}^{(\ell)}\right\|_b}$, $\cdots$,
$\frac{u_{N_{H,h},h}^{(\ell)}}{\left\|u_{N_{H,h},h}^{(\ell)}\right\|_b}$ constitute an orthonormal 
basis for the space $V_{H,h}^{(\ell)}$
in the sense of the inner product $b_h(\cdot, \cdot)$.

Combining (\ref{Upper_Bound_Result}) and (\ref{Equality_4_i_aug}) leads to the following inequality
\begin{eqnarray}\label{Equality_5_k_aug}
\left\|(I-F_{k,h}^{(\ell)})\mathcal P_{H,h}^{(\ell)}\bar u_{i,h}\right\|_{a,h}^2
&\leq&\frac{\bar\mu_{k+1,h}}{(\delta_{k,i,h}^{(\ell)})^2}\eta_a^2(W_H)\left\|(I-\mathcal P_{H,h}^{(\ell)})\bar u_{i,h}\right\|_{a,h}^2.
\end{eqnarray}
From (\ref{Equality_5_k_aug}) and the orthogonal property
$a_h((I-\mathcal P_{H,h}^{(\ell)})\bar u_{i,h}, (I-F_{k,h}^{(\ell)})\mathcal P_{H,h}^{(\ell)}\bar u_{i,h})=0$,
it follows that
\begin{eqnarray*}
\left\|\bar u_{i,h}-F_{k,h}^{(\ell)}\bar u_{i,h}\right\|_{a,h}^2&=&\left\|\bar u_{i,h}-\mathcal P_{H,h}^{(\ell)}\bar u_{i,h}\right\|_{a,h}^2
+\left\|(I-F_{k,h}^{(\ell)})\mathcal P_{H,h}^{(\ell)}\bar u_{i,h}\right\|_{a,h}^2\nonumber\\
&\leq&\left(1+\frac{\bar\mu_{k+1,h}}{(\delta_{k,i,h}^{(\ell)})^2}\eta_a^2(W_H)\right)
\left\|(I-\mathcal P_{H,h}^{(\ell)})\bar u_{i,h}\right\|_{a,h}^2.
\end{eqnarray*}
This is the desired result (\ref{Error_Estimate_Inverse_aug}).
%-----------------------------------------------------------------------------------------------------
%---------------

Similarly, with the help of (\ref{Upper_Bound_Result}), (\ref{Orthogonal_Decomposition_k_aug}),
(\ref{Alpha_Estimate_aug}), (\ref{Equality_u_j_aug}) and (\ref{Orthonormal_Basis_aug}),
we have the following estimates
\begin{eqnarray*}
&&\left\|(I-F_{k,h}^{(\ell)})\mathcal P_{H,h}^{(\ell)}\bar u_{i,h}\right\|_b^2 = \left\|\sum_{j=k+1}^{N_{H,h}}\alpha_ju_{j,h}^{(\ell)}\right\|_b^2
= \sum_{j=k+1}^{N_{H,h}}\alpha_j^2\left\|u_{j,h}^{(\ell)}\right\|_b^2\nonumber\\
&&=\sum_{j=k+1}^{N_{H,h}} \left(\frac{1}{\bar \mu_{i,h}-\mu_{j,h}^{(\ell)}}\right)^2
b_h\big(\bar u_{i,h}-\mathcal P_{H,h}^{(\ell)}\bar u_{i,h},u_{j,h}^{(\ell)}\big)^2\left\|u_{j,h}^{(\ell)}\right\|_{b,h}^2\nonumber\\
&&\leq\frac{1}{(\delta_{k,i,h}^{(\ell)})^2}\sum_{j=k+1}^{N_{H,h}}\left\|u_{j,h}^{(\ell)}\right\|_{b,h}^4\
b_h\left(\bar u_{i,h}-\mathcal P_{H,h}^{(\ell)}\bar u_{i,h},\frac{u_{j,h}^{(\ell)}}{\left\|u_{j,h}^{(\ell)}\right\|_{b,h}}\right)^2\nonumber\\
&&=\frac{1}{(\delta_{k,i,h}^{(\ell)})^2}\sum_{j=k+1}^{N_{H,h}}(\mu_{j,h}^{(\ell)})^2
b_h\left(\bar u_{i,h}-\mathcal P_{H,h}^{(\ell)}\bar u_{i,h},\frac{u_{j,h}^{(\ell)}}{\left\|u_{j,h}^{(\ell)}\right\|_{b,h}}\right)^2
\leq \frac{(\mu_{k+1,h}^{(\ell)})^2}{(\delta_{k,i,h}^{(\ell)})^2}\left\|\bar u_{i,h}-\mathcal P_{H,h}^{(\ell)}\bar u_{i,h}\right\|_{b,h}^2\nonumber\\
&&
\leq \frac{\bar \mu_{k+1,h}^2}{(\delta_{k,i,h}^{(\ell)})^2}\left\|\bar u_{i,h}-\mathcal P_{H,h}^{(\ell)}\bar u_{i,h}\right\|_{b,h}^2,
\end{eqnarray*}
which leads to the inequality
\begin{eqnarray}\label{Equality_8_k_aug}
\left\|(I-F_{k,h}^{(\ell)})\mathcal P_{H,h}^{(\ell)}\bar u_{i,h}\right\|_{b,h} \leq \frac{\bar\mu_{k+1,h}}{\delta_{k,i,h}^{(\ell)}}
\left\|\bar u_{i,h}-\mathcal P_{H,h}^{(\ell)}\bar u_{i,h}\right\|_{b,h}.
\end{eqnarray}
From (\ref{L2_Energy_Estimate}), (\ref{Equality_8_k_aug}) and the triangle inequality, we have  the
following error estimates for the eigenvector approximations in the $\left\|\cdot\right\|_{b,h}$-norm
\begin{eqnarray*}\label{Inequality_11}
&&\left\|\bar u_{i,h}-F_{k,h}^{(\ell)}\bar u_{i,h}\right\|_{b,h}\leq
\left\|\bar u_{i,h}-\mathcal P_{H,h}^{(\ell)}\bar u_{i,h}\right\|_{b,h}
+ \left\|(I-F_{k,h}^{(\ell)})\mathcal P_{H,h}^{(\ell)}\bar u_{i,h}\right\|_{b,h}\nonumber\\
&&\leq\left(1+\frac{\bar\mu_{k+1,h}}{\delta_{k,i,h}^{(\ell)}}\right)
\left\|(I-\mathcal P_{H,h}^{(\ell)})\bar u_{i,h}\right\|_{b,h}\leq \left(1+\frac{\bar\mu_{k+1,h}}{\delta_{k,i,h}^{(\ell)}}\right)
\eta_a(W_H)\left\|(I-\mathcal P_{H,h}^{(\ell)})\bar u_{i,h}\right\|_{a,h}\nonumber\\
&&\leq \left(1+\frac{\bar\mu_{k+1,h}}{\delta_{k,i,h}^{(\ell)}}\right)\eta_a(W_H)
\left\|\bar u_{i,h}-F_{k,h}^{(\ell)}\bar u_{i,h}\right\|_{a,h}.
\end{eqnarray*}
This is the second desired result (\ref{L2_Error_Estimate_Algorithm_1_aug}) and the proof is completed.
\end{proof}

\begin{theorem}\label{Theorem_Error_Estimate_k}
Under the conditions of Lemma \ref{Lemma_Error_Estimate_Subspace_k},
Algorithm \ref{Algorithm_k} has the following error estimate for $\ell \geq 1$
\begin{eqnarray}\label{Error_Estimate_Inverse}
\left\|\bar u_{i,h} -F_{k,h}^{(\ell+1)}\bar u_{i,h} \right\|_{a,h} \leq \gamma
\left\|\bar u_{i,h} - F_{k,h}^{(\ell)}\bar u_{i,h} \right\|_{a,h},
\end{eqnarray}
where
\begin{eqnarray}\label{Definition_Gamma}
\gamma = \bar\lambda_{i,h} \sqrt{1+\frac{\eta_a^2(W_H)}{\bar\lambda_{k+1,h}\big(\delta_{k,i,h}^{(\ell+1)}\big)^2}}
\left(1+\frac{\bar\mu_{k+1,h}}{\delta_{k,i,h}^{(\ell)}}\right)\eta_a^2(W_H).
\end{eqnarray}
\end{theorem}
%--------------------------------------------------------------------------------------------------
%------------
\begin{proof}
%First, let us consider the error estimate $\left\|\bar u_{i,h}- F_{k,h}^{(\ell)}\bar u_{i,h}\right\|_{b,h}$.
%Due to Algorithm \ref{Algorithm_k},
%we know that the approximations $u_{1,h}^{(\ell)}, \cdots, u_{k,h}^{(\ell)}$ come
%from (\ref{Aug_Eigenvalue_Problem_k_1}) (the case that $\ell=1$) or (\ref{Aug_Eigenvalue_Problem_k})
%(the case that $\ell >1$). For both cases, from Lemma  \ref{Lemma_Error_Estimate_Subspace_k},
%there exist exact eigenfunctions $\bar u_{1,h},\cdots, \bar u_{k,h}$
%such that the following error estimates for the eigenvector approximations
%$u_{1,h}^{(\ell)},\cdots,  u_{k,h}^{(\ell)}$ hold for $i=1, \cdots, k$
%\begin{eqnarray}\label{Inequality_13}
%\left\|\bar u_{i,h}- F_{k,h}^{(\ell)}\bar u_{i,h}\right\|_{b,h} &\leq \left(1+\frac{\mu_{k+1}}{\delta_{k,i}^{(\ell)}}\right)
%\eta_a(V_{H,h}^{(\ell)})\left\|\bar u_{i,h}- F_{k,h}^{(\ell)}\bar u_{i,h}\right\|_{a,h} \nonumber\\
%&\leq \left(1+\frac{\mu_{k+1}}{\delta_{k,i}^{(\ell)}}\right)\eta_a(V_H)\left\|\bar u_{i,h}- F_{k,h}^{(\ell)}\bar u_{i,h}\right\|_{a,h},
%\end{eqnarray}
%where we have used the inequality $\eta_a(V_{H,h}^{(\ell)})\leq \eta_a(V_H)$ since $V_H\subset V_{H,h}^{(\ell)}$.

From Algorithm \ref{Algorithm_k}, it is easy to know that $u_{1,h}^{(\ell)}, \cdots, u_{k,h}^{(\ell)}$
is the orthogonal basis for the space
${\rm span}\{u_{1,h}^{(\ell)}, \cdots, u_{k,h}^{(\ell)}\}$. We define the $b_h(\cdot,\cdot)$-orthogonal
projection operator $\pi_{k,h}^{(\ell)}$ to the space ${\rm span}\{u_{1,h}^{(\ell)}$, $\cdots$, $u_{k,h}^{(\ell)}\}$.
Then there exist $k$ real numbers $q_1, \cdots, q_k \in \mathbb R$ such that $\pi_{k,h}^{(\ell)}\bar u_{i,h}$
has the following expansion
\begin{eqnarray}\label{Expansion_L2_aug}
\pi_{k,h}^{(\ell)}\bar u_{i,h} = \sum_{j=1}^k q_ju_{j,h}^{(\ell)}.
\end{eqnarray}
From (\ref{Projection_aug}) and the definition of $V_{H,h}^{(\ell+1)}$ in Step 3 of Algorithm \ref{Algorithm_k},
we obtain the orthogonal property of the projection operator $\mathcal P_{H,h}^{(\ell+1)}$, together with
(\ref{Linear_Equation_k}), (\ref{L2_Energy_Estimate}), (\ref{L2_Error_Estimate_Algorithm_1_aug})
and (\ref{Expansion_L2_aug}),
the following inequalities hold
\begin{eqnarray}\label{Inequality_16}
&&\left\|\bar u_{i,h} - \mathcal P_{H,h}^{(\ell+1)}\bar u_{i,h}\right\|_{a,h}^2 =
a_h\left(\bar u_{i,h} - \mathcal P_{H,h}^{(\ell+1)}\bar u_{i,h}, \bar u_{i,h} - \mathcal P_{H,h}^{(\ell+1)}\bar u_{i,h}\right)\nonumber\\
&&=a_h\left(\bar u_{i,h} - \sum_{j=1}^k\bar\lambda_{i,h}\frac{q_j}{\lambda_{j,h}^{(\ell)}}
\widehat u_{j,h}^{(\ell+1)}, \bar u_{i,h} - \mathcal P_{H,h}^{(\ell+1)}\bar u_{i,h}\right)\nonumber\\
&&=\bar\lambda_{i,h} b_h\left(\bar u_{i,h} - \sum_{j=1}^k\frac{q_j}{\lambda_{j,h}^{(\ell)}}
\lambda_{j,h}^{(\ell)}u_{j,h}^{(\ell)}, \bar u_{i,h} - \mathcal P_{H,h}^{(\ell+1)}\bar u_{i,h}\right)\nonumber\\
&&=\bar\lambda_{i,h} b_h\left(\bar u_{i,h} - \sum_{j=1}^kq_ju_{j,h}^{(\ell)},
\bar u_{i,h} - \mathcal P_{H,h}^{(\ell+1)}\bar u_{i,h}\right)
=\bar\lambda_{i,h} b_h\left(\bar u_{i,h} - \pi_{k,h}^{(\ell)}\bar u_{i,h}, \bar u_{i,h}
- \mathcal P_{H,h}^{(\ell+1)}\bar u_{i,h}\right)\nonumber\\
&&\leq \bar\lambda_{i,h}\left\|\bar u_{i,h} - \pi_{k,h}^{(\ell)}\bar u_{i,h}\right\|_{b,h}\left\|\bar u_{i,h}
- \mathcal P_{H,h}^{(\ell+1)}\bar u_{i,h}\right\|_{b,h}\nonumber\\
&&\leq \bar\lambda_{i,h}\left\|\bar u_{i,h} - F_{k,h}^{(\ell)}\bar u_{i,h}\right\|_{b,h}\left\|\bar u_{i,h}
- \mathcal P_{H,h}^{(\ell+1)}\bar u_{i,h}\right\|_{b,h}\nonumber\\
&&\leq  \bar\lambda_{i,h} \bar\eta_a(W_{H})\left\|\bar u_{i,h}- F_{k,h}^{(\ell)}\bar u_{i,h}\right\|_{a,h}
\eta_a(W_H)\left\|\bar u_{i,h} - \mathcal P_{H,h}^{(\ell+1)}\bar u_{i,h}\right\|_{a,h},
\end{eqnarray}
where $\bar\eta_a(W_{H})$ is defined in Lemma \ref{Lemma_Error_Estimate_Subspace_k}.
%Corollary \ref{Error_Estimate_Corollary_k}.
Then from (\ref{Inequality_16}), it follows that
\begin{eqnarray}\label{Inequality_18_k_aug}
\left\|\bar u_{i,h} - \mathcal P_{H,h}^{(\ell+1)}\bar u_{i,h}\right\|_{a,h} \leq \bar\lambda_{i,h}
\bar\eta_a(W_{H})
\eta_a(W_H)\left\|\bar u_{i,h}- F_{k,h}^{(\ell)}\bar u_{i,h}\right\|_{a,h}.
\end{eqnarray}
Since $u_{1,h}^{(\ell+1)}, \cdots, u_{k,h}^{(\ell+1)}$ only come from (\ref{Aug_Eigenvalue_Problem_k}),  and
Lemma \ref{Lemma_Error_Estimate_Subspace_k}, we have for $i=1, \cdots, k$
\begin{eqnarray*}
\left\|\bar u_{i,h}- F_{k,h}^{(\ell+1)}\bar u_{i,h}\right\|_{a,h}
%&\leq& \sqrt{1+\frac{\eta_a^2(V_{H,h}^{(\ell+1)})}{\lambda_{k+1}\big(\delta_{k,i}^{(\ell+1)}\big)^2}}
%\left\|(I-\mathcal P_{H,h}^{(\ell+1)})\bar u_{i,h}\right\|_{a,h} \\
&\leq&   \sqrt{1+\frac{\eta_a^2(W_H)}{\bar\lambda_{k+1,h}\big(\delta_{k,i,h}^{(\ell+1)}\big)^2}}
\left\|(I-\mathcal P_{H,h}^{(\ell+1)})\bar u_{i,h}\right\|_{a,h}.
\end{eqnarray*}
Together with (\ref{Inequality_18_k_aug}), we arrive at
\begin{eqnarray*}\label{Inequality_19}
\left\|\bar u_{i,h} -F_{k,h}^{(\ell+1)}\bar u_{i,h}\right\|_{a,h} \leq \bar\lambda_{i,h} \sqrt{1+\frac{\eta_a^2(W_H)}
{\bar\lambda_{k+1,h}
\big(\delta_{k,i,h}^{(\ell+1)}\big)^2}} \bar\eta_a(W_{H})
\eta_a(W_H)\left\|\bar u_{i,h}- F_{k,h}^{(\ell)}\bar u_{i,h}\right\|_{a,h},
\end{eqnarray*}
which is the desired result (\ref{Error_Estimate_Inverse}) and the proof is completed.
\end{proof}
%---------------------------------------------------------------------------------------------
\begin{remark}
According to Theorem \ref{Theorem_Error_Estimate_k}, the augmented subspace techniques 
have a second order convergence rate, as indicated by the convergence 
result (\ref{Error_Estimate_Inverse}). 
Furthermore, we ought to lower the term $\eta_a(W_H)$, which is dependent on 
the coarse conforming linear finite element space $W_H$, 
in order to speed up the convergence rate.
In other words, the convergence can be accelerated by expanding the subspace $W_H$.
\end{remark}
%---------------------------------------------------------------------------------------------
\begin{remark}\label{Remark_Eigenvalue}
Since the error estimates for the eigenvalue approximation can be simply inferred from the following error expansion, we are only concerned with the error estimates for the eigenvector approximation in this paper
\begin{eqnarray*}\label{rayexpan}
0\leq \widehat{\lambda}_i-\bar\lambda_{i,h}
=\frac{a_h(\bar u_{i,h}-\psi,\bar u_{i,h}-\psi)}{b(\psi,\psi)}-\bar\lambda_{i,h}
\frac{b(\bar u_{i,h}-\psi,\bar u_{i,h}-\psi)}{b(\psi,\psi)}+2\frac{a_h(\bar u_{i,h},\psi)-b(\bar \lambda_{i,h}\bar u_{i,h}, \psi)}{b(\psi, \psi)},
\end{eqnarray*}
where $\psi$ is the eigenfunction approximation for the exact eigenfunction $\bar u_{i,h}$ and
\begin{eqnarray*}
\widehat{\lambda}_i=\frac{a_h(\psi,\psi)}{b(\psi,\psi)}.
\end{eqnarray*}
\end{remark}
%---------------------------------------------------------------------------------------------
Since each linear equation can be solved separately, it follows that Step 2 
of Algorithm \ref{Algorithm_k} can be performed using the parallel 
computing approach. Nevertheless, a kind of parallel methods for eigenvalue 
problems can be designed using the augmented subspace approach. The eigenvalue problem (\ref{Aug_Eigenvalue_Problem_k}) is solved in Step 3 of Algorithm \ref{Algorithm_k}. 
However, we must perform the inner products of the $k$ vectors 
in the high dimensional space $V_h$ in order to generate 
the matrices for (\ref{Aug_Eigenvalue_Problem_k}). 
This is a very low scalable procedure for the parallel 
computing \cite{LiXieXuYouZhang,XuXieZhang, ZhangLiXieXuYou}.
That is to say, a bottleneck for parallel computing does exist in the inner 
product calculation for many high dimensional vectors. We provide an additional 
version of the augmented subspace technique for a single (possibly non-smallest) 
eigenpair, which represents the single process version of this kind of parallel schemes, 
to get around this crucial bottleneck.
Algorithm \ref{Algorithm_1} defines the relevant numerical approach.
This idea in relation to the conforming finite element technique has already 
been put out and examined in \cite{XuXieZhang}.

In Algorithm \ref{Algorithm_1},  we assume that the given eigenpair
approximation $(\lambda_{i,h}^{(\ell)}, u_{i,h}^{(\ell)})\in\mathbb R\times V_h$
with different superscripts is the closest
to an exact eigenpair $(\bar\lambda_{i,h}, \bar u_{i,h})$ of (\ref{Weak_Eigenvalue_Discrete}).
Based on this setting, we can give the following convergence result for the augmented
subspace method defined by Algorithm \ref{Algorithm_1}.
\begin{algorithm}[hbt!]
\caption{Augmented subspace method for one eigenpair}\label{Algorithm_1}
\begin{enumerate}
\item For $\ell=1$, we define $\widehat u_{i,h}^{(\ell)}=u_{i,h}^{(\ell)}$, and
the augmented subspace $V_{H,h}^{(\ell)} = W_H +{\rm span}\{\widehat u_{i,h}^{(\ell)}\}$.
Then solve the following eigenvalue problem:
Find $(\lambda_{i,h}^{(\ell)},u_{i,h}^{(\ell)})\in \mathbb R\times V_{H,h}^{(\ell)}$
such that $a_h(u_{i,h}^{(\ell)},u_{i,h}^{(\ell)})=1$ and
\begin{equation}\label{Aug_Eigenvalue_Problem_1_1}
a_h(u_{i,h}^{(\ell)},v_{H,h}) = \lambda_{i,h}^{(\ell)}b_h(u_{i,h}^{(\ell)},v_{H,h}),
\ \ \ \ \ \forall v_{H,h}\in V_{H,h}^{(\ell)}.
\end{equation}

\item Solve the following linear boundary value problem:
Find $\widehat{u}_{i,h}^{(\ell+1)}\in V_h$ such that
\begin{equation}\label{Linear_Equation}
a_h(\widehat{u}_{i,h}^{(\ell+1)},v_h) = \lambda_{i,h}^{(\ell)}b_h(u_{i,h}^{(\ell)},v_h),
\ \  \forall v_h\in V_h.
\end{equation}
\item Define the augmented subspace $V_{H,h}^{(\ell+1)} = W_H +
{\rm span}\{\widehat{u}_{i,h}^{(\ell+1)}\}$ and solve the following eigenvalue problem:
Find $(\lambda_{i,h}^{(\ell+1)},u_{i,h}^{(\ell+1)})\in \mathbb R\times V_{H,h}^{(\ell+1)}$
such that $a_h(u_{i,h}^{(\ell+1)},u_{i,h}^{(\ell+1)})=1$ and
\begin{equation}\label{Aug_Eigenvalue_Problem_1_2}
a_h(u_{i,h}^{(\ell+1)},v_{H,h}) = \lambda_{i,h}^{(\ell+1)}b_h(u_{i,h}^{(\ell+1)},v_{H,h}),
\ \ \ \ \ \forall v_{H,h}\in V_{H,h}^{(\ell+1)}.
\end{equation}
Solve (\ref{Aug_Eigenvalue_Problem_1_2}) and the output $(\lambda_{i,h}^{(\ell+1)},u_{i,h}^{(\ell+1)})$
is chosen such that $u_{i,h}^{(\ell+1)}$ has the largest component in ${\rm span}\{\widehat{u}_{i,h}^{(\ell+1)}\}$
among all eigenfunctions of (\ref{Aug_Eigenvalue_Problem_1_2}).
\item Set $\ell=\ell+1$ and go to Step 2 for the next iteration until convergence.
\end{enumerate}
\end{algorithm}

For each $\ell$, it is easy to know, the eigenvalue problem (\ref{Aug_Eigenvalue_Problem_1_1})
and  (\ref{Aug_Eigenvalue_Problem_1_2}) also have the following eigenvalues  \cite{BabuskaOsborn_1989,BabuskaOsborn_Book},
$$0< \lambda_{1,h}^{(\ell)}\leq  \lambda_{2,h}^{(\ell)}
\leq \cdots \leq \lambda_{k,h}^{(\ell)} \leq \cdots \leq  \lambda_{N_{H,h},h}^{(\ell)},$$
and corresponding eigenfunctions
\begin{eqnarray}\label{Discrete_Eigenfunctions_aug}
u_{1,h}^{(\ell)},  u_{2,h}^{(\ell)}, \cdots,  u_{k,h}^{(\ell)}, \cdots,  u_{N_{H,h},h}^{(\ell)},
\end{eqnarray}
where $N_{H,h}={\rm dim}V_{H,h}^{(\ell)} = N_H+1$ and $a_h(u_{i,h}^{(\ell)}, u_{j,h}^{(\ell)})=\delta_{ij}$, $1\leq i,j\leq N_{H,h}$.

It is simple to understand that the WG finite element space $V_h$ is a subspace of the low dimensional augmented subspace $V_{H,h}^{(\ell)}$ in Algorithm \ref{Algorithm_1}.
Then, Algorithm \ref{Algorithm_1}'s error estimates are comparable to those of Algorithm \ref{Algorithm_k}. We also utilize the definitions (\ref{Projection_aug}) and (\ref{eta_a_H_Def}) for the sake of simplicity in notation. Next, we apply the property (\ref{Upper_Bound_Result}) and error estimates (\ref{Energy_Error_Estimate}), and finally, we employ (\ref{L2_Energy_Estimate}) for the eigenvalue problems (\ref{Aug_Eigenvalue_Problem_1_1}) and (\ref{Aug_Eigenvalue_Problem_1_2}).

\begin{lemma}\label{Lemma_Error_Estimate_Subspace_1}
Let  $(\bar\lambda_h,\bar u_h)$ denote an exact eigenpair of
the eigenvalue problem (\ref{Weak_Eigenvalue_Discrete}).
Assume the eigenpair approximation $(\lambda_{i,h}^{(\ell)}, u_{i,h}^{(\ell)})$ 
has the property that
$\mu_{i,h}^{(\ell)}=1/\lambda_{i,h}^{(\ell)}$ is closest to $\bar\mu_h=1/{\bar \lambda_h}$.
The spectral projector $E_{i,h}^{(\ell)}: V_h\mapsto {\rm span}\{u_{i,h}^{(\ell)}\}$
according to the eigenpair approximation $(\lambda_{i,h}^{(\ell)},u_{i,h}^{(\ell)})
\in\mathbb R\times V_{H,h}^{(\ell)}$ is defined as follows
\begin{eqnarray*}
a_h(E_{i,h}^{(\ell)}w, u_{i,h}^{(\ell)}) = a_h(w, u_{i,h}^{(\ell)}),\ \  \ \ {\rm for}\  w\in V_h.
\end{eqnarray*}
Then the eigenpair approximation $(\lambda_{i,h}^{(\ell)},u_{i,h}^{(\ell)})
\in\mathbb R\times V_{H,h}^{(\ell)}$ produced by
Algorithm \ref{Algorithm_1} satisfies the following error estimates
\begin{eqnarray}
\left\|\bar u_h-E_{i,h}^{(\ell)}\bar u_h\right\|_{a,h} &\leq&\bar\lambda_{i,h}
\sqrt{1+\frac{\eta_a^2(W_H)}{\bar\lambda_{1,h}\big(\delta_{\lambda,h}^{(\ell)}\big)^2}}
\left\|\bar u_h- \mathcal P_{H,h}^{(\ell)} \bar u_h\right\|_{a,h},\ \ \ \ \ \label{Estimate_h_1_a_aug}\\
\left\|\bar u_h-E_{i,h}^{(\ell)}\bar u_h\right\|_{b,h}&\leq&\bar\eta_a(W_H)
\left\|\bar u_h-E_{i,h}^{(\ell)}\bar u_h\right\|_{a,h},\label{Estimate_h_1_b_aug}
\end{eqnarray}
where $\delta_{\lambda,h}$ and $\bar\eta_a(W_H)$ are defined as follows
\begin{eqnarray}\label{Definition_Delta_k_1_aug}
\delta_{\lambda,h}^{(\ell)} = \min_{j\neq i}\left|\frac{1}{\lambda_{j,h}^{(\ell)}}-\frac{1}{\bar\lambda_h}\right|,
\ \ \
\bar\eta_a(W_{H}) = \left(1+\frac{1}{\bar\lambda_{1,h}\delta_{\lambda,h}^{(\ell)}}\right)\eta_a(W_H).
\end{eqnarray}
\end{lemma}
\begin{proof}
Since $(I-E_{i,h}^{(\ell)})\mathcal P_{H,h}^{(\ell)}\bar u_h\in V_{H,h}^{(\ell)}$ and
$(I-E_{i,h}^{(\ell)})\mathcal P_{H,h}^{(\ell)}\bar u_h\in {\rm span}\{u_{1,h}^{(\ell)},\cdots, u_{i-1,h}^{(\ell)}, u_{i+1,h}^{(\ell)}$,
$\cdots$,  $u_{N_{H,h},h}^{(\ell)}\}$,
the following orthogonal expansion holds
\begin{eqnarray}\label{Orthogonal_Decomposition_1_aug}
(I-E_{i,h}^{(\ell)})\mathcal P_{H,h}^{(\ell)}\bar u_h=\sum_{j\neq i}\alpha_ju_{j,h}^{(\ell)},
\end{eqnarray}
where  $\alpha_j=a_h(\mathcal P_{H,h}^{(\ell)}\bar u_h,u_{j,h}^{(\ell)})$. From Lemma \ref{Strang_Lemma}, we have
the same equality (\ref{Alpha_Estimate_aug}).

From the orthogonal property of eigenfunctions $u_{1,h}^{(\ell)},\cdots, \bar u_{N_{H,h},h}^{(\ell)}$, we acquire
\begin{eqnarray*}
1 = a_h(u_{j,h}^{(\ell)},u_{j,h}^{(\ell)}) = \lambda_{j,h}^{(\ell)}b_h(u_{j,h}^{(\ell)},u_{j,h}^{(\ell)})
= \lambda_{j,h}^{(\ell)}\left\|u_{j,h}^{(\ell)}\right\|_{b,h}^2,
\end{eqnarray*}
which leads to the following property
\begin{eqnarray}\label{Equality_u_1_aug}
\left\|u_{j,h}^{(\ell)}\right\|_{b,h}^2=\frac{1}{\lambda_{j,h}^{(\ell)}}= \mu_{j,h}^{(\ell)}.
\end{eqnarray}
Because of (\ref{Aug_Eigenvalue_Problem_1_1}),  (\ref{Aug_Eigenvalue_Problem_1_2}) 
and the definition of eigenfunctions
$u_{1,h}^{(\ell)},\cdots, u_{N_{H,h},h}^{(\ell)}$, we obtain the following equalities
\begin{eqnarray}\label{Orthonormal_Basis_aug_1}
a_h(u_{j,h}^{(\ell)}, u_{k,h}^{(\ell)})=\delta_{jk},
\ \ \ \ \ b_h\left(\frac{u_{j,h}^{(\ell)}}{\left\|u_{j,h}^{(\ell)}\right\|_{b,h}},
\frac{u_{k,h}^{(\ell)}}{\left\|u_{k,h}^{(\ell)}\right\|_{b,h}}\right)=\delta_{jk},\ \ \ 1\leq j,k\leq N_{H,h}.
\end{eqnarray}
Then due to  (\ref{Upper_Bound_Result}), (\ref{Alpha_Estimate_aug}),  (\ref{Orthonormal_Basis_aug}),
(\ref{Orthogonal_Decomposition_1_aug}) and  (\ref{Equality_u_1_aug}), we have following estimates
\begin{eqnarray}\label{Equality_4_1_aug}
&&\left\|(I-E_{i,h}^{(\ell)})\mathcal P_{H,h}^{(\ell)}\bar u_h\right\|_{a,h}^2
= \left\|\sum_{j\neq i}\alpha_ju_{j,h}^{(\ell)}\right\|_{a,h}^2
= \sum_{j\neq i}\alpha_j^2
=\sum_{j\neq i} \left(\frac{1}{\bar \mu_h- \mu_{j,h}^{(\ell)}}\right)^2 b_h\big(\bar u_h-\mathcal P_{H,h}^{(\ell)}\bar u_h,u_{j,h}^{(\ell)}\big)^2\nonumber\\
&&
\leq\frac{1}{(\delta_{\lambda,h}^{(\ell)})^2}\sum_{j\neq i}\left\|u_{j,h}^{(\ell)}\right\|_{b,h}^2
b_h\left(\bar u_h-\mathcal P_{H,h}^{(\ell)}\bar u_h,\frac{u_{j,h}^{(\ell)}}{\left\|u_{j,h}^{(\ell)}\right\|_{b,h}}\right)^2\nonumber\\
&&=\frac{1}{(\delta_{\lambda,h}^{(\ell)})^2}\sum_{j\neq i}\mu_{j,h}^{(\ell)}
b_h\left(\bar u_h-\mathcal P_{H,h}^{(\ell)}\bar u_h,\frac{u_{j,h}^{(\ell)}}{\left\|u_{j,h}^{(\ell)}\right\|_b}\right)^2\nonumber\\
&&
\leq \frac{\mu_{1,h}^{(\ell)}}{(\delta_{\lambda,h}^{(\ell)})^2}\sum_{j\neq i}
b_h\left(\bar u_h-\mathcal P_{H,h}^{(\ell)}\bar u_h,\frac{u_{j,h}^{(\ell)}}{\left\|u_{j,h}^{(\ell)}\right\|_{b,h}}\right)^2
\leq \frac{\mu_{1,h}^{(\ell)}}{(\delta_{\lambda,h}^{(\ell)})^2}\left\|\bar u_h-\mathcal P_{H,h}^{(\ell)}\bar u_h\right\|_{b,h}^2,
\end{eqnarray}
where the last inequality holds since $\frac{u_{1,h}^{(\ell)}}{\left\|u_{1,h}^{(\ell)}\right\|_b}$, $\cdots$,
$\frac{u_{N_{H,h},h}^{(\ell)}}{\left\|u_{N_{H,h},h}^{(\ell)}\right\|_b}$ 
constitute an orthonormal  basis for the space $V_{H,h}^{(\ell)}$
in the sense of the inner product $b_h(\cdot, \cdot)$.

Combining (\ref{Upper_Bound_Result}) and (\ref{Equality_4_1_aug}) leads to the following inequality
\begin{eqnarray}\label{Equality_5_1_aug}
\left\|(I-E_{i,h}^{(\ell)})\mathcal P_{H,h}^{(\ell)}\bar u_h\right\|_{a,h}^2
&\leq&\frac{\bar\mu_{1,h}}{(\delta_{\lambda,h}^{(\ell)})^2}\eta_a^2(W_H)\left\|(I-\mathcal P_{H,h}^{(\ell)})\bar u_h\right\|_{a,h}^2.
\end{eqnarray}
From (\ref{Equality_5_1_aug}) and the orthogonal property
$a_h((I-\mathcal P_{H,h}^{(\ell)})\bar u_h, (I-E_{i,h}^{(\ell)})\mathcal P_{H,h}^{(\ell)}\bar u_h)=0$,
it follows that
\begin{eqnarray*}
\left\|\bar u_h-E_{i,h}^{(\ell)}\bar u_h\right\|_{a,h}^2&=&\left\|\bar u_h-\mathcal P_{H,h}^{(\ell)}\bar u_h\right\|_{a,h}^2
+\left\|(I-E_{i,h}^{(\ell)})\mathcal P_{H,h}^{(\ell)}\bar u_h\right\|_{a,h}^2\nonumber\\
&\leq&\left(1+\frac{\bar\mu_{1,h}}{(\delta_{\lambda,h}^{(\ell)})^2}\eta_a^2(W_H)\right)
\left\|(I-\mathcal P_{H,h}^{(\ell)})\bar u_h\right\|_{a,h}^2.
\end{eqnarray*}
This is the desired result (\ref{Estimate_h_1_a_aug}).
%-----------------------------------------------------------------------------------------------------

Similarly, with the help of (\ref{Upper_Bound_Result}), (\ref{Alpha_Estimate_aug}),
(\ref{Orthogonal_Decomposition_1_aug}), (\ref{Equality_u_1_aug}) and (\ref{Orthonormal_Basis_aug_1}),
we have following estimates
\begin{eqnarray*}
&&\left\|(I-E_{i,h}^{(\ell)})\mathcal P_{H,h}^{(\ell)}\bar u_h\right\|_{b,h}^2
= \left\|\sum_{j\neq i}\alpha_ju_{j,h}^{(\ell)}\right\|_{b,h}^2
= \sum_{j\neq i}\alpha_j^2\left\|u_{j,h}^{(\ell)}\right\|_{b,h}^2\nonumber\\
&&=\sum_{j\neq i} \left(\frac{1}{\bar \mu_h-\mu_{j,h}^{(\ell)}}\right)^2
b_h\big(\bar u_h-\mathcal P_{H,h}^{(\ell)}\bar u_h, u_{j,h}^{(\ell)}\big)^2\left\|u_{j,h}^{(\ell)}\right\|_{b,h}^2\nonumber\\
&&\leq\frac{1}{(\delta_{\lambda,h}^{(\ell)})^2}\sum_{j\neq i}\left\|u_{j,h}^{(\ell)}\right\|_{b,h}^4\
b_h\left(\bar u_h-\mathcal P_{H,h}^{(\ell)}\bar u_h, \frac{u_{j,h}^{(\ell)}}{\left\|u_{j,h}^{(\ell)}\right\|_{b,h}}\right)^2\nonumber\\
&&=\frac{1}{(\delta_{\lambda,h}^{(\ell)})^2}\sum_{j\neq i}(\mu_{j,h}^{(\ell)})^2
b_h\left(\bar u_h-\mathcal P_{H,h}^{(\ell)}\bar u_h,\frac{u_{j,h}^{(\ell)}}{\left\|u_{j,h}^{(\ell)}\right\|_{b,h}}\right)^2
\leq \frac{(\mu_{1,h}^{(\ell)})^2}{(\delta_{\lambda,h}^{(\ell)})^2}\left\|\bar u_h-\mathcal P_{H,h}^{(\ell)}
\bar u_h\right\|_{b,h}^2\nonumber\\
&&
\leq \frac{\bar \mu_{1,h}^2}{(\delta_{\lambda,h}^{(\ell)})^2}\left\|\bar u_h-\mathcal P_{H,h}^{(\ell)}
\bar u_h\right\|_{b,h}^2,
\end{eqnarray*}
which leads to the inequality
\begin{eqnarray}\label{Equality_8_1_aug}
\left\|(I-E_{i,h}^{(\ell)})\mathcal P_{H,h}^{(\ell)}\bar u_h\right\|_{b,h} \leq \frac{\bar\mu_{1,h}}{\delta_{\lambda,h}^{(\ell)}}
\left\|\bar u_h-\mathcal P_{H,h}^{(\ell)}\bar u_h\right\|_{b,h}.
\end{eqnarray}
From (\ref{L2_Energy_Estimate}), (\ref{Equality_8_1_aug}) and the triangle inequality, we have  the
following error estimates for the eigenvector approximations in the $\left\|\cdot\right\|_{b,h}$-norm
\begin{eqnarray*}
&&\left\|\bar u_h-E_{i,h}^{(\ell)}\bar u_h\right\|_{b,h}\leq
\left\|\bar u_h-\mathcal P_{H,h}^{(\ell)}\bar u_h\right\|_{b,h}
+ \left\|(I-E_{i,h}^{(\ell)})\mathcal P_{H,h}^{(\ell)}\bar u_h\right\|_{b,h}\nonumber\\
&&\leq\left(1+\frac{\bar\mu_{1,h}}{\delta_{\lambda,h}^{(\ell)}}\right)
\left\|(I-\mathcal P_{H,h}^{(\ell)})\bar u_h\right\|_{b,h}
\leq \left(1+\frac{\bar\mu_{1,h}}{\delta_{\lambda,h}^{(\ell)}}\right)
\eta_a(W_H)\left\|(I-\mathcal P_{H,h}^{(\ell)})\bar u_h\right\|_{a,h}\nonumber\\
&&\leq \left(1+\frac{\bar\mu_{1,h}}{\delta_{\lambda,h}^{(\ell)}}\right)\eta_a(W_H)
\left\|\bar u_h-E_{i,h}^{(\ell)}\bar u_h\right\|_{a,h}.
\end{eqnarray*}
This is the second desired result (\ref{Estimate_h_1_b_aug}) and the proof is completed.
\end{proof}

%--------------------------------------------------------------------------------------------
\begin{theorem}\label{Theorem_Error_Estimate_1}
Under the conditions of Lemma \ref{Lemma_Error_Estimate_Subspace_1},
Algorithm \ref{Algorithm_1} has the following error estimate for $\ell\geq 1$
\begin{eqnarray}
\left\|\bar u_h-E_{i,h}^{(\ell+1)}\bar u_h\right\|_{a,h} \leq\bar\lambda_{i,h}
\sqrt{1+\frac{\eta_a^2(W_H)}{\bar \lambda_{1,h}\big(\delta_{\lambda,h}^{(\ell+1)}\big)^2}}
\left(1+\frac{1}{\bar\lambda_{1,h}\delta_{\lambda,h}^{(\ell)}}\right)\eta_a^2(W_H)
\left\|\bar u_h-E_{i,h}^{(\ell)}\bar u_h\right\|_{a,h}.\ \ \ \ \ \label{Estimate_h_1_a}
\end{eqnarray}
\end{theorem}
%----------------------------------------------------------------------------------------------
\begin{proof}
We define the $b(\cdot,\cdot)$-orthogonal projection operator
$\pi_h^{(\ell)}$ to the space ${\rm span}\{u_{i,h}^{(\ell)}\}$.
Then there exists a real number $q\in\mathbb R$
such that $\pi_h^{(\ell)}\bar u_h = q u_{i,h}^{(\ell)}$.
Then from  the orthogonal property of the projection operator $\mathcal P_{H,h}^{(\ell+1)}$,
(\ref{L2_Energy_Estimate}), (\ref{Linear_Equation}) and (\ref{Estimate_h_1_b_aug}), we obtain
\begin{eqnarray}\label{Inequality_16_2}
&&\left\|\bar u_h - \mathcal P_{H,h}^{(\ell+1)}\bar u_h\right\|_{a,h}^2
= a_h\left(\bar u_h - \mathcal P_{H,h}^{(\ell+1)}\bar u_h, \bar u_h 
- \mathcal P_{H,h}^{(\ell+1)}\bar u_h\right)\nonumber\\
&&=a_h\left(\bar u_h -\frac{\bar\lambda_{i,h}}{\lambda_{i,h}^{(\ell)}}q\widehat u_{i,h}^{(\ell+1)}, 
\bar u_h - \mathcal P_{H,h}^{(\ell+1)}\bar u_h\right)\nonumber\\
&&=a_h\left(\bar u_h,  \bar u_h - \mathcal P_{H,h}^{(\ell+1)}\bar u_h\right) - \frac{\bar\lambda_{i,h}}{\lambda_{i,h}^{(\ell)}}q
a_h\left(\widehat u_{i,h}^{(\ell+1)}, \bar u_h - \mathcal P_{H,h}^{(\ell+1)}\bar u_h\right)\nonumber\\
&&=\bar\lambda_hb_h\left(\bar u_h,  \bar u_h - \mathcal P_{H,h}^{(\ell+1)}\bar u_h\right) 
- \bar\lambda_{i,h}
b_h\left(q u_{i,h}^{(\ell)}, \bar u_h - \mathcal P_{H,h}^{(\ell+1)}\bar u_h\right)\nonumber\\
&&=\bar\lambda_hb_h\left(\bar u_h 
- \pi_h^{(\ell)}\bar u_h, \bar u_h - \mathcal P_{H,h}^{(\ell+1)}\bar u_h\right)
\leq \bar\lambda_h\left\|\bar u_h 
- \pi_h^{(\ell)}\bar u_h\right\|_{b,h}\left\|\bar u_h - \mathcal P_{H,h}^{(\ell+1)}\bar u_h\right\|_{b,h}\nonumber\\
&&
\leq \bar\lambda_h\left\|\bar u_h - E_{i,h}^{(\ell)}\bar u_h\right\|_{b,h}\left\|\bar u_h - \mathcal P_{H,h}^{(\ell+1)}\bar u_h\right\|_{b,h}\nonumber\\
&&\leq  \bar\lambda_h\left(1+\frac{1}{\bar\lambda_{1,h}
\delta_{\lambda,h}^{(\ell)}}\right) \eta_a(W_H)\left\|\bar u_h- E_{i,h}^{(\ell)}\bar u_h\right\|_{a,h}
\eta_a(W_H)\left\|\bar u_h - \mathcal P_{H,h}^{(\ell+1)}\bar u_h\right\|_{a,h}\nonumber\\
&&\leq  \bar\lambda_h\left(1+\frac{1}{\bar\lambda_{1,h}\delta_{\lambda,h}^{(\ell)}}\right)
\eta_a^2(W_H)\left\|\bar u_h- E_{i,h}^{(\ell)}\bar u_h\right\|_{a,h}\left\|\bar u_h - \mathcal P_{H,h}^{(\ell+1)}\bar u_h\right\|_{a,h}.
\end{eqnarray}
Since the approximation $u_{i,h}^{(\ell+1)}$ only comes from (\ref{Aug_Eigenvalue_Problem_1_1})
or (\ref{Aug_Eigenvalue_Problem_1_2}),  together with Lemma \ref{Lemma_Error_Estimate_Subspace_1},
we have
\begin{eqnarray}\label{Inequality_17_111}
\left\|\bar u_h- E_{i,h}^{(\ell+1)}\bar u_h\right\|_{a,h}
&\leq&  \sqrt{1+\frac{\eta_a^2(W_H)}{\bar\lambda_{1,h}\big(\delta_{\lambda,h}^{(\ell+1)}\big)^2}}
\left\|(I-\mathcal P_{H,h}^{(\ell+1)})\bar u_h\right\|_{a,h}.
\end{eqnarray}
From (\ref{Inequality_16_2}), there holds
\begin{eqnarray}\label{Inequality_17_2}
\left\|\bar u_h - \mathcal P_{H,h}^{(\ell+1)}\bar u_h\right\|_{a,h} \leq \bar\lambda_h
\left(1+\frac{1}{\bar\lambda_{1,h}\delta_{\lambda,h}^{(\ell)}}\right)\eta_a^2(W_H)\left\|\bar u_h- E_{i,h}^{(\ell)}\bar u_h\right\|_{a,h}.
\end{eqnarray}
Combining (\ref{Inequality_17_111}) with (\ref{Inequality_17_2}), we have the following estimate
\begin{eqnarray*}\label{Inequality_18}
\left\|\bar u_h-E_{i,h}^{(\ell+1)}\bar u_h\right\|_{a,h} \leq \bar\lambda_h
\sqrt{1+\frac{\eta_a^2(W_H)}{\bar\lambda_{1,h}\big(\delta_{\lambda,h}^{(\ell+1)} \big)^2}}
\left(1+\frac{1}{\bar\lambda_{1,h}\delta_{\lambda,h}^{(\ell)}}\right) \eta_a^2(W_H)\left\|\bar u_h-E_{i,h}^{(\ell)}\bar u_h\right\|_{a,h}.
\end{eqnarray*}
This is the desired result (\ref{Estimate_h_1_a}) and the proof is complete.
\end{proof}
%-------------------------------------------------------------------------------------

\section{Applications to Laplace eigenvalue problem}\label{Section_4}
This section will demonstrate the applications of  
the augmented subspace techniques introduced in Section \ref{Section_3} 
to the Laplace eigenvalue problem and provide the associated convergence rates.
It is noteworthy that the finest WG finite element space has little bearing 
on the coarse mesh $\mathcal T_H$ mesh size selection in augmented subspace techniques.
Compared to the two-grid WG finite element technique 
\cite{ZhaiHuZhang,ZhaiXieZhangZhang_TwoGrid}, 
wherein the choices of coarse and fine meshes are not free each other, 
this represents a significant distinction.

Here, we are concerned with the following standard Laplace eigenvalue problem: 
Find $(\lambda,u)\in\mathbb R\times H_0^1(\Omega)$ such that
\begin{eqnarray}\label{Test_Problem}
\left\{
\begin{array}{rcl}
-\Delta u&=& \lambda u,\ \ \ {\rm in}\ \Omega,\\
u&=&0,\ \ \ \ {\rm on}\ \partial\Omega,\\
|u|_1^2&=&1,
\end{array}
\right.
\end{eqnarray}
where $|\cdot|_1$ represents $H^1$-type semi-norm and the computing domain is set 
to be the unit square $\Omega=(0,1)\times (0,1)$.
Then, in (\ref{weak_eigenvalue_problem}), the bilinear forms $a(\cdot, \cdot)$ 
and $b(\cdot,\cdot)$ are defined as follows
\begin{eqnarray*}
a(u,v)=\int_{\Omega}\nabla u\cdot\nabla vd\Omega, &&
b(u,v)=\int_{\Omega}uvd\Omega.
\end{eqnarray*}
Additionally, the norms $\left\|\cdot\right\|_{a,h}$ and $\left\|\cdot\right\|_{b,h}$
defined in (\ref{Nomr_a}) and (\ref{Nomr_b}) are equivalent to the $H^1$-type 
semi-norm $|\cdot|_1$ and $L^2$ norm $\left\|\cdot\right\|_0$, respectively.
In order to use the WG finite element discretization method, 
we employ the meshes defined in Section 2.

Here, the problem (\ref{Test_Problem}) is treated using the augmented subspace 
techniques specified by Algorithms \ref{Algorithm_k} and \ref{Algorithm_1}.
In this section, the regular refinement is used to create the fine mesh $\mathcal T_h$ 
from the coarse mesh $\mathcal T_H$.  The WG finite element space on the 
fine mesh $\mathcal T_h$ is set to $V_h$, and the coarse conforming linear 
finite element space on the coarse mesh $\mathcal T_H$ is set to $W_H$. 
We consider the computational domain $\Omega$ to be convex for the sake of simplicity.

In order to give the explicit convergence rate of the augmented subspace 
methods defined by Algorithms \ref{Algorithm_k} and \ref{Algorithm_1},  
we  need to estimate the quantity $\eta_a(W_H)$ in (\ref{eta_a_H_Def}).
For this aim, we define the conforming linear finite element projection 
$\mathcal P_H: H_0^1(\Omega)\mapsto W_H$
as follows
\begin{eqnarray}
a(\mathcal P_Hw, v_H) = a(w, v_H),\ \ \ \forall v_H\in W_H,\ \ {\rm for}\ w\in H_0^1(\Omega).
\end{eqnarray}
It is well known that the following error estimate holds
\begin{eqnarray}
\|Tf-\mathcal P_HTf\|_1 \leq CH \|Tf\|_2 \leq CH\|f\|_{b,h},
\end{eqnarray}
where $T: L^2(\Omega)\mapsto H_0^1(\Omega)$ is defined as follows
\begin{eqnarray}
a(Tf, v) = b(f,v),\ \ \ \ \forall v\in H_0^1(\Omega).
\end{eqnarray}

In order to deduce the estimate for the term $\eta_a(W_H)$, 
we define the norm $\|\cdot\|_{1,h}$ as follows
\begin{eqnarray*}
\|v\|_{1,h} ^2= \sum_{K\in\mathcal T_h}
\left(\|\nabla v_0\|_{0,K}^2 +  h_K^{-1}\|v_0-v_b\|_{\partial K}^2\right).
\end{eqnarray*}
Obviously, the norm $\|\cdot\|_{1,h}$ coincides with $\|\cdot\|_1$ 
on the Sobolev space $H_0^1(\Omega)$. Furthermore, there is the following 
equivalence between $\|\cdot\|_{1,h}$ and $\|\cdot\|_{a,h}$ on the 
WG finite element space $V_h$.
\begin{lemma}(\cite{MuWangWangYe})
For any $v_h\in V_h$, the following inequalities hold
\begin{eqnarray}
C_7 \|v_h\|_{1,h} \leq \|v_h\|_{a,h} \leq C_8 \|v_h\|_{1,h},
\end{eqnarray}
where $C_7$ and $C_8$ are two constants independent of the mesh size $h$.
\end{lemma}

Then $\|T_hf-\mathcal P_HTf\|_{a,h} $ has following inequalities
\begin{eqnarray}\label{Inequality_45}
&&\|T_hf-\mathcal P_HTf\|_{a,h} \leq \|T_hf-Q_hTf\|_{a,h} + \|Q_hTf-\mathcal P_HTf\|_{a,h}\nonumber\\
&&\leq  \|T_hf-Q_hTf\|_{a,h} + C\|Q_hTf-\mathcal P_HTf\|_{1,h}\nonumber\\
&&\leq  \|T_hf-Q_hTf\|_{a,h} + C\|Q_hTf-Tf\|_{1,h}+ C\|Tf-\mathcal P_HTf\|_{1,h}\nonumber\\
&&\leq  \|T_hf-Q_hTf\|_{a,h} + C\|Q_hTf-Tf\|_{1,h}+ C\|Tf-\mathcal P_HTf\|_1\nonumber\\
&&\leq C(h+h+H)\|Tf\|_2 \leq CH\|f\|_{b,h},
\end{eqnarray}
where the constant depends on the shape of the mesh $\mathcal T_H$.

From the definition of $\eta_a(W_H)$ in (\ref{eta_a_H_Def}), and (\ref{Inequality_45}),
we can obtain the following estimates
\begin{eqnarray}
\eta_a(W_H)&\leq& \sup_{\substack{ f\in L^2(\Omega)\\ \left\|f\right\|_{b,h}=1}}
\inf_{w_H\in W_H}\left\|T_hf-w_H\right\|_{a,h}
\leq \sup_{\substack{ f\in L^2(\Omega)\\ \left\|f\right\|_{b,h}=1}}
\left\|T_hf-\mathcal P_HTf\right\|_{a,h} \nonumber\\
&\leq& \sup_{\substack{ f\in L^2(\Omega)\\ \left\|f\right\|_{b,h}=1}} CH\|f\|_{b,h} 
= CH.
\end{eqnarray}

Based on Theorems \ref{Theorem_Error_Estimate_k} and \ref{Theorem_Error_Estimate_1},
the convergence results for the augmented subspace method can be concluded 
with the following inequalities
\begin{eqnarray}
&&\left\|\bar u_{i,h} -F_{k,h}^{(\ell+1)}\bar u_{i,h}\right\|_{a,h}
\leq C\big(CH\big)^{2\ell}\left\|\bar u_{i,h} - F_{k,h}^{(1)}\bar u_{i,h}\right\|_{a,h},
\ \ \ \ i=1, \cdots, k,\label{Test_1_1}\\
&&\left\|\bar u_{i,h} -F_{k,h}^{(\ell+1)}\bar u_{i,h} \right\|_{b,h}
\leq CH \left\|\bar u_{i,h} -F_{k,h}^{(\ell+1)}\bar u_{i,h}\right\|_{a,h},\ \ \ \ i=1, 
\cdots, k,\label{Test_1_0}
\end{eqnarray}
and
\begin{eqnarray}
\left\|\bar u_{i,h}-E_{i,h}^{(\ell+1)}\bar u_{i,h}\right\|_{a,h} &\leq& C\big(CH\big)^{2\ell}
\left\|\bar u_{i,h}-E_{i,h}^{(1)}\bar u_{i,h}\right\|_{a,h},\label{Test_2_1}\\
\left\|\bar u_h-E_{i,h}^{(\ell+1)}\bar u_h\right\|_{b,h}
&\leq& CH\left\|\bar u_{i,h}-E_{i,h}^{(\ell+1)}\bar u_{i,h}\right\|_{a,h}.\label{Test_2_0}
\end{eqnarray}
The goal of this section is to validate these convergence findings using a 
few numerical examples.The exact WG finite element eigenfunction can 
be found by directly solving the eigenvalue problem on the fine WG 
finite element space $V_h$. Let this be noted.
To aid with comprehension, the nomenclature in all of the following 
figures denotes the exact WG finite element eigenfunctions and 
the augmented subspace approximations, respectively, with and 
without the ``{\tt dir}" superscript.

\subsection{Augmented subspace method for $P_0/P_0$ WG finite element space}
For the WG finite element space $P_0/P_0$, we examine the performance of 
the augmented subspace approach described by Algorithms \ref{Algorithm_k} 
and \ref{Algorithm_1} in the first subsection. Here, $W_H$ is defined 
as the conforming linear finite element space on the coarse 
mesh $\mathcal T_H$ in all numerical cases.
The $P_0/P_0$ WG finite element space $V_h$ defined on the finer mesh $\mathcal T_h$ 
can be written as follows 
\begin{eqnarray*}
V_h = \Big\{ v: v|_{K_0} \in \mathcal P_0(K_0)\ {\rm for}\ K\in \mathcal T_h;
  v|_e \in \mathcal P_0(e) \ {\rm for}\ e\in\mathcal E_h, \ {\rm and}\ v|_e = 0\ {\rm for}\
  e\in \mathcal E_h\cap\partial\Omega\Big\}.
\end{eqnarray*}
The fine mesh $\mathcal T_h$ is obtained from the coarse mesh $\mathcal T_H$ 
by the regular refinement. Here, we set the size $h=\sqrt{2}/256$ 
for the fine mesh $\mathcal T_h$.

We also verify the convergence results for the conforming linear finite element space $W_H$ 
with various sizes $H$ by examining the numerical errors corresponding to the results in (\ref{Test_1_1})-(\ref{Test_2_0}).  The goal is to determine how the mesh size 
$H$ affects the convergence rate. In this case, the regular type of 
quasiuniform mesh $\mathcal T_H$ is also specified as the coarse mesh.

\begin{figure}[http!]
\centering
\includegraphics[width=7cm,height=4.5cm]{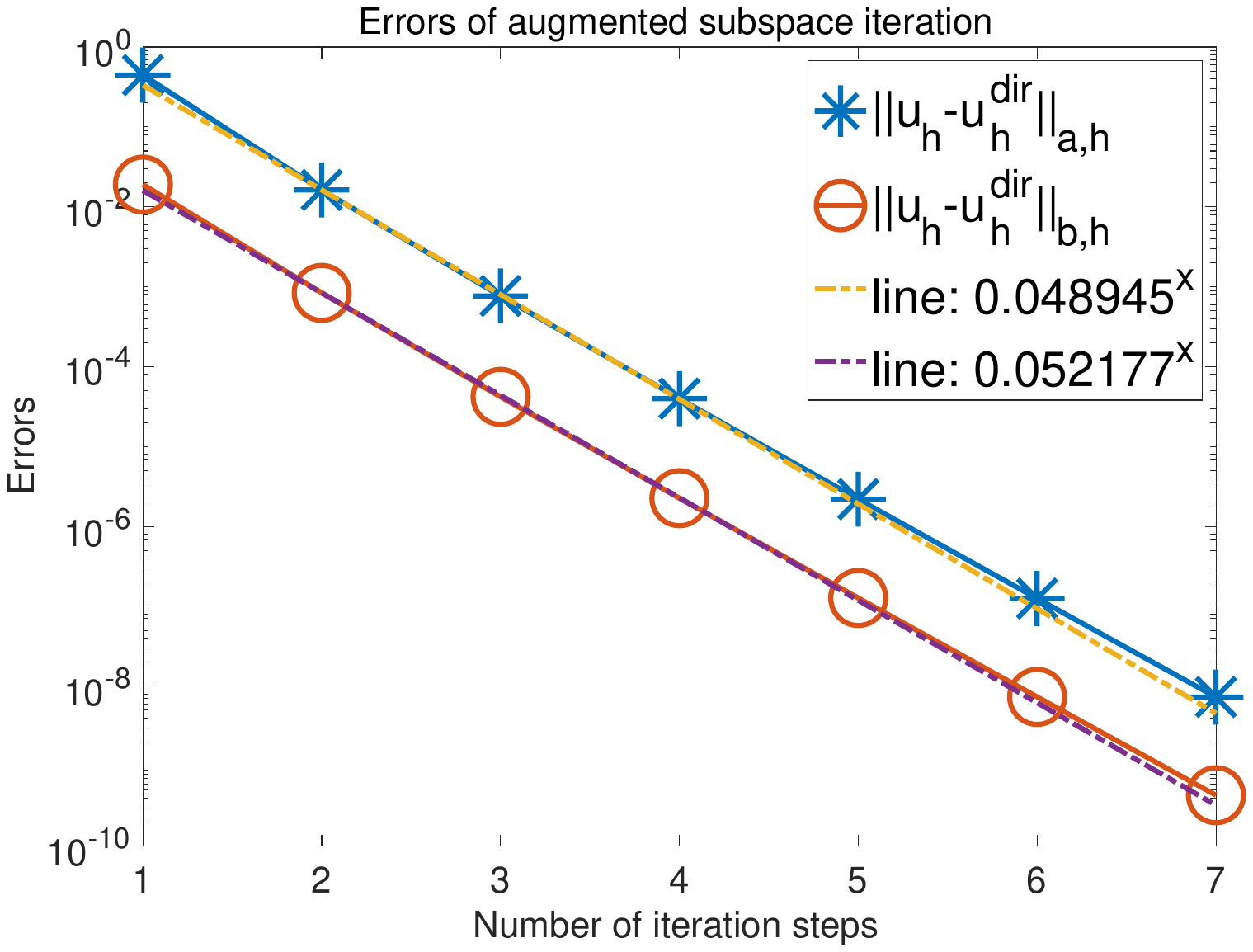}
\includegraphics[width=7cm,height=4.5cm]{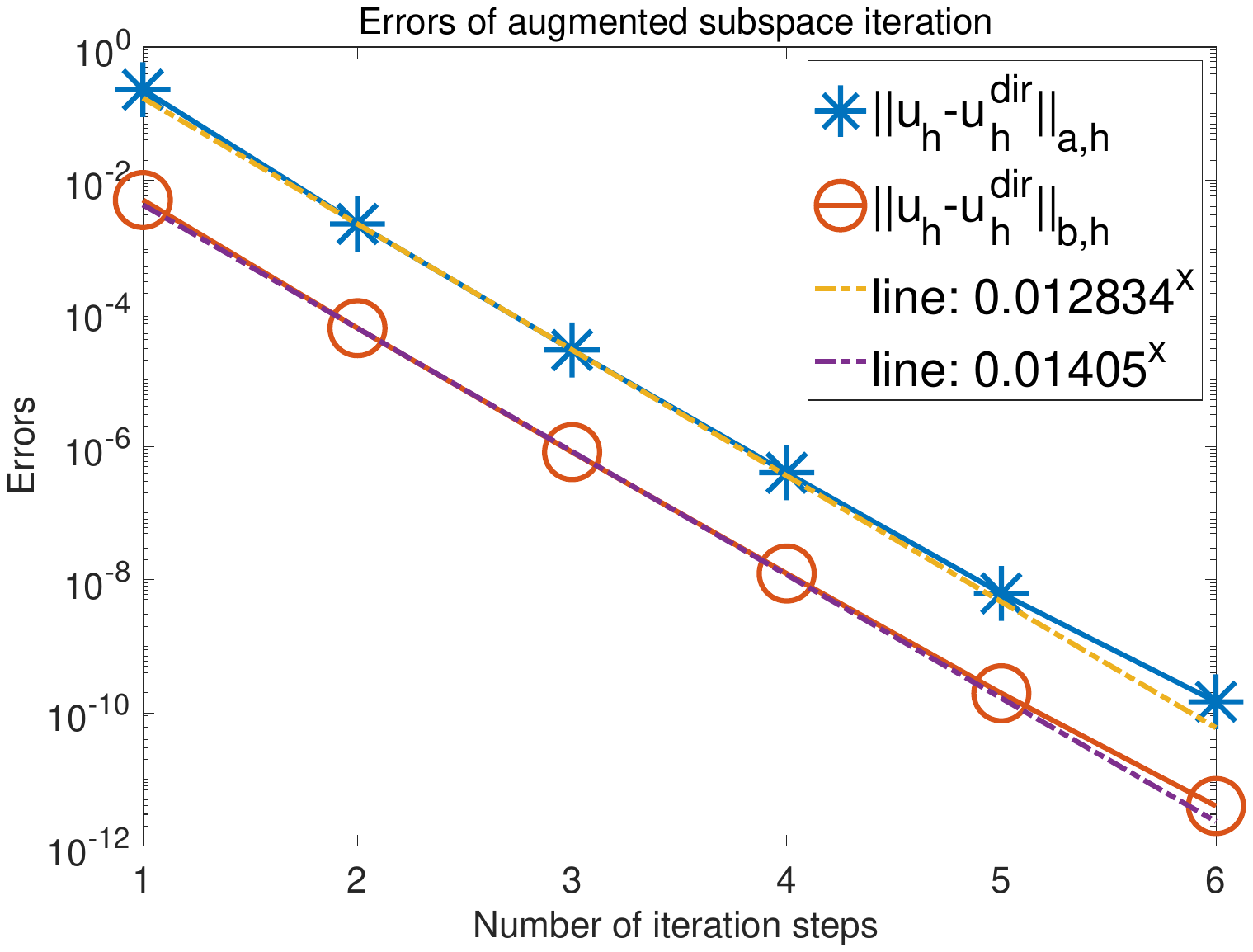}\\
\includegraphics[width=7cm,height=4.5cm]{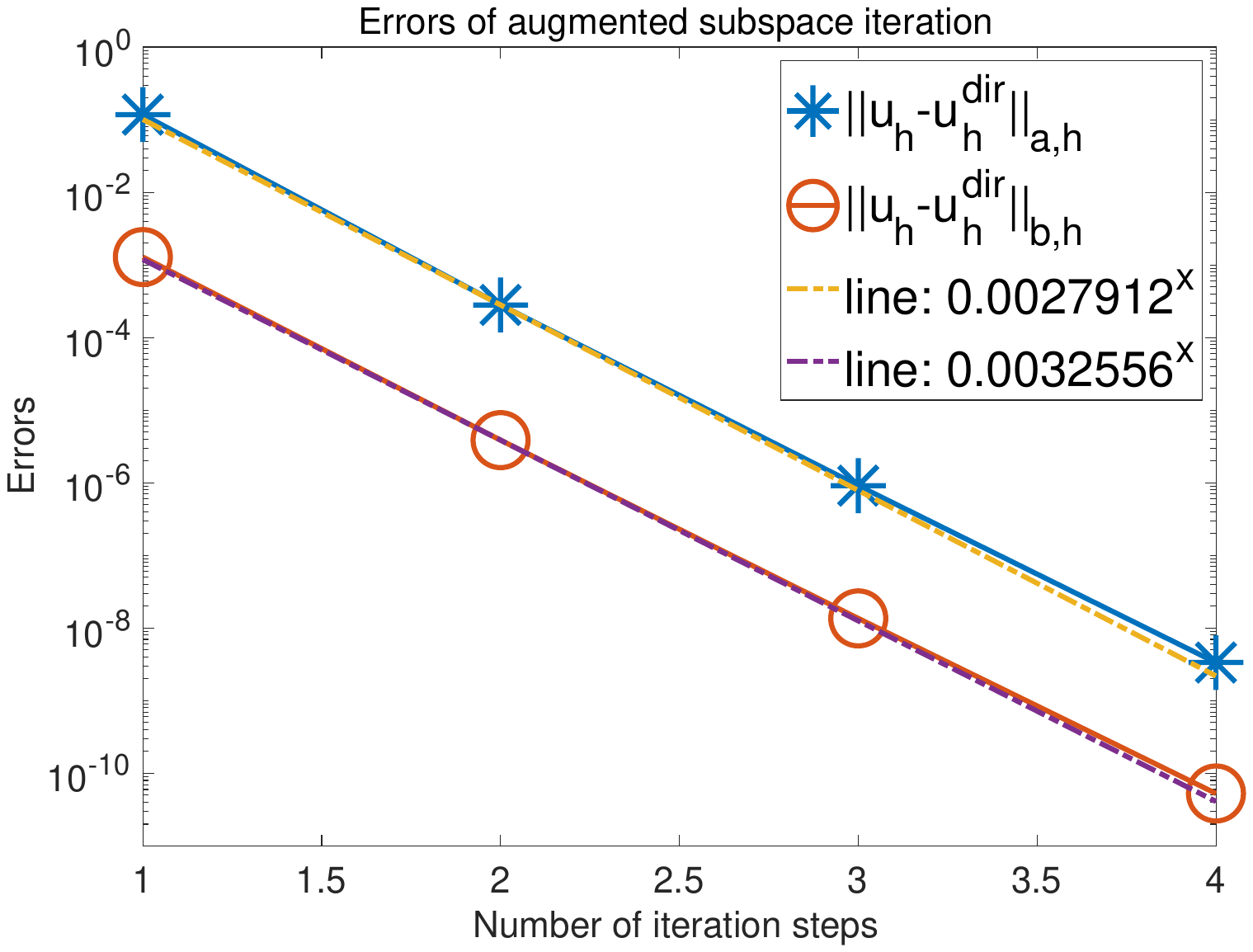}
\includegraphics[width=7cm,height=4.5cm]{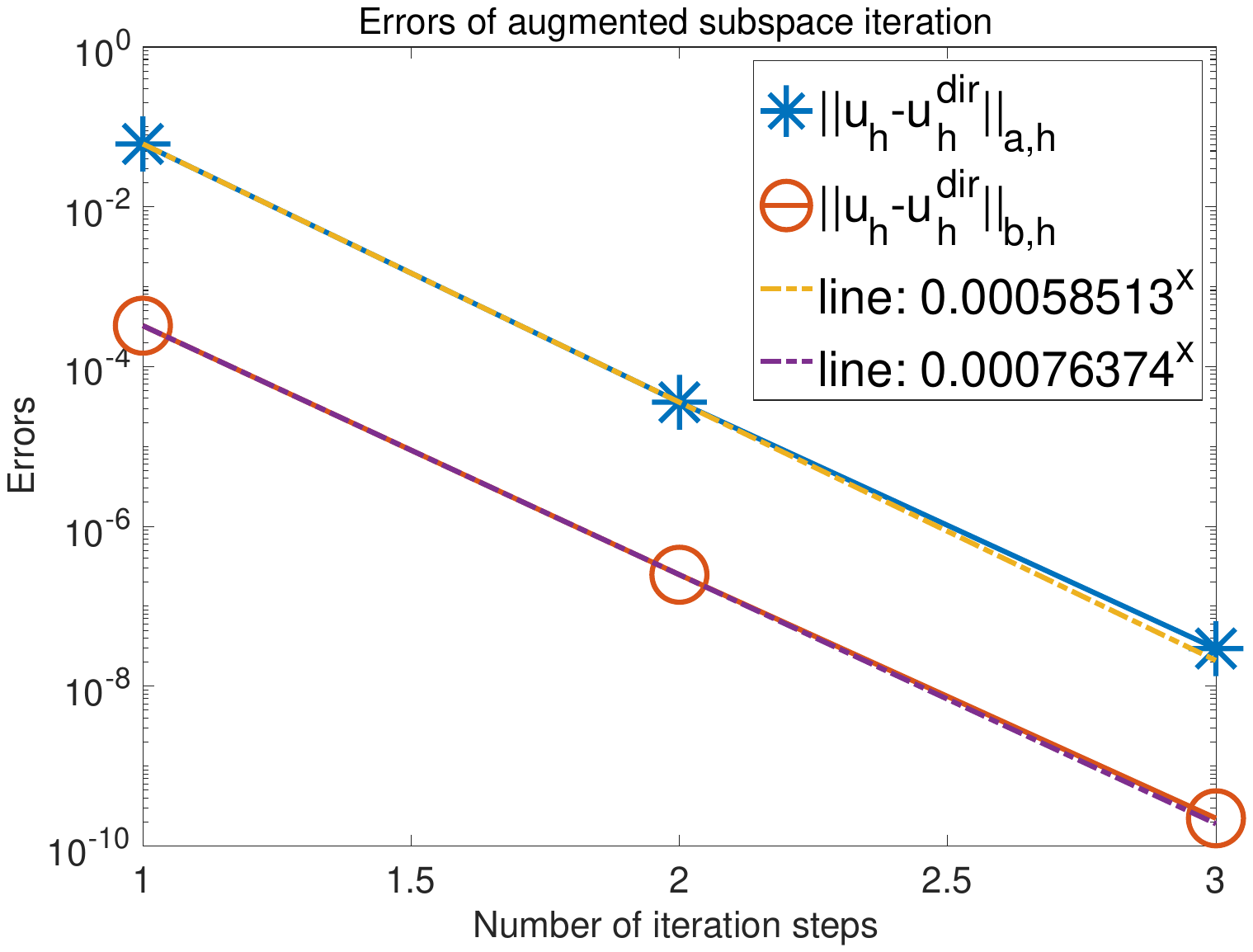}
\caption{The convergence behaviors for the first eigenfunction by Algorithm \ref{Algorithm_k}
corresponding to the $P_0/P_0$ WG finite element method and the coarse mesh 
size $H=\sqrt{2}/8$, $\sqrt{2}/16$, $\sqrt{2}/32$ and $\sqrt{2}/64$, 
respectively.}\label{Result_Coarse_Mesh}
\end{figure}

Under the boundary condition restriction, the initial eigenfunction approximation 
is specified to be rand vectors in this case. Next, we employ the augmented subspace 
approach, as specified by Algorithms \ref{Algorithm_k} and \ref{Algorithm_1}, 
to carry out the iteration steps.
The convergence behaviors for the first eigenfunction using the augmented subspace 
techniques are displayed in Figure \ref{Result_Coarse_Mesh}, 
and they correspond to the coarse mesh sizes $H=\sqrt{2}/8$, $\sqrt{2}/16$, 
$\sqrt{2}/32$, and $\sqrt{2}/64$, respectively.
The rates of convergence associated with $\|\cdot\|_{a,h}$ 
and $\left\|\cdot\right\|_{b,h}$ are, respectively, $0.048945$, $0.012834$, $0.00279122$, 
 $0.00058513$ and $0.052177$, $0.01405$, $0.0032556$, $0.00076374$.
As a consequence, the results (\ref{Test_1_1})-(\ref{Test_2_0}) hold and 
validate the second order convergence speed of the augmented subspace 
technique described by Algorithms \ref{Algorithm_k} and \ref{Algorithm_1}.

\begin{figure}[http!]
\centering
\includegraphics[width=7cm,height=4.5cm]{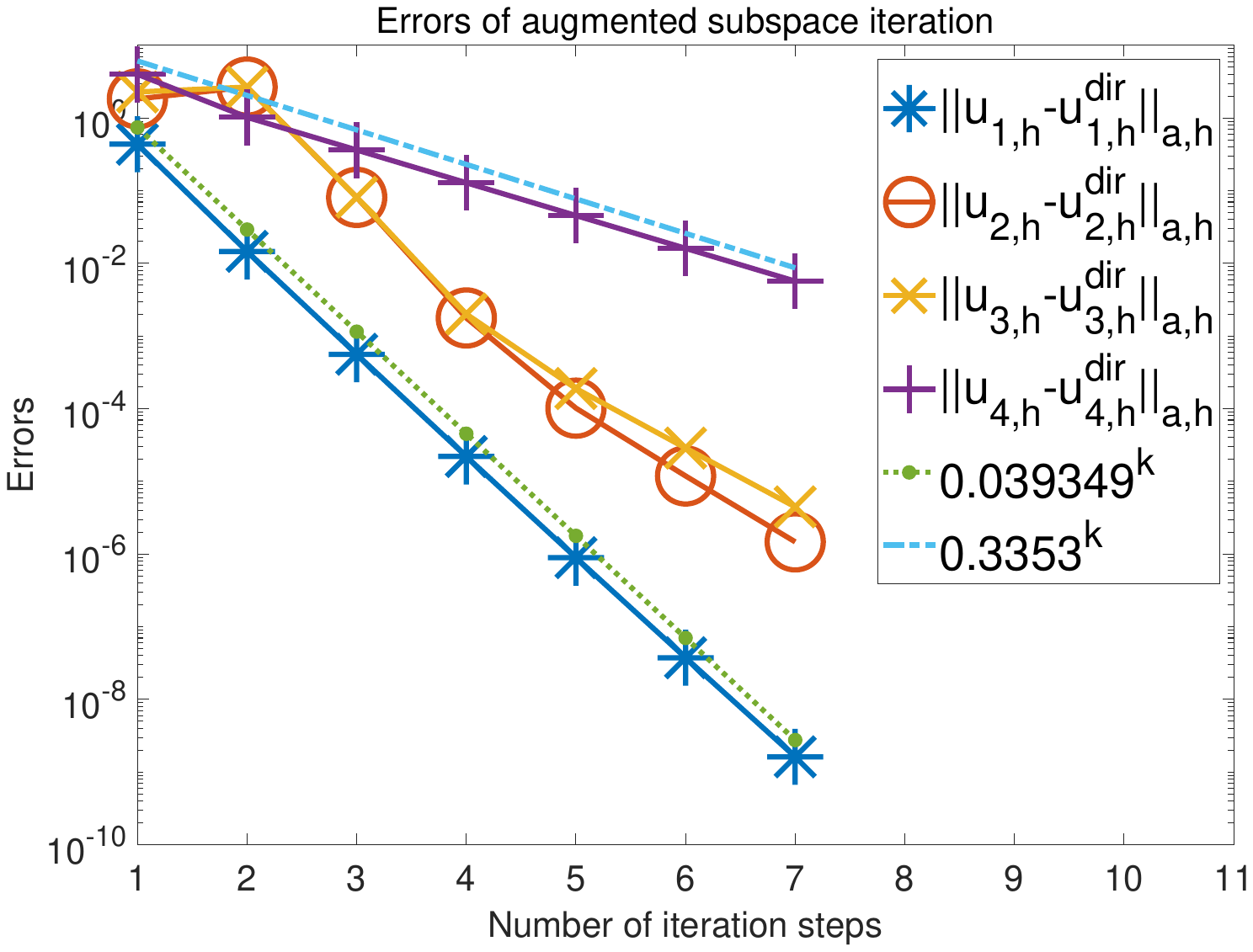}
\includegraphics[width=7cm,height=4.5cm]{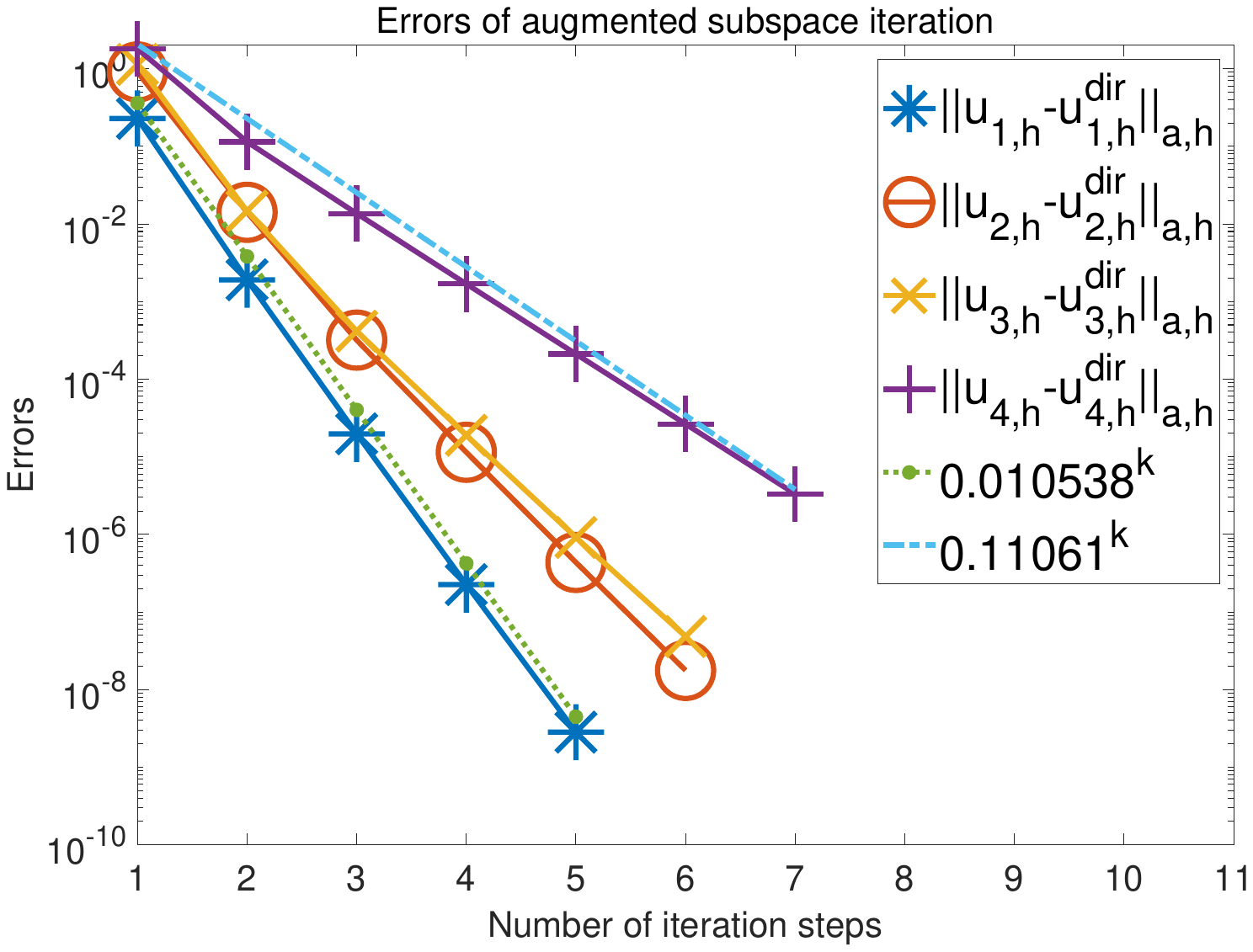}\\
\includegraphics[width=7cm,height=4.5cm]{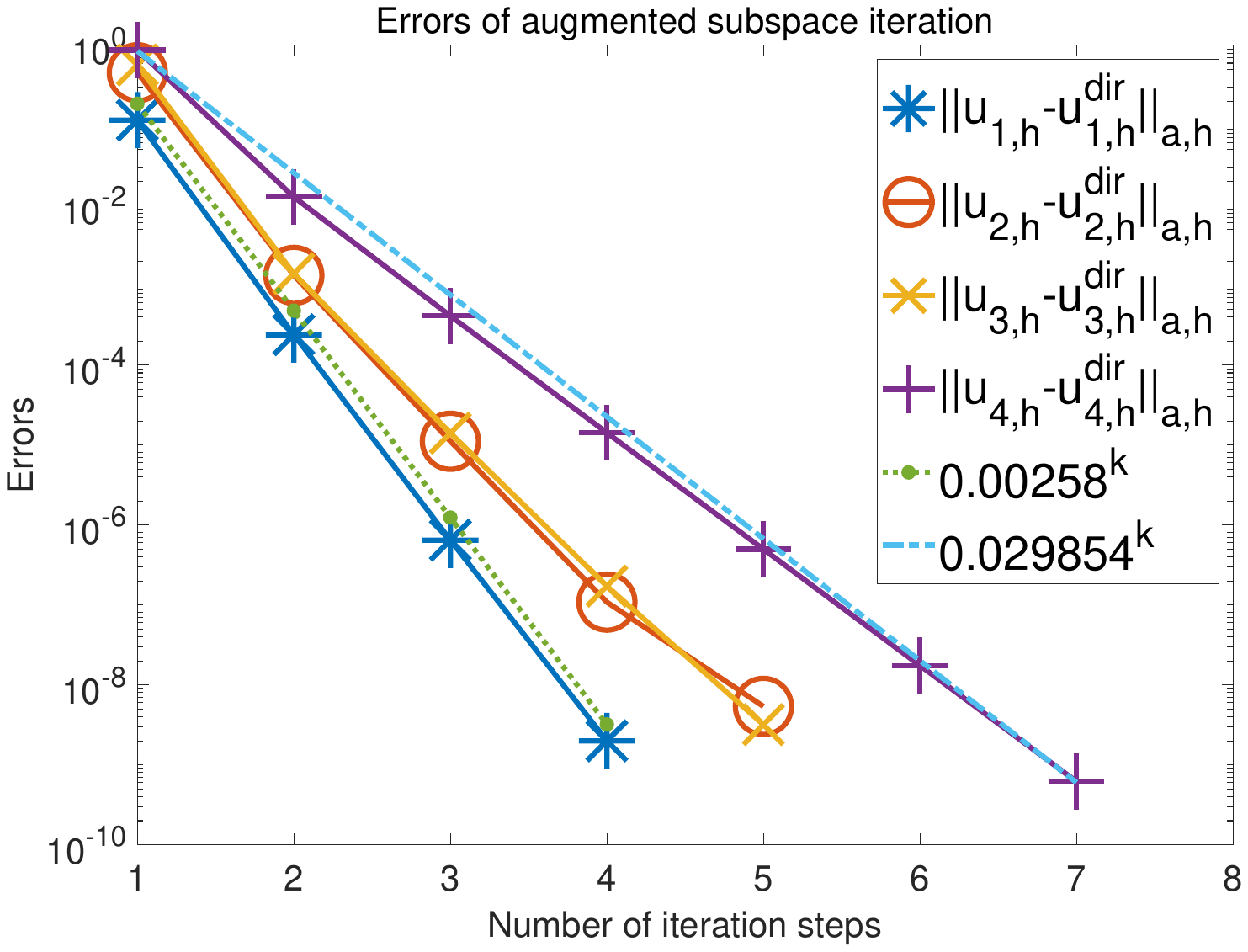}
\includegraphics[width=7cm,height=4.5cm]{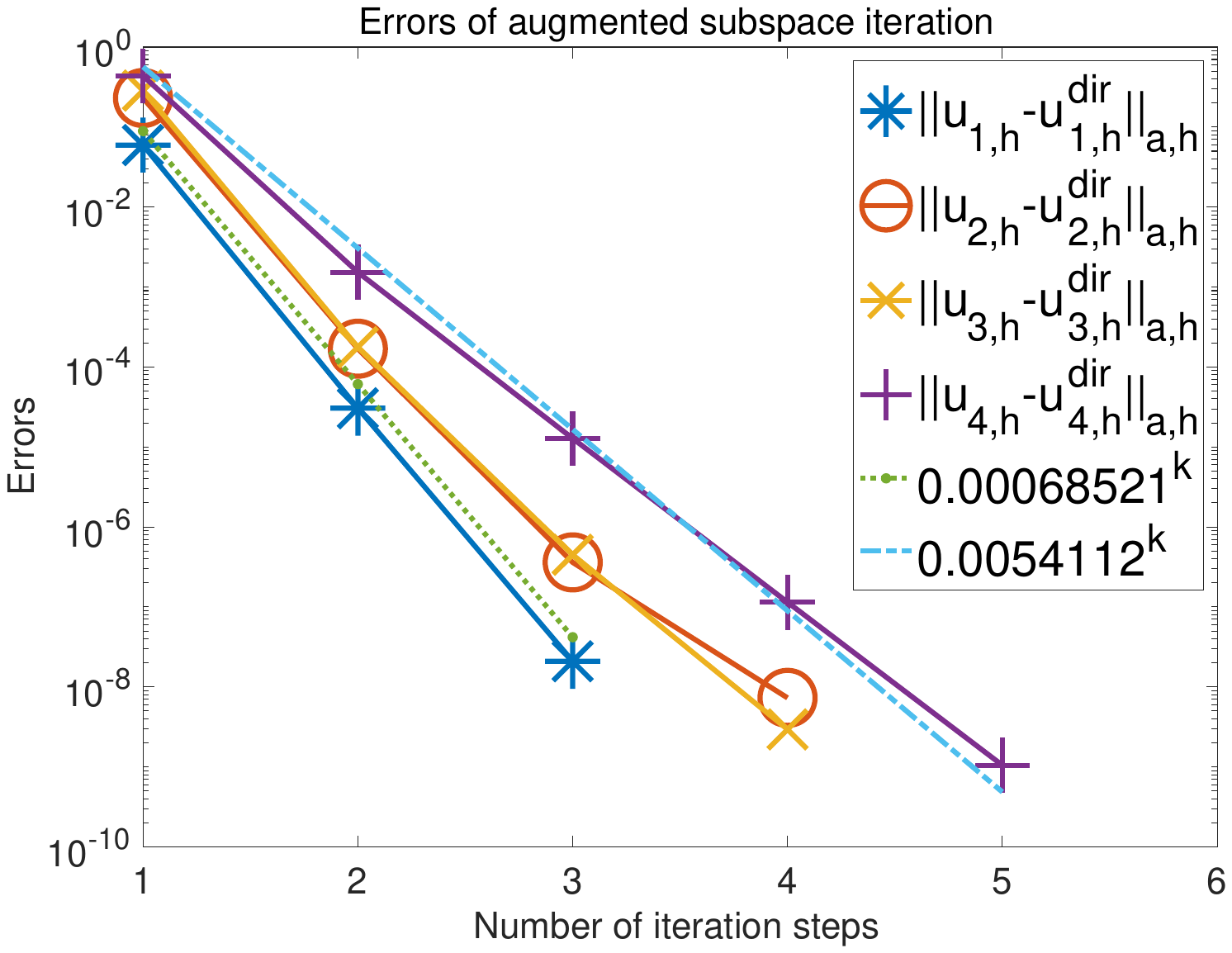}
\caption{The convergence behaviors for the smallest $4$ eigenfunctions by Algorithm \ref{Algorithm_k}
with the $P_0/P_0$ WG finite element method and the coarse space being the linear finite element space on the
mesh with size $H=\sqrt{2}/8$, $\sqrt{2}/16$, $\sqrt{2}/32$ and $\sqrt{2}/64$, respectively.}\label{Result_Coarse_Mesh_4}
\end{figure}

Next, we evaluate Algorithm \ref{Algorithm_k} in terms of its ability to compute the first $4$ eigenpairs.
The corresponding convergence behaviors for the smallest $4$ eigenfunctions by Algorithm \ref{Algorithm_k} are displayed in Figure \ref{Result_Coarse_Mesh_4}. The conforming linear finite element space on the mesh with sizes $H=\sqrt{2}/8$, $\sqrt{2}/16$, $\sqrt{2}/32$, and $\sqrt{2}/64$, respectively, forms the coarse space $W_H$.
Employing the $4$-th eigenfunction as an example, we can determine the related convergence rates, 
which indicate the second  convergence order of the algorithm given by Algorithm \ref{Algorithm_k}, to be $0.3353$, $0.11061$, $0.029854$, and $0.0054112$.
Furthermore, we are able to observe from Figure \ref{Result_Coarse_Mesh_4} that the $4$-th eigenfunction's convergence rate is slower than the $1$-st eigenfunction's, which is in accordance with Theorem \ref{Theorem_Error_Estimate_k}.

\begin{figure}[http!]
\centering
\includegraphics[width=7cm,height=4.5cm]{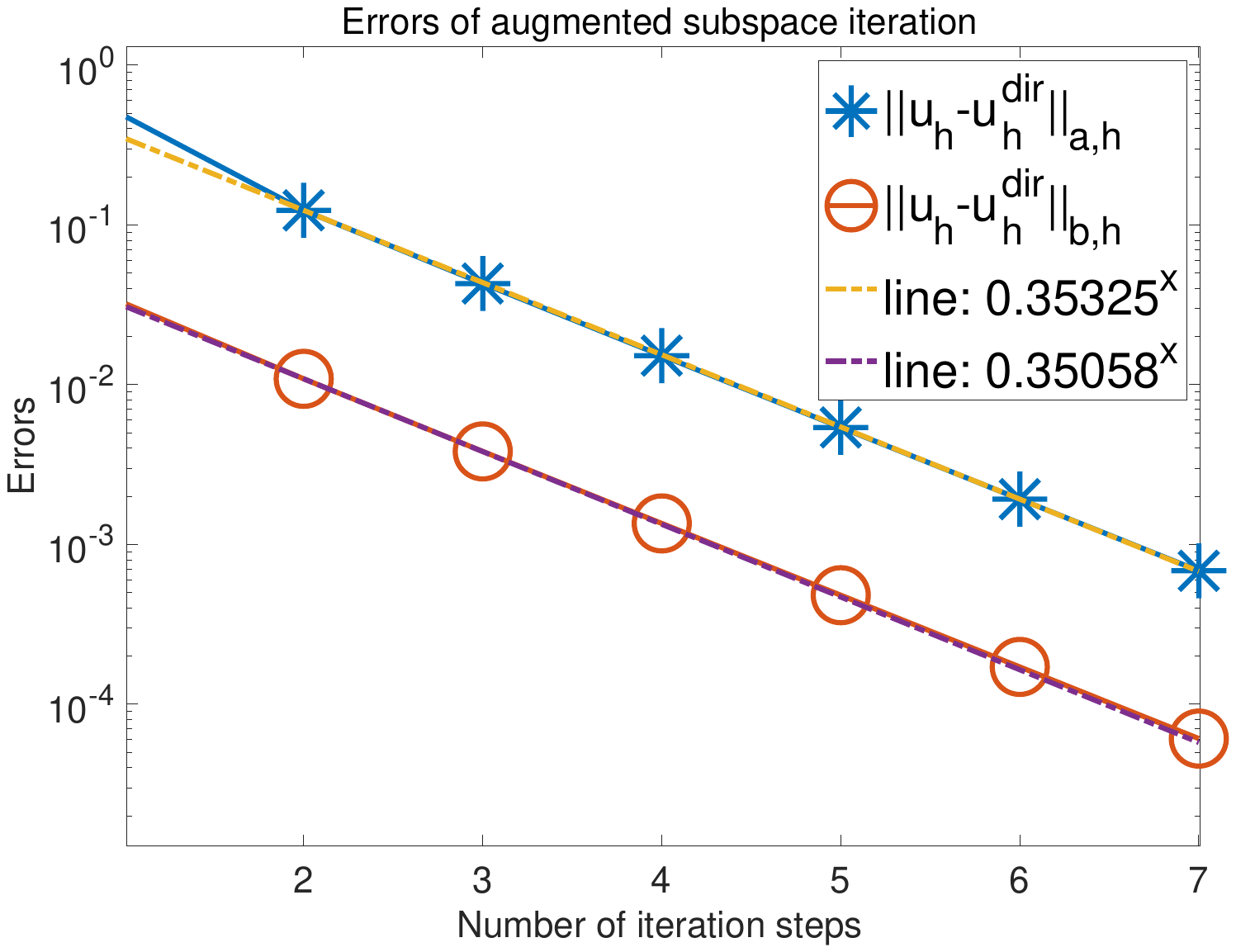}
\includegraphics[width=7cm,height=4.5cm]{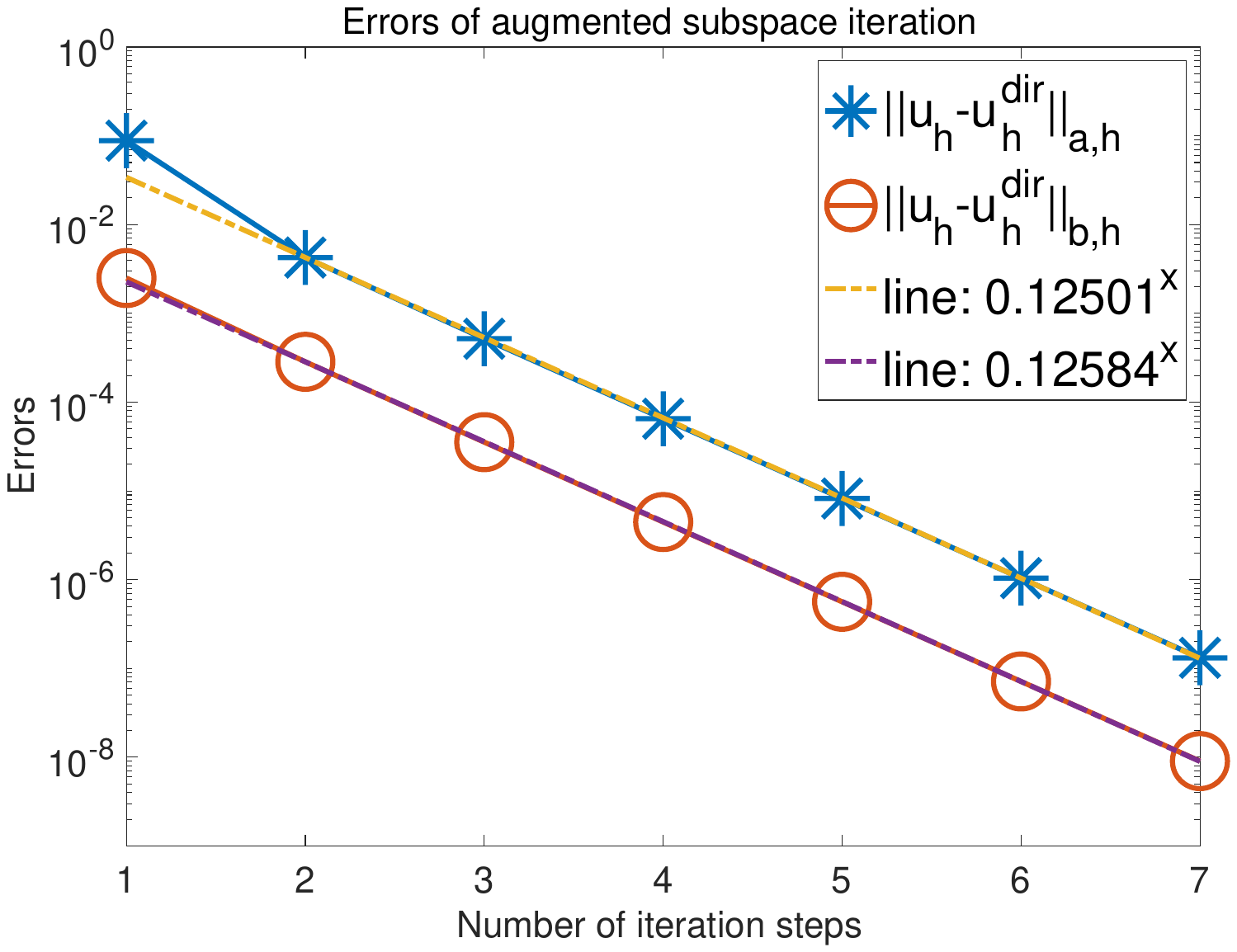}\\
\includegraphics[width=7cm,height=4.5cm]{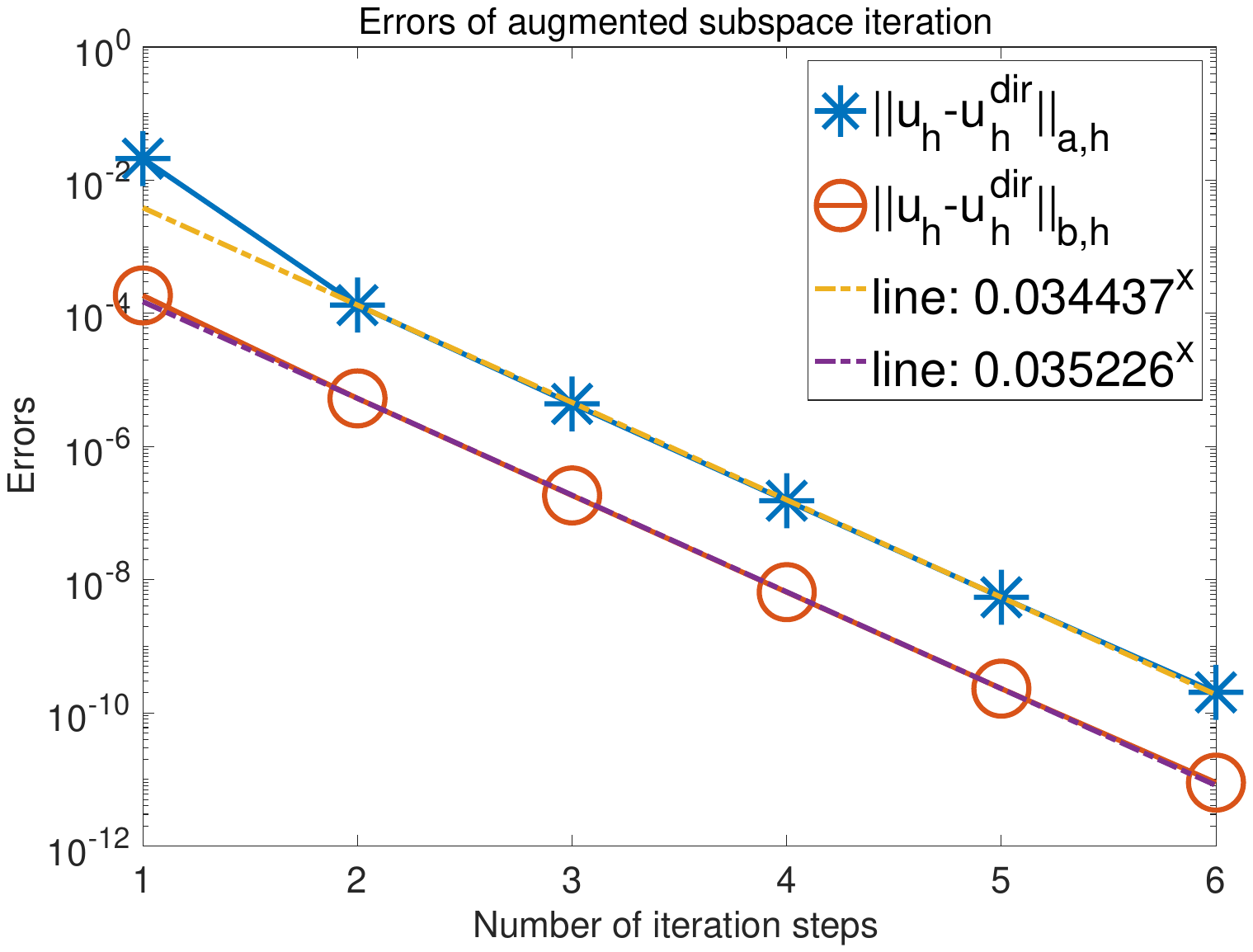}
\includegraphics[width=7cm,height=4.5cm]{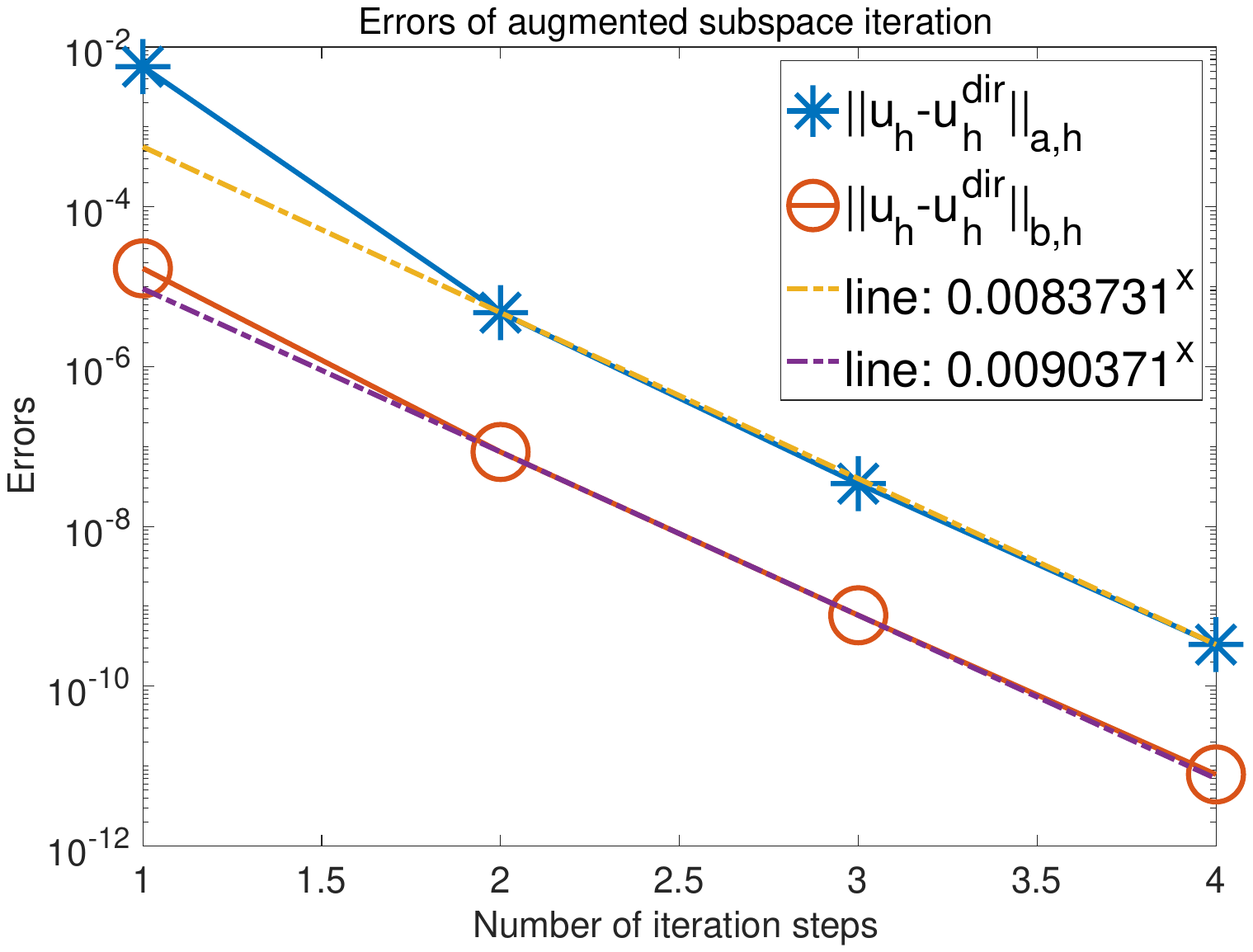}
\caption{The convergence behaviors for the only $4$-th eigenfunction by Algorithm \ref{Algorithm_1}
with the $P_0/P_0$ WG finite element method and  the coarse space being the linear finite element space on
the mesh with size $H=\sqrt{2}/8$, $\sqrt{2}/16$,
$\sqrt{2}/32$ and $\sqrt{2}/64$, respectively.}\label{Result_Coarse_Mesh_4_Only}
\end{figure}

Assessing Algorithm \ref{Algorithm_1}'s performance in determining 
the single $4$-th eigenpair is the next objective.
Since smallest eigenpairs is not the goal, the eigenvalue problem (\ref{weak_eigenvalue_problem}) 
is solved on the coarse WG finite element space $V_H$ to provide the initial eigenfunction approximation.
The augmented subspace approach, which is specified by Algorithm \ref{Algorithm_1}, 
is then used to carry out the iteration phases.
The coarse space was the linear finite element space on the mesh with sizes 
$H=\sqrt{2}/8$, $\sqrt{2}/16$, $\sqrt{2}/32$, and $\sqrt{2}/64$, respectively. 
The corresponding convergence behaviors for the only $4$-th eigenfunction by 
Algorithm \ref{Algorithm_1} are depicted in Figure \ref{Result_Coarse_Mesh_4_Only}.
The norms $\|\cdot\|_{a,h}$ and $\left\|\cdot\right\|_{b,h}$ 
in Figure \ref{Result_Coarse_Mesh_4_Only} correspond to the convergence rates, 
which are $0.35325$, $0.12501$, $0.034437$ and $0.0083731$, 
and $0.35058$, $0.12584$, $0.035226$ and $0.0090371$, respectively. 
According to these findings, the augmented subspace approach described by 
Algorithm \ref{Algorithm_1} has a second order speed of convergence, validating the findings of (\ref{Test_2_1})-(\ref{Test_2_0}).

\subsection{Augmented subspace method for $P_1/P_1$ WG finite element space}
We examine the augmented subspace method's performance for 
the WG finite element space $P_1/P_1$, as described by 
Algorithms \ref{Algorithm_k} and \ref{Algorithm_1}, in the second subsection.
Also, $W_H$ is designated as the conforming linear finite element 
space on the coarse mesh $\mathcal T_H$ in these numerical tests.
Here, $V_h$ is the $P_1/P_1$ WG finite element space defined 
on the fine mesh $\mathcal T_h$, which is generated by 
the regular refinement from the coarse mesh $\mathcal T_H$.

Here, we set the size $h=\sqrt{2}/256$ for the fine mesh $\mathcal T_h$ and the WG finite element 
space $V_h$ is defined as follows 
\begin{eqnarray*}
V_h = \Big\{ v: v|_{K_0} \in \mathcal P_1(K_0)\ {\rm for}\ K\in \mathcal T_h;
v|_e \in \mathcal P_1(e) \ {\rm for}\ e\in\mathcal E_h, \ {\rm and}\ v|_e = 0\ {\rm for}\
e\in \mathcal E_h\cap\partial\Omega\Big\}.
\end{eqnarray*}

We also check the numerical errors corresponding to the conforming linear finite 
element space $W_H$ with different sizes $H$. This helps to confirm the 
convergence results for the $P_1/P_1$ WG finite element 
technique described in (\ref{Test_1_1})-(\ref{Test_2_0}).
Here, also determining how the convergence rate varies with mesh size $H$ is a goal.
In this case, the regular type of quasiuniform mesh $\mathcal T_H$ 
is also specified as the coarse mesh.

\begin{figure}[http!]
\centering
\includegraphics[width=7cm,height=4.5cm]{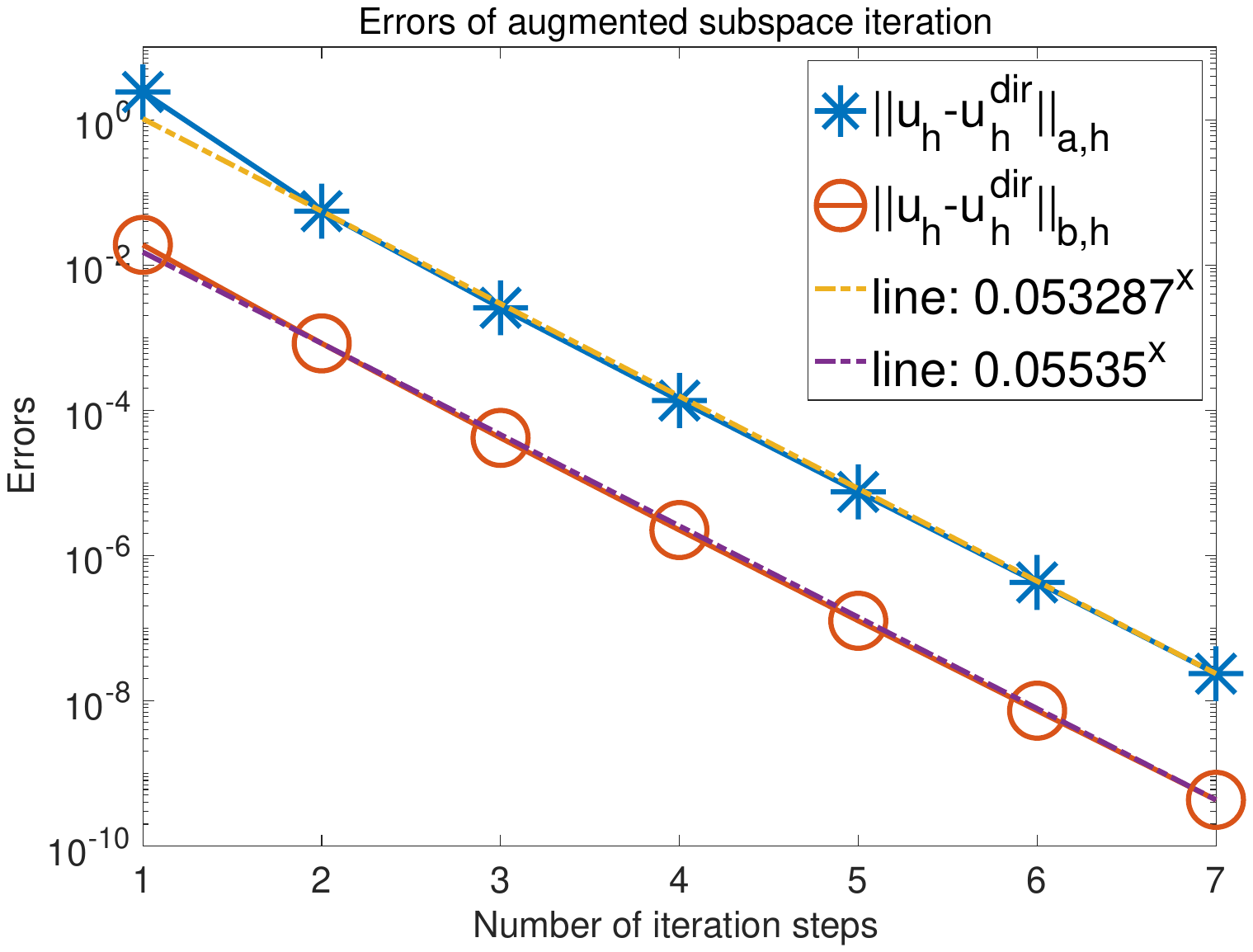}
\includegraphics[width=7cm,height=4.5cm]{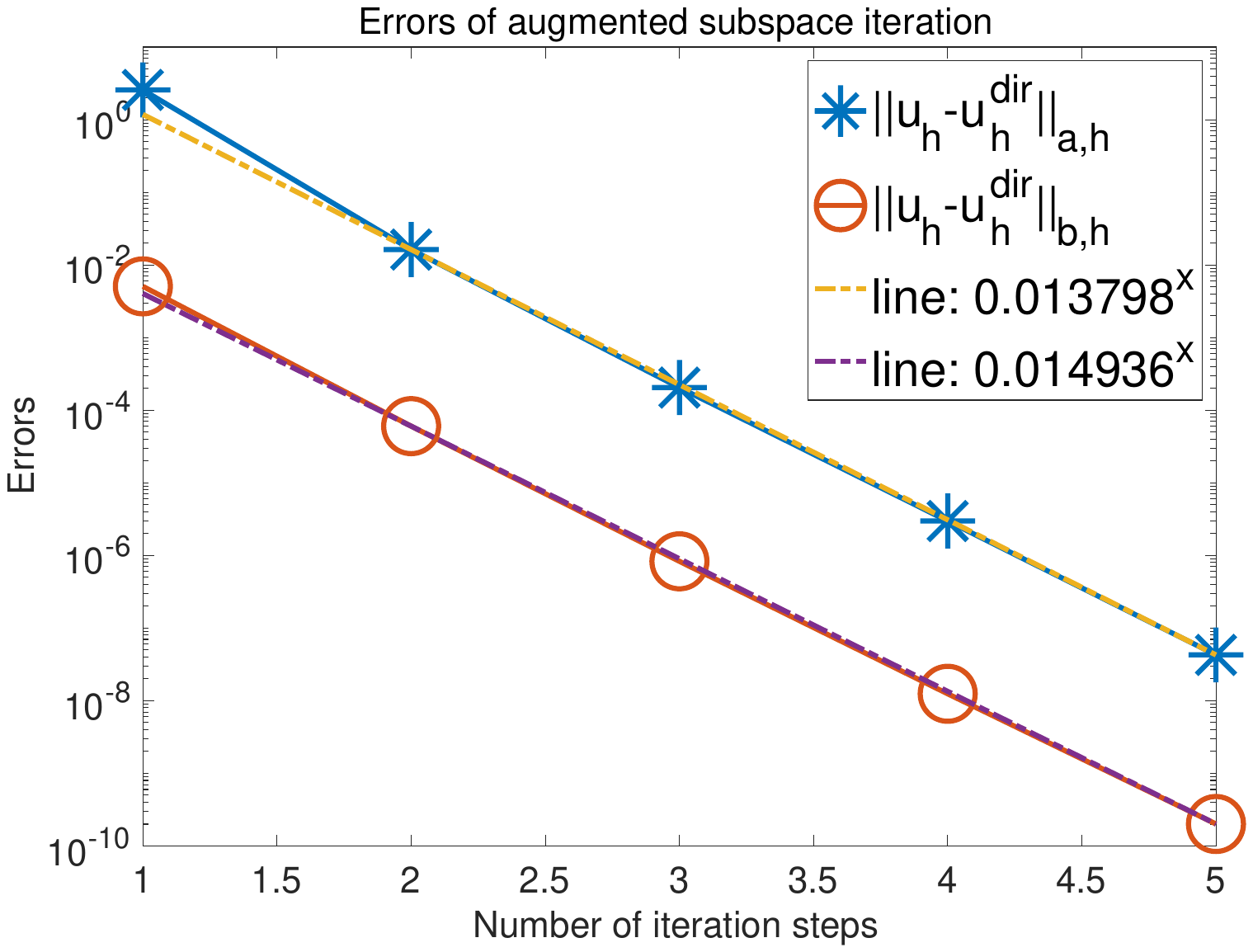}\\
\includegraphics[width=7cm,height=4.5cm]{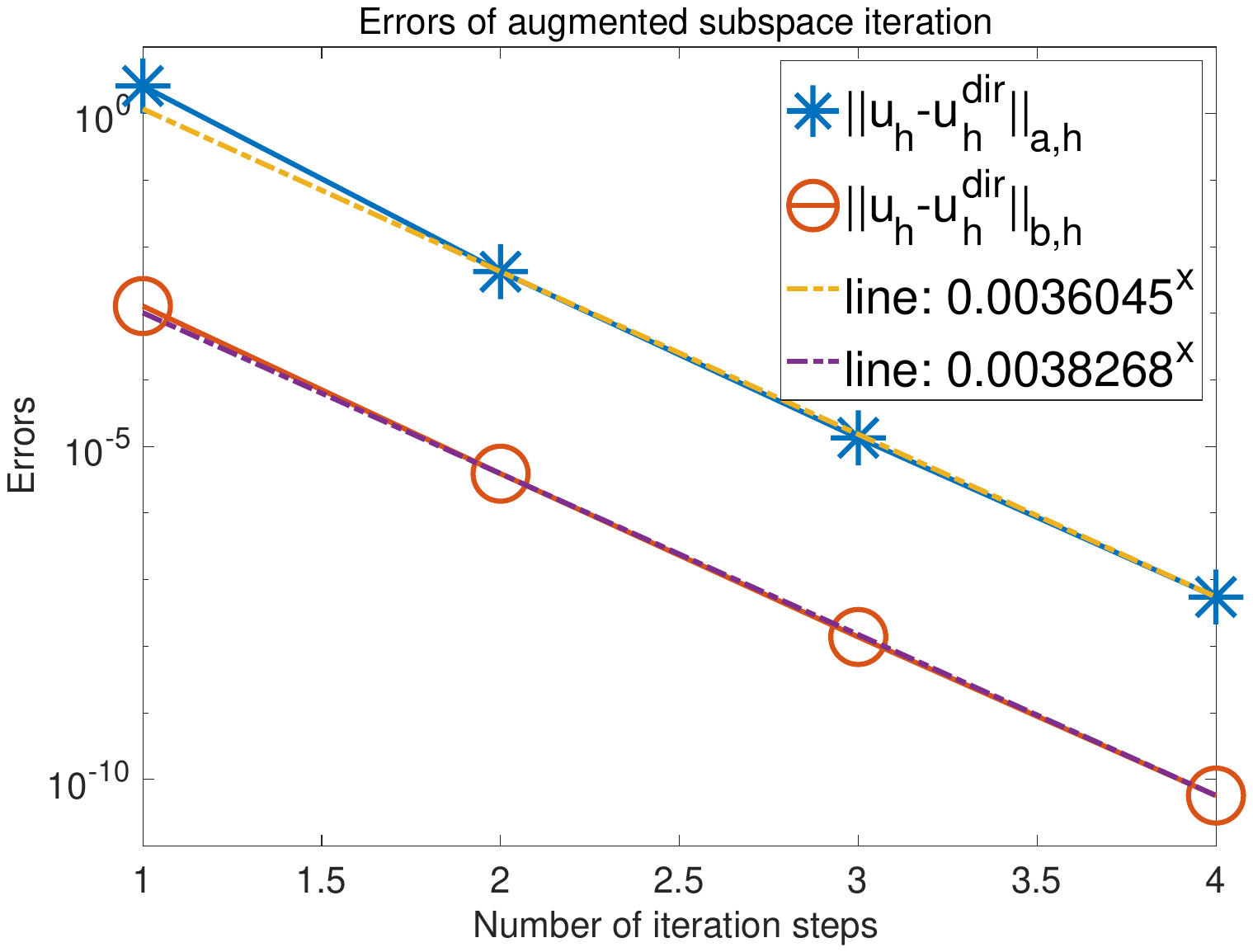}
\includegraphics[width=7cm,height=4.5cm]{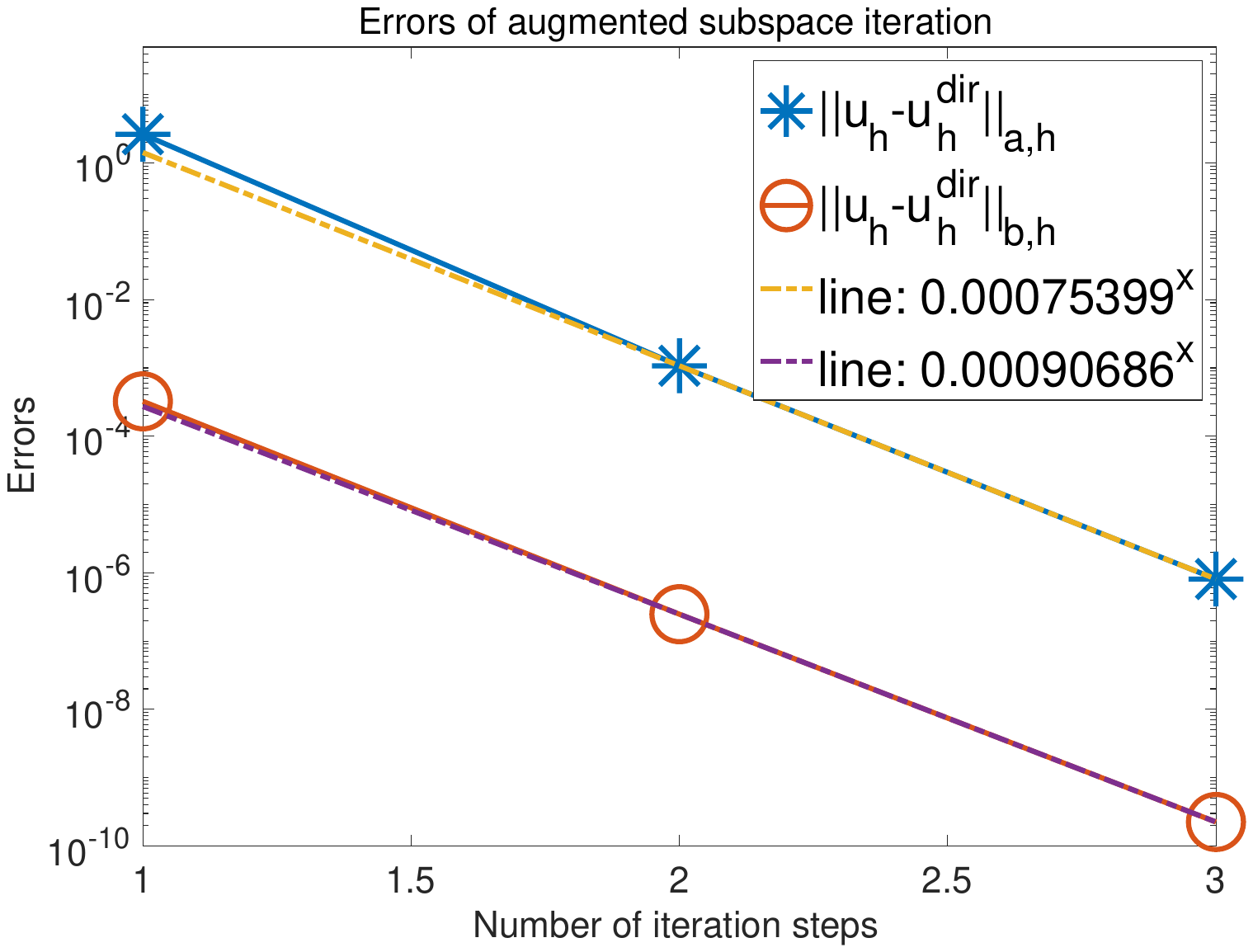}
\caption{The convergence behaviors for the first eigenfunction by Algorithm \ref{Algorithm_k}
corresponding to the $P_1/P_1$ WG finite element method and the coarse mesh size 
$H=\sqrt{2}/8$, $\sqrt{2}/16$, $\sqrt{2}/32$ and $\sqrt{2}/64$, 
respectively.}\label{Result_Coarse_Mesh_2}
\end{figure}

In a similar vein, under the boundary condition restriction, 
the initial eigenfunction approximation is also made to be rand vectors.
The convergence characteristics for the first eigenfunction using the 
augmented subspace techniques are displayed in Figure \ref{Result_Coarse_Mesh_2}, 
which corresponds to the coarse mesh sizes $H=\sqrt{2}/8$, $\sqrt{2}/16$, $\sqrt{2}/32$, 
and $\sqrt{2}/64$, respectively.
$\|\cdot\|_{a,h}$ and $\left\|\cdot\right\|_{b,h}$ have respective convergence 
rates of $0.053287$, $0.013798$, $0.0036045$, $0.00075399$ 
and $0.05535$, $0.014936$, $0.0038268$, $0.00090686$.
The findings support the results (\ref{Test_1_1})-(\ref{Test_2_0}) 
by demonstrating the second order convergence speed of the augmented 
subspace technique specified in Algorithms \ref{Algorithm_k} and \ref{Algorithm_1}.

\begin{figure}[http!]
\centering
\includegraphics[width=7cm,height=4.5cm]{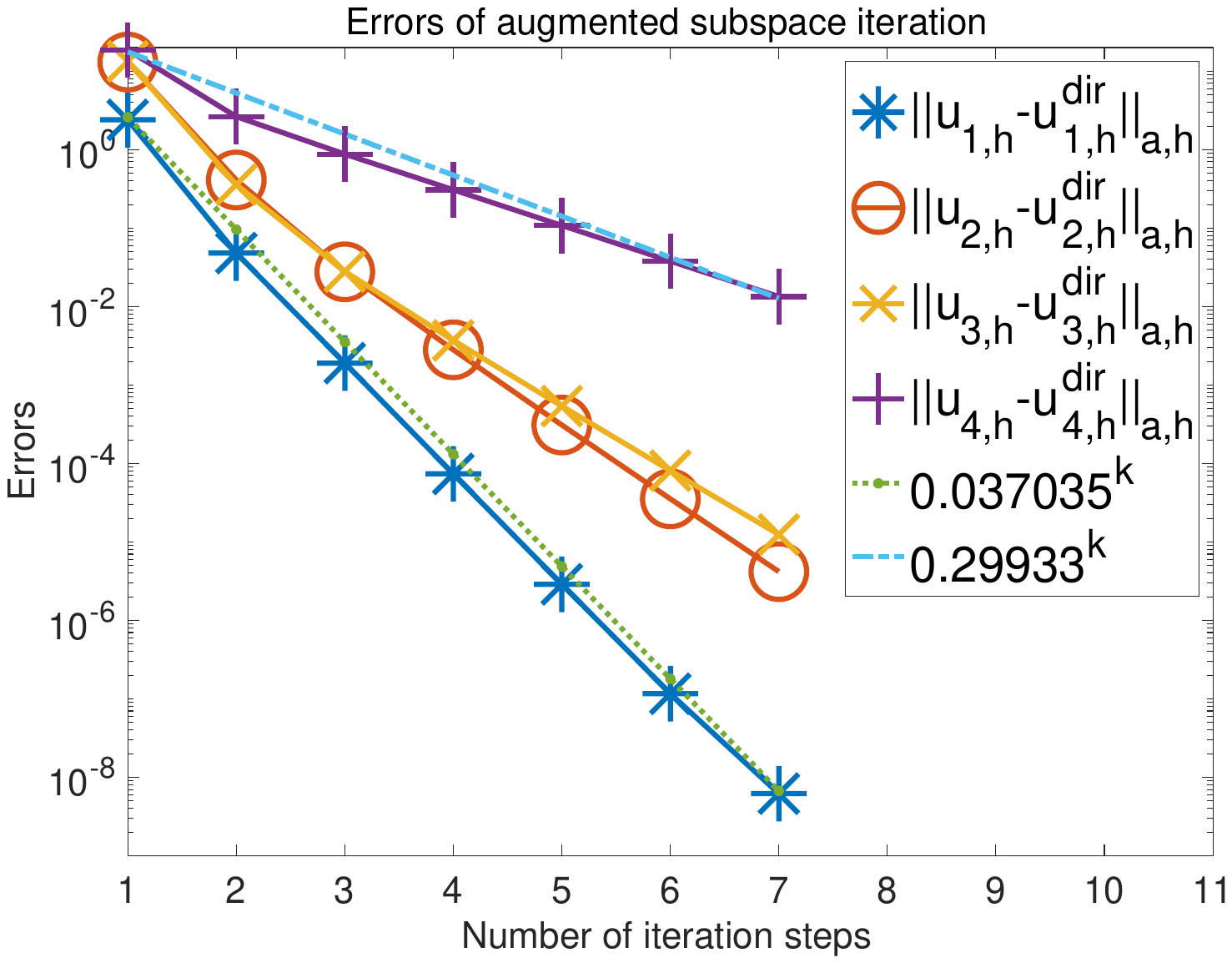}
\includegraphics[width=7cm,height=4.5cm]{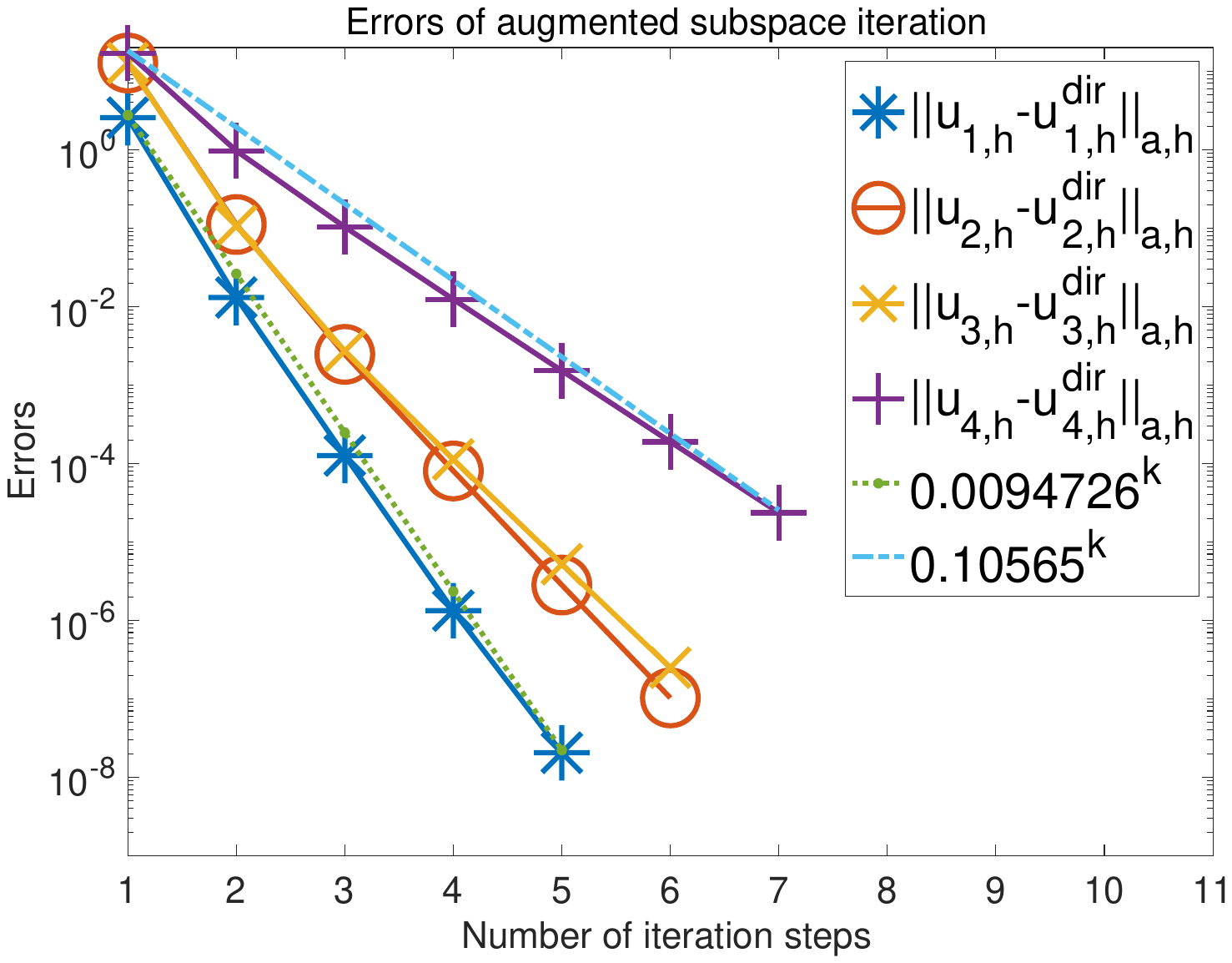}\\
\includegraphics[width=7cm,height=4.5cm]{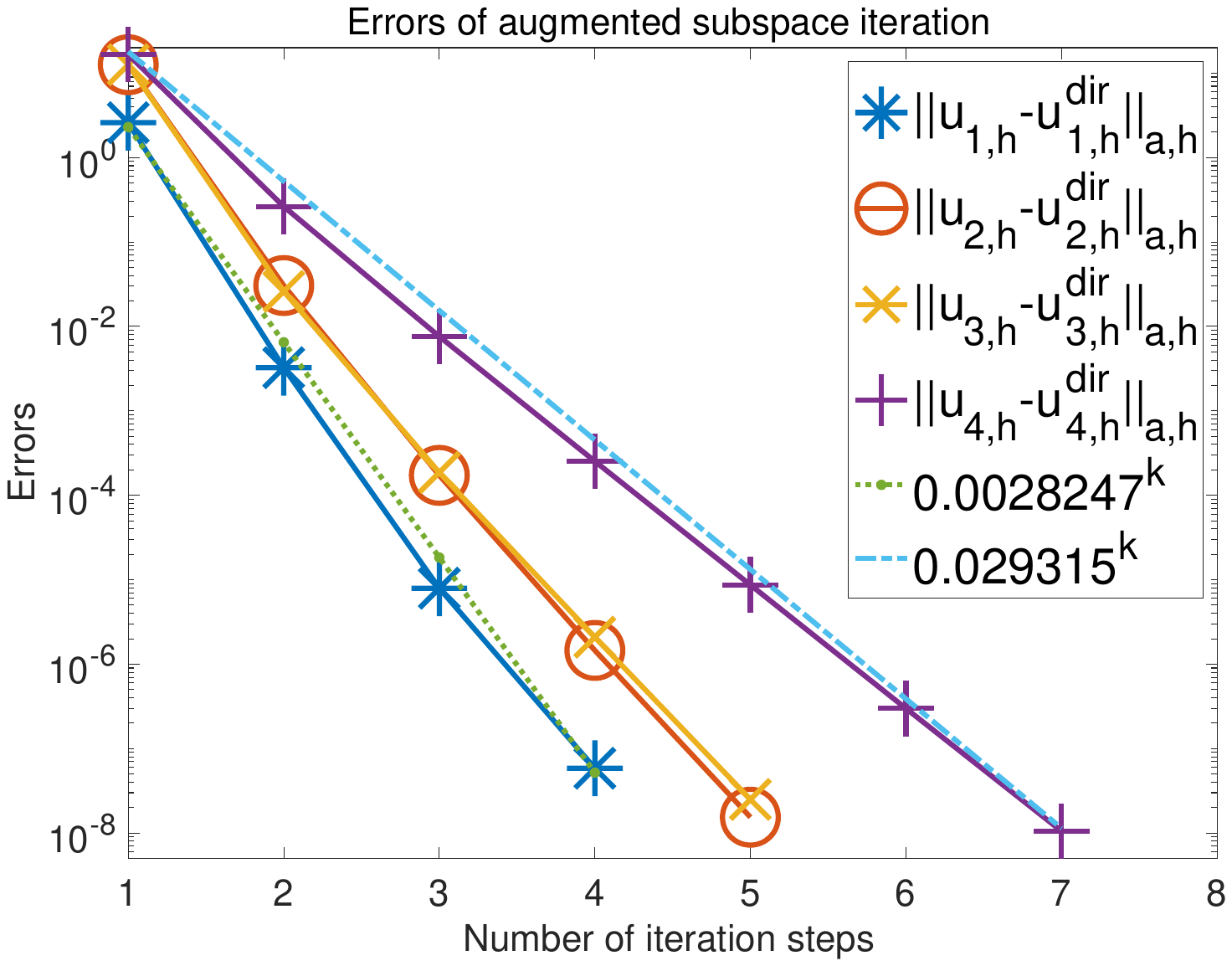}
\includegraphics[width=7cm,height=4.5cm]{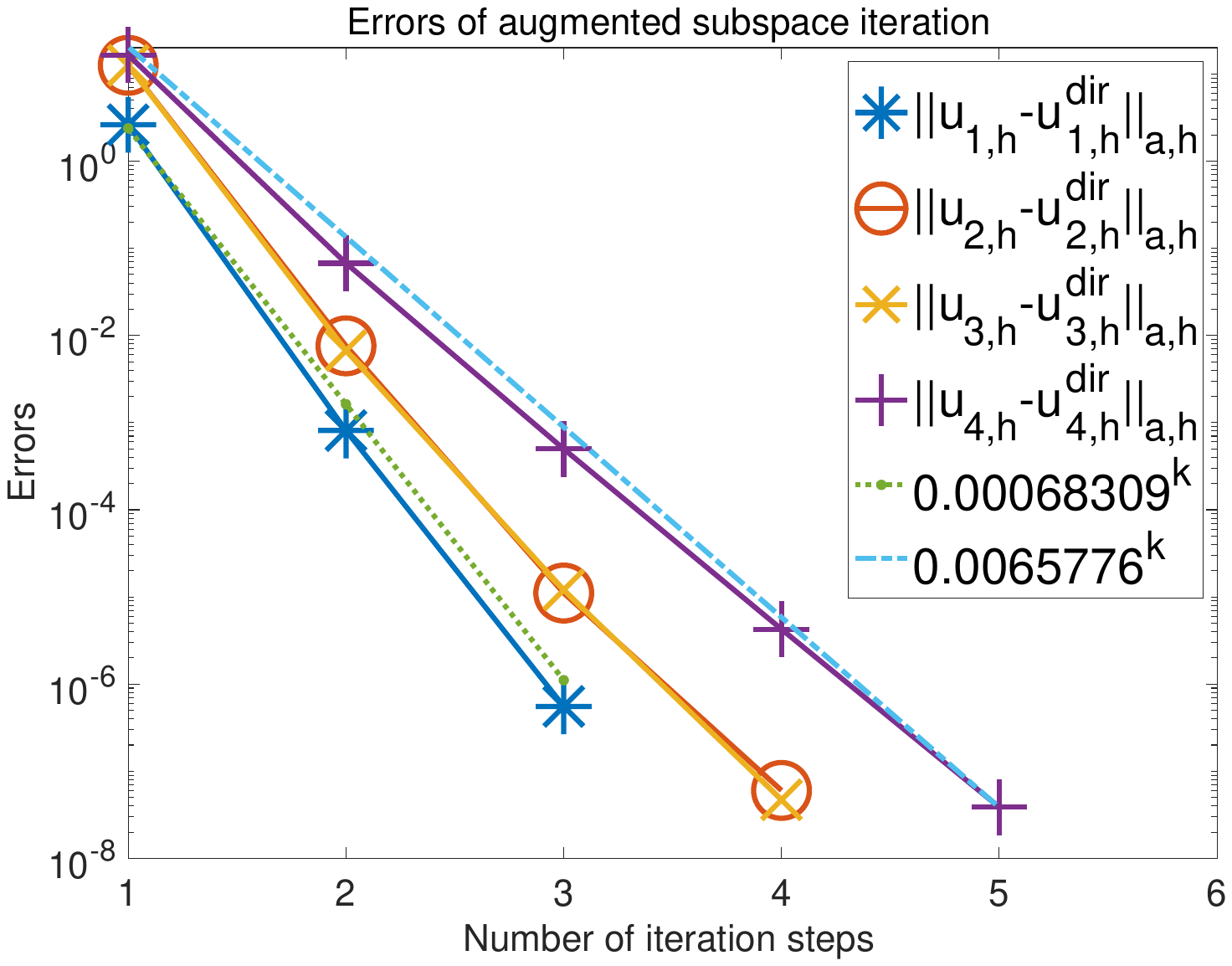}
\caption{The convergence behaviors for the smallest $4$ eigenfunctions 
by Algorithm \ref{Algorithm_k} with the $P_1/P_1$ WG finite element method 
and the coarse space being the linear finite element space on the mesh with 
size $H=\sqrt{2}/8$, $\sqrt{2}/16$, $\sqrt{2}/32$ and $\sqrt{2}/64$, 
respectively.}\label{Result_Coarse_Mesh_4_2}
\end{figure}

Next, we additionally examine Algorithm \ref{Algorithm_k}'s 
performance in terms of computing the first $4$ eigenpairs.
The corresponding convergence behaviors for the smallest $4$ eigenfunctions by 
Algorithm \ref{Algorithm_k} are presented in Figure \ref{Result_Coarse_Mesh_4_2}. 
The conforming linear finite 
element space on the mesh with sizes $H=\sqrt{2}/8$, $\sqrt{2}/16$, $\sqrt{2}/32$, 
and $\sqrt{2}/64$, respectively, constitutes the coarse space.
By employing the $4$-th eigenfunction as an example, 
we can also get the related convergence rates $0.29933$, $0.10565$, $0.029315$, and $0.0065776$, 
which reflect second convergence order of Algorithm \ref{Algorithm_k}.

\begin{figure}[http!]
\centering
\includegraphics[width=7cm,height=4.5cm]{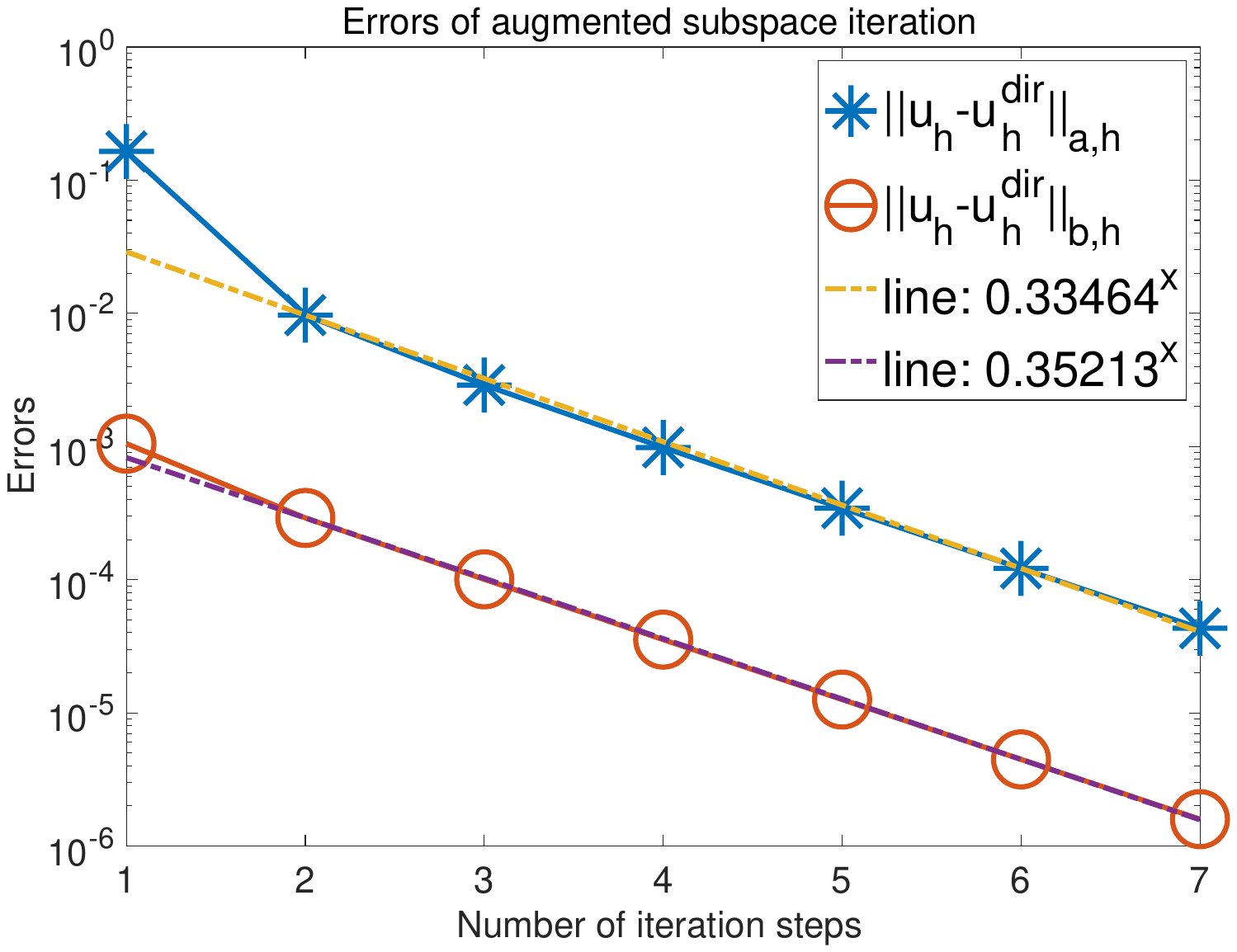}
\includegraphics[width=7cm,height=4.5cm]{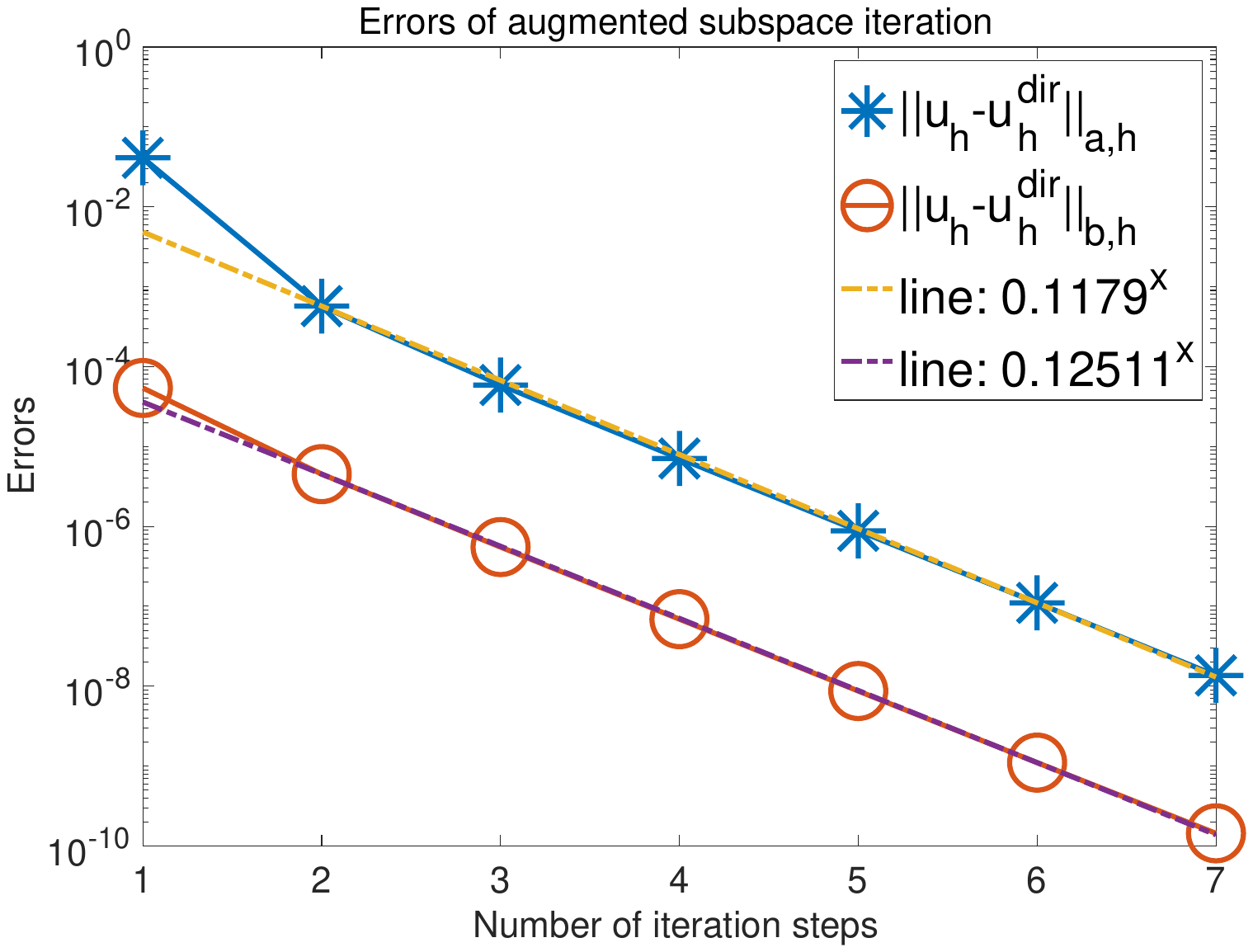}\\
\includegraphics[width=7cm,height=4.5cm]{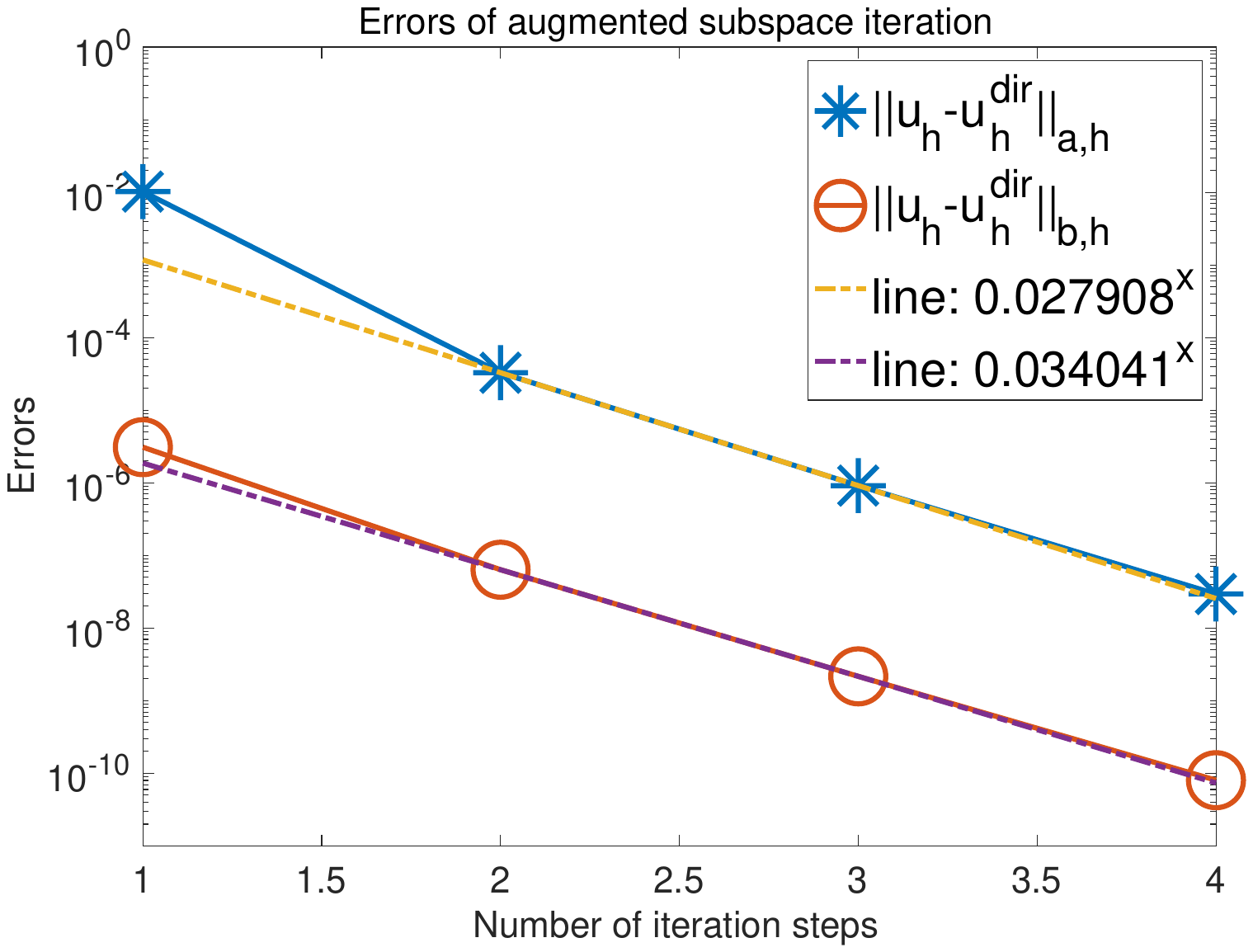}
\includegraphics[width=7cm,height=4.5cm]{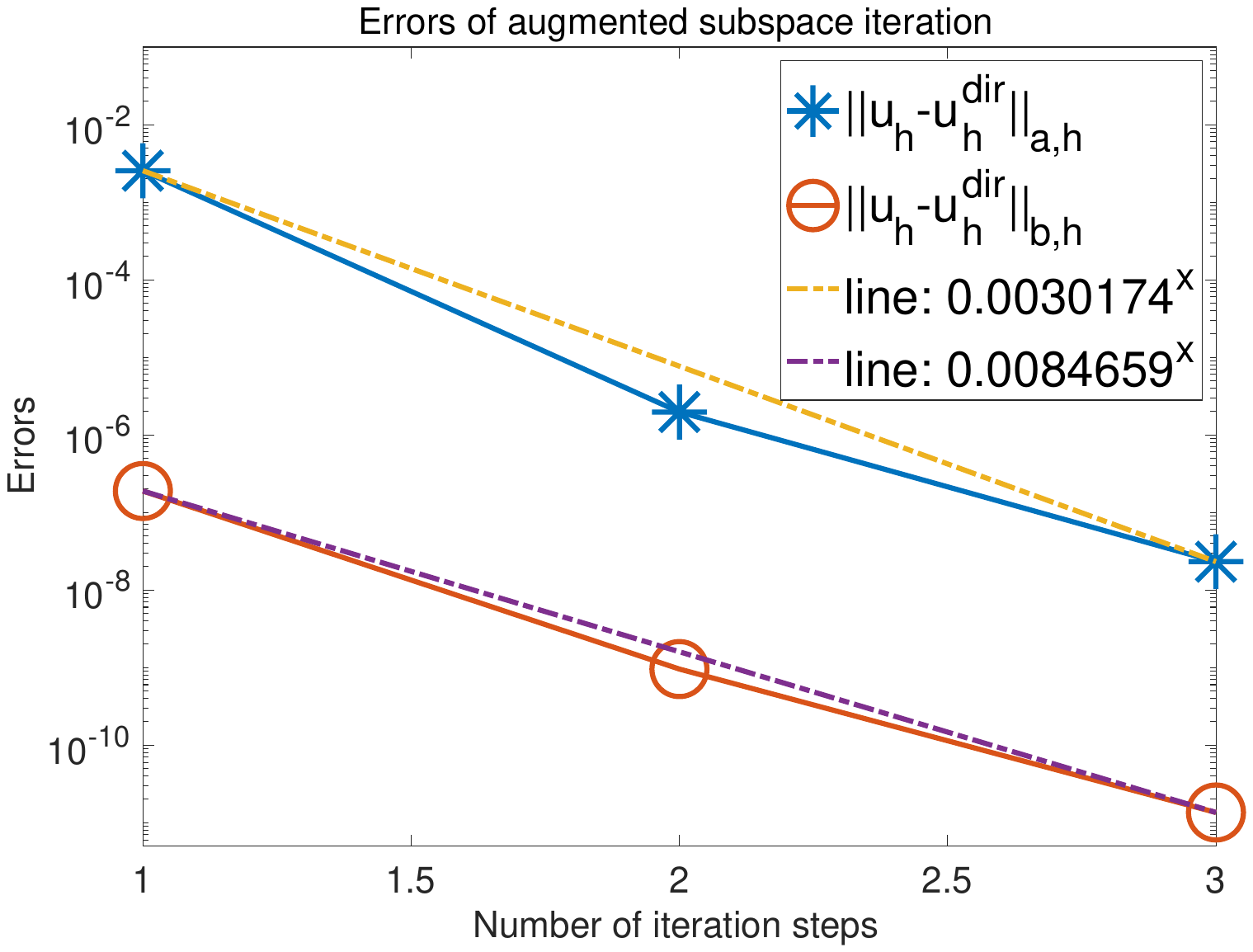}
\caption{The convergence behaviors for the only $4$-th eigenfunction by Algorithm \ref{Algorithm_1}
with the $P_1/P_1$ WG finite element method and the coarse space being the linear finite element
space on the mesh with size $H=\sqrt{2}/8$, $\sqrt{2}/16$,
$\sqrt{2}/32$ and $\sqrt{2}/64$, respectively.}\label{Result_Coarse_Mesh_4_Only_2}
\end{figure}

The final objective is evaluating the efficiency of Algorithm \ref{Algorithm_1} 
in determining the only $4$-th eigenpair.
Similarly, the coarse WG finite element space $V_H$ is used to solve the eigenvalue problem (\ref{weak_eigenvalue_problem}) to get the initial eigenfunction approximation.
The corresponding convergence behaviors for the only $4$-th eigenfunction 
by Algorithm \ref{Algorithm_1} are displayed in 
Figure \ref{Result_Coarse_Mesh_4_Only_2}. 
The conforming linear finite element space on the mesh with sizes 
$H=\sqrt{2}/8$, $\sqrt{2}/16$, $\sqrt{2}/32$, and $\sqrt{2}/64$, 
respectively, is the coarse space.
The convergence rates associated with $\|\cdot\|_{a,h}$ 
and $\left\|\cdot\right\|_{b,h}$ are $0.33464$, $0.1179$, $0.027908$, 
$0.0030174$ and $0.35213$, $0.12511$, $0.034041$, $0.0084659$, 
respectively, as depicted in Figure \ref{Result_Coarse_Mesh_4_Only_2}.
The results (\ref{Test_2_1})-(\ref{Test_2_0}) are likewise validated by these findings.

\section{Concluding remarks}
In this study, two augmented subspace strategies for addressing the eigenvalue 
problems using the WG finite element method are proposed, with the assistance 
of conforming  linear finite element space on the coarse mesh.  
We construct the associated error estimates, which demonstrate that the WG method's 
augmented subspace scheme has a second convergence order in relation to the coarse mesh size.

We can develop a sort of eigensolver for algebraic eigenvalue problems, 
which originate from the discretization of the differential eigenvalue 
problem using the WG finite element technique, based on these provided 
augmented subspace approaches. Moreover, the methods presented here provide 
a means of designing the parallel eigensolver for the WG finite element 
discretization technique, which will be the subject of our next research project.

\section*{Acknowledgements}
This work was partly supported by 
%the National Key Research and Development Program of China (No. 2019YFA0709601), 
Beijing Natural Science Foundation (No. Z200003), 
National Natural Science Foundation of China (No. 1233000214, 12301475, 12301465), 
the National Center for Mathematics and Interdisciplinary Science, CAS, 
and by the Research Foundation for Beijing University of Technology New Faculty (No. 006000514122516).

%\section*{References}
%\bibliographystyle{elsarticle-harv}
%%\bibliographystyle{elsarticle-num-names}
%%\biboptions{square,numbers,sort&compress}
%\bibliography{mybibfile}

\end{document}